\documentclass{amsart}
\usepackage{amsmath, amssymb, amsthm, amscd, graphicx}

\theoremstyle{plain}
\newtheorem{prop}{Proposition}[section]
\newtheorem{coro}[prop]{Corollary}

\newtheorem{conj}[prop]{Conjecture}
\newtheorem{lemm}[prop]{Lemma}
\newtheorem{ques}[prop]{Question}

\theoremstyle{definition}
\newtheorem*{defi}{Definition}

\newtheorem{exam}[prop]{Example}
\newtheorem{rema}[prop]{Remark}
\newtheorem*{remn}{Remark on notation}
\newtheorem*{ackn}{Acknowledgment}

\numberwithin{equation}{section}

\def\Reff#1; #2; #3; #4; #5; #6; #7\par{%
\bibitem{#1} #2, {\it #3}, #4 {\bf #5} (#6) #7}

\def\Ref#1; #2; #3; #4\par{%
\bibitem{#1} #2, {\it #3}, #4}

\def\adj{\vrule width0pt height4pt}
\def\HH#1{w_{\adj#1}}
\def\HHH#1{w^*_{\adj#1}}
\newcommand{\ev}{f}
\def\ee{e}
\def\eee{E}

\def\PP#1#2{\gc_{\adj#1, #2}}
\def\PPP#1#2{\gs_{\adj#1, #2}}

\def\gA{A}
\def\gC{C}
\def\gI{I}
\def\gJ{J}
\def\gS{S}

\def\ga{a}
\def\gc{c}

\def\gs{s}
\def\gss{\sigma}

\def\aa#1{\ga_{#1}}
\def\AA#1{\gA_{#1}}

\def\bb#1{x_{#1}}

\def\BB#1{X_{#1}}
\def\BBB#1{Y_{#1}}
\def\cc#1{\gc_{#1}}
\def\CC#1{\gC_{#1}}
\def\ss#1{\gs_{#1}}
  
\def\sss#1{\gss_{#1}}

\def\SS#1{\gS_{#1}}

\def\Ga{\boldsymbol{\ga}}
\def\GA{\boldsymbol{\gA}}

\def\Gc{\boldsymbol{\gc}}
\def\GC{\boldsymbol{\gC}}
\def\Gs{\boldsymbol{\gs}}
\def\GS{\boldsymbol{\gS}}
\def\Gss{\boldsymbol{\gss}}

\def\smT{{\scriptscriptstyle T}}

\def\AAm#1{\gA_{#1}\inv}
\def\AAp#1{\gA_{#1}}
\def\AApm#1{\gA_{#1}^{\pm1}}

\def\CCp#1{\gC_{#1}}

\def\CCtm#1{(\gC^\smT_{#1})\inv}
\def\CCtp#1{\gC^\smT_{#1}}
\def\CCtpm#1{\gC^\smT_{#1}{}^{\pm1}}
\def\IIdm#1{\gI_{#1}\inv}
\def\IIdp#1{\gI_{#1}}
\def\IIdpm#1{\gI_{#1}^{\pm1}}

\def\JJdpm#1{\gJ_{#1}^{\pm1}}

\def\SSp#1{\gS_{#1}}
\def\SSpm#1{\gS_{#1}^{\pm1}}
\def\SStm#1{(\gS^\smT_{#1})\inv}
\def\SStp#1{\gS^\smT_{#1}}
\def\SStpm#1{\gS^\smT_{#1}{}^{\pm1}}
\def\SSps#1{\gS_{#1}^{^\#}}
\def\SSpL#1{\gS_{#1}^{T_L}}

\def\Ass{\mathcal{A}}
\def\Com{\mathcal{C}}
\def\Comt{\mathcal{C}^\smT}
\def\Sem{\mathcal{S}}
\def\Semt{\mathcal{S}^\smT}
\def\IId{\mathcal{I}}
\def\JJd{\mathcal{J}}

\def\Gr(#1;#2){\langle#1\,; #2\rangle}
\def\Mon(#1;#2){\langle#1\,; #2\rangle^\smp}

\def\Bi{B_\infty}
\newcommand{\Bs}{B_\bullet}
\def\bx{x}
\def\by{y}
\def\bz{z}

\def\FG{G}
\def\Fp{F^\smp}

\def\GM#1{{\mathcal G}(#1)}

\def\GG#1{G(#1)}
\def\GGp#1{G^\smp(#1)}

\def\GGA{\GG\Ass}

\def\GGAC{\GG{\Ass,\Com}}
\def\GGAS{\GG{\Ass,\Sem}}

\def\GGIJ{\GG{\IId, \JJd, \pp}}

\def\GMA{\GM\Ass}

\def\GMAC{\GM{\Ass,\Com}}
\def\GMAS{\GM{\Ass,\Sem}}
\def\GMACT{\GM{\Ass,\Comt}}
\def\GMAST{\GM{\Ass,\Semt}}
\def\GMIJ{\GM{\IId, \JJd, \pp}}

\def\GMACL{\GM{\Ass,\Com^{\LDS}}}

\def\GMASL{\GM{\Ass,\Sem^{\LDS}}}
\def\GGASL{\GG{\Ass,\Sem^{\LDS}}}

\newcommand{\GMASF}{\GM{\Ass,\Sem^{\conj(\FG)}}}
\newcommand{\GGASF}{\GG{\Ass,\Sem^{\conj(\FG)}}}

\newcommand{\GMASG}{\GM{\Ass,\Sem^{\conj(G)}}}
\newcommand{\GGASG}{\GG{\Ass,\Sem^{\conj(G)}}}

\def\GGASB{\GG{\Ass,\Sem^{\Bs}}}

\def\GMASs{\GM{\Ass,\Sem^{^{\!\#}}}}
\def\GGASs{\GG{\Ass,\Sem^{^{\!\#}}}}

\def\Set{T}
\def\Disc{S}
\newcommand{\es}{t}
\newcommand{\esep}{s}
\newcommand{\ess}{\es'}
\newcommand{\esss}{\es''}

\def\act{\mathbin{\scriptscriptstyle\bullet}}
\def\actt{\mathbin{\underline{\scriptscriptstyle\bullet}}}
\def\ld{\LD{\hbox{-}}{\hbox{-}}}

\def\minineg{\hskip-0.3pt}
\def\LD#1#2{{#1}\minineg\lceil\minineg{#2}\rceil}

\def\LD#1#2{{}^{#1}{#2}}

\def\LD#1#2{{#1}\mathord{^{\wedge}}{#2}}

\def\LD#1#2{{#1}\minineg[\minineg{#2}\minineg]}

\def\op{\mathord{^{\wedge}}}
\def\op{\relax}
\let\OP=\circ

\let\opp=\cdot
\def\oppol{\mathord{\circ}}

\def\eq#1{\equiv_{\scriptscriptstyle#1}}

\def\eqA{\equiv_{\RA}}

\def\eqacs{\equiv_{\Racs}}
\def\eqACS{\equiv_{\RACS}}
\def\eqas{\equiv_{\Ras}}

\def\eqG{\equiv_{\scriptscriptstyle\square}}
\def\eqH{\equiv_{\includegraphics[scale=0.3]{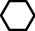}}}
\def\eqP{\equiv_{\includegraphics[scale=0.3]{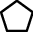}}}
\def\eqIH{\equiv_{(IH)}}

\def\RR{\boldsymbol{r}}
\def\RR{R}

\def\Ra{\RR_{\!\aa¥}}
\def\RA{\RR_{\!\AA¥}}
\def\Rac{\RR_{\!\aa¥\cc¥}}
\def\RAC{\RR_{\!\AA¥\CC¥}}
\def\Racs{\RR_{\!\aa¥\cc¥\ss¥}}
\def\RACS{\RR_{\!\AA¥\CC¥\SS¥}}
\def\Ras{\RR_{\!\aa¥\ss¥}}
\def\RAS{\RR_{\!\AA¥\SS¥}}
\def\Rass{\RR_{\!\aa¥\sss¥}}

\def\RGA{\square_{\gA}}
\def\RGC{\square_{\gC}}
\def\RGS{\square_{\gS}}
\def\RGperp{\square_\perp}
\def\RH{\boldsymbol{h}}
\def\RP{\boldsymbol{p}}
\def\RH{\includegraphics[scale=0.5]{SmallHexagon.eps}}
\def\RP{\includegraphics[scale=0.5]{SmallPentagon.eps}}

\def\shift{\partial}
\def\dd{\shift}
\def\ddp#1{\shift(#1)}
\def\ddd#1{\shift^{#1}}
\def\dddp#1#2{\shift^{#1}(#2)}
\def\DD#1{\shift_{#1}}

\newcommand{\vine}[1]{\langle#1\rangle}
\def\et{\mathord\bullet}

\def\Tree{T}
\def\Ti{\Tree_{\scriptscriptstyle\NNN}}
\def\Tii{\Ti'}
\def\To{\Tree_{\scriptscriptstyle\emptyset}}

\def\subst{\sigma}
\def\substi{\subst_1}
\def\substii{\subst_2}

\def\hh{h}
\def\II{I}
\def\nn{n}
\def\tt{t}

\def\IIi{\II'}
\def\IIii{\II''}
\def\IIiii{\II'''}
\def\nni{\nn'}
\def\nnii{\nn''}
\def\nniii{\nn'''}
\def\tti{\tt'}
\def\ttii{\tt''}
\def\ttiii{\tt'''}
\def\ttt{\tt_0}
\def\ttti{\ttt'}
\def\tttii{\ttt'}

\def\IIi{\II_1}
\def\IIii{\II_2}
\def\IIiii{\II_3}
\def\nni{\nn_1}
\def\nnii{\nn_2}
\def\nniii{\nn_3}
\def\tti{\tt_1}
\def\ttii{\tt_2}
\def\ttiii{\tt_3}
\def\ttiiii{\tt_4}
\def\ttt{\tt'}
\def\ttti{\ttt_1}
\def\tttii{\ttt_2}

\catcode`\¥=\active\def¥{\relax}
\def\0{0}
\def\1{1}

\let\a=\alpha
\newcommand{\Aut}{\mathrm{Aut}}

\let\b=\beta

\def\card#1{{\mathtt\#}#1}
\def\cl#1{\overline{\vrule width 0pt height 5pt #1}}
\def\conj{\mathrm{conj}}

\def\defect#1#2{\delta_{#1}(#2)}
\def\Dom{{\rm Dom}}

\def\e{\varepsilon}
\def\eadd{{\scriptstyle\phi}}
\newcommand{\ea}{x}
\newcommand{\eb}{y}

\def\ff{g}

\let\g=\gamma
\def\ge{\geqslant}
\def\gg{g}
\def\ggi{g'}
\def\ggii{g''}
\def\ggiii{g'''}

\def\id{{\rm id}}
\def\ie{{\it i.e.}}
\def\ince{\subseteq}
\newcommand{\inv}{^{-1}}

\def\JJ{J}

\def\KK{K}

\def\LDS{L}
\def\le{\leqslant}

\newcommand{\NNN}{\mathbf{N}}

\def\PAi{(P\!A_1)}
\def\PAii{(P\!A_2)}
\def\PAiii{(P\!A_3)}

\let\pp=\ldots

\def\Q#1{\quad#1\quad}

\def\qqquad{\quad\qquad}

\def\resp{{\it resp.\ }}

\def\smm{{\scriptscriptstyle -}}
\def\smp{{\scriptscriptstyle +}}

\def\SUP(#1, #2){\card{#1_{#2\rightarrow}}}
\def\Sy{{\mathfrak S}}
\def\Syi{\Sy_\infty}
\def\Sys{\Sy_\bullet}

\newcommand{\ww}{w}
\newcommand{\wwi}{\ww_1}
\newcommand{\wwii}{\ww_2}

\newcommand{\wwo}{\ww_0}

\def\WW#1{W(#1)}

\def\xx{x}
\def\XX{X}

\def\yy{y}

\def\zz{z}

\begin{document}

\author{Patrick DEHORNOY}
\address{Laboratoire de Math\'ematiques Nicolas
Oresme UMR 6139\\ Universit\'e de Caen,
14032~Caen, France}
\email{dehornoy@math.unicaen.fr}
\urladdr{//www.math.unicaen.fr/\!\!\!\hbox{$\sim$}dehornoy}

\title{Geometric presentations for Thompson's 
groups}



\begin{abstract}
Starting from the observation that Thompson's
groups~$F$ and~$V$ are the geometry groups
respectively of associativity, and of
associativity together with commutativity, we
deduce new presentations of these groups. These
presentations naturally lead to introducing a
new subgroup~$\Sys$ of~$V$ and a torsion free
extension~$\Bs$ of~$\Sys$. We prove that $\Sys$
and $\Bs$ are the geometry groups of
associativity together with the law $x(yz) =
y(xz)$, and of associativity together with a
twisted version of this law involving
self-distributivity, respectively.
\end{abstract}

\maketitle

Previous work showed that associating
to an algebraic law a so-called geometry group
that captures some specific geometrical
features gives useful information about that law:
the approach proved instrumental for studying
exotic laws like self-distributivity $x(yz) =
(xy)(xz)$~\cite{Dgd}
or~$x(yz)=(xy)(yz)$~\cite{Dgj}. In the case of
associativity~\cite{Dfg}, the geometry group turns
out to be Thompson's group~$F$, not a surprise as
the connection of the latter with associativity
has been known for long time~\cite{McT}.

In this paper, we develop a rather general method
for constructing geometry groups and, chiefly,
finding presentations for these groups, and we
apply this method in the case of
associati\-vity---thus finding presentations
of~$F$---and of associativity plus
commutativity---thus finding new presentations of
Thompson's group~$V$, as the latter happens to be
the involved geometry group. 

In the case of~$F$, the new presentation,
which is centered around MacLane's pentagon
relation, is more symmetric than the usual ones
and it leads to an interesting lattice structure
connected with Stasheff's associahedra; this
structure will be investigated in~\cite{Dhz}. In
the case of~$V$, on which we concentrate here,
we describe several new presentations
corresponding to various choices of the
generators. In each case, once some preliminary
combinatorial results are established, proving
that a candidate list of relations actually makes
a presentation is a straightforward application of
our general method and a very simple argument.

Perhaps the main merit of the above presentations
of~$V$ is to naturally lead to introducing two new
groups which seem interesting in themselves.
Indeed, one of these presentations explicitly
includes the Coxeter presentation of the
symmetric group~$\Syi$ (direct limit of
the~$\Sy_n$'s), thus emphasizing the existence of
a copy of~$\Syi$ inside~$V$. When we extract
those generators and relations that correspond
to~$F$ and to that copy of~$\Syi$, we obtain a
subgroup~$\Sys$ of~$V$, and, when we remove
the torsion relations $s_i^2 = 1$ in the
involved Coxeter presentation, we obtain an
extension~$\Bs$ of~$\Sys$: the connection
between~$\Bs$ and~$\Sys$ is the same as the one
between Artin's braid group~$\Bi$ and~$\Syi$. 

The algebraic and geometric properties of the
groups~$\Sys$ and, specially,~$\Bs$ are very rich.
In the current paper, we address these groups only
from the viewpoint of geometry groups, and we
prove two results: on the one hand, the
group~$\Sys$ is itself a geometry group, namely
that of associativity  together with the left
semi-commutativity law~$x(yz) = y(xz)$; on the
other hand, in some convenient sense, $\Bs$ is the
geometry group for associativity together with a
twisted version of semi-commutativity in which
$x(yz) = y(xz)$ is weakened into~$x(yz) =
x[y](xz)$, where $x, y \mapsto x[y]$ is a second
binary operation obeying a self-distributivity
condition.

The groups~$\Sys$ and~$\Bs$ to which our approach
leads turn out to be (isomorphic to) the
groups~$\widehat{V}$ and
$\widehat{BV}$ recently introduced and
investigated by M.\,Brin in~\cite{Bri, Bri1,
Bri2}. The current work can be seen as an
independent rediscovery of these groups. Let us
mention still another approach to~$\Bs$ as a group
of so-called parenthesized braids:
see~\cite{Dhe}, which contains a thorough study
of~$\Bs$. Various groups connecting Thompson's
groups and braids, some of them close to~$\Bs$,
also appear in~\cite{GrS, FuK, KaS}.

The paper is organized as follows. In
Section~1, we describe in a general context the
method that is used several times in the
paper for identifying a presentation of a
group. In Section~2, we investigate the (easy)
case of associativity and Thompson's group~$F$ as
a warm-up. In Section~3, we address the more
interesting case of associativity together with
commutativity, and obtain in this way several new
presentations of~$V$. In Section~4, we consider
the case of semi-commutativity, and of the
corresponding group~$\Sys$. Finally, Section~5
is devoted to the group~$\Bs$ and its
connection with twisted semi-commutativity and
self-distributive operations---in this section,
some algebraic results about~$\Bs$ are borrowed
from~\cite{Dhe}.  

\begin{remn}
This paper involves both Thompson's groups and
braid groups. Different notational conventions
exist. As our approach is mainly oriented toward
the group~$\Bs$, and also for the reasons listed
in~\cite{BrG}, we choose the braid conventions,
hence using actions on the right---so
$xy$ means ``$x$ then~$y$''---and numbering the
generators from~$1$. To avoid confusion, we use
a specific notation, namely~$\aa i$, for the
generators of~$F$, so that our~$\aa i$
corresponds to the standard
generator~$x_{i-1}\inv$ or $X_{i-1}\inv$
of~\cite{CFP}.
\end{remn}

\begin{ackn}
The author thanks Sean Cleary for drawing his
attention to~\cite{Bri} after a first version of
this text was written, as well as Matthew Brin
and Mark Lawson for helpful comments and
suggestions, and Charles-Antoine Lou\"et for
some corrections.
\end{ackn}

\section{A method for finding presentations}
\label{S:Crit}

Throughout the paper, $\NNN$ denotes the set
of all positive integers ($0$ excluded). 

In the sequel, we address the problem of finding
a presentation of a group several times, and we
solve it using the same argument. So it makes
sense to describe this common method first.
Although perhaps never described explicitly,  the
latter was already used in~\cite{Dgd}.

\subsection{Partial group actions}

The situation we investigate is essentially that of a
group action. However, our framework is both weaker and
stronger than the standard one. The weakening is that
the actions we consider are partial in that
every element of the group need not act on
every element; the strengthening is that our actions
satisfy a strong freeness hypothesis, namely the
existence of elements with a trivial stabilizer.

Several weak forms of group action may be thought 
of. The one convenient here is as follows. It is
essentially equivalent to the one
investigated in~\cite{KeL} (in the case of
groups)---see also~\cite{Law}---and
in~\cite{MeS} (in the case of monoids).

\begin{defi}
Let $G$ be a group, or a monoid. We define a
{\it partial (right) action} of~$G$ on a
set~$\Set$ to be a mapping~$\phi$ of~$G$ into
the partial injections of~$\Set$ into itself
such that, writing
$\es \act \gg$ for the image of~$\es$
under~$\phi(\gg)$, the following conditions are
satisfied:

$\PAi$ For every~$\es$ in~$\Set$, we have $\es \act 1 =
\es$; 

$\PAii$ For all~$\gg, \hh$ in~$G$ and $\es$ in~$\Set$,
if $\es \act \gg$ is defined, then $(\es \act \gg) \act
\hh$ is defined if and only if $\es \act\gg\hh$ is,
and, in this case, they are equal;

$\PAiii$ For each finite family $\gg_1, \pp, \gg_n$
in~$G$, there exists at least one element~$\es$
in~$\Set$ such that $\es \act \gg_1$, \pp, $\es \act
\gg_n$ are defined.
\end{defi}

Note that, in the case of a partial action, $\es
\act \gg\hh$ being defined does not guarantee
that $\es \act \gg$ is. However, the following
is easy:

\begin{lemm} \label{L:Actg}
Assume that $\phi$ is a partial action of a
group~$G$ on a set~$\Set$. Then

(i) The relations $\ess = \es \act \gg$ and  $\es = \ess
\act \gg\inv$ are equivalent;

(ii) The relation $(\exists \gg \in G)(\ess = \es \act
\gg)$ is an equivalence relation on~$\Set$.

(iii) The stabilizer of each element of~$\Set$ is
a subgroup of~$G$.
\end{lemm}

\begin{proof}
As $\gg \gg\inv$ is~$1$, $\PAii$ implies
that, if $\es \act \gg$ is defined, then $(\es \act \gg) 
\act \gg\inv$ is defined if and only if $\es \act 1$ is,
which is true by~$\PAi$. Then we
find $(\es \act \gg) \act \gg\inv = \es \act 1 =
\es$. Hence the relation of~$(ii)$ is
symmetric; $\PAi$ implies that it is reflexive, and
$\PAii$ that it is transitive. Finally,
by~$\PAii$ and~$(i)$, $\es \act \gg = \es \act
\ggi = \es$ implies that 
$\es \act (\gg\ggi)$ and $\es \act
\gg\inv$ are defined and equal~$\es$.
\end{proof}

Thus, like an ordinary (total) action, a partial
action of a group on a set~$\Set$ defines a
partition of~$\Set$ into disjoint orbits.
In the sequel we often use presentations and
expressions of the elements of a group by
words. We fix the following notation. 

\begin{defi}
Assume that $G$ is a group and that $\XX$ is a
subset of~$G$. We denote by~$\WW\XX$ the
set of all words built using letters from~$\XX
\cup \XX\inv$, \ie, all finite sequence of such
letters. For~$\ww$ in~$\WW\XX$, we usually
denote by~$\cl\ww$ the evaluation of~$\ww$
in~$G$.
\end{defi}

In the case of a partial group action, it
will be convenient to extend the action to
words:

\begin{defi}
Assume that $G$ is a group with a partial
action on~$\Set$, and $\XX$ is a subset of~$G$.
For~$\es$ in~$\Set$ and~$\ww$ in~$\WW\XX$, we
define $\es \actt \ww$ to be~$\es \act \cl\ww$
whenever $\es \act \cl\wwo$ is defined for each
prefix~$\wwo$ of~$\ww$, and to be undefined
otherwise.
\end{defi}

Note that different words representing the same
element of the group may act differently: for
instance, for every~$\xx$ in~$\XX$, the
word~$\xx \xx\inv$ and the empty word~$\e$
represent~$1$ in~$G$, but, for~$\es$ in~$\Set$,
the we always have $\es \actt \e = \es$, while
$\es \actt \xx\xx\inv = \es$ is true only if
$\es \act \xx$ is defined. However,
applying~$\PAiii$ to the (finite) family consisting
of all prefixes of~$\ww$ gives

\begin{lemm} \label{L:WordAction}
Assume that the group~$G$ has a partial action
on~$\Set$ and $\XX$ is a subset of~$G$. Then,
for each word~$\ww$ in~$\WW\XX$, there exists
at least one element~$\es$ of~$\Set$ such that
$\es \actt \ww$ is defined.
\end{lemm}

\subsection{An injectivity criterion}

Our criterion for recognizing presentations is
based on the following easy remark.

\begin{prop} \label{P:InjCrit}
Let $\pi: \widetilde{G} \to G$ be a
surjective group homomorphism. Assume that $G$
has a partial action on~$\Set$ and there exists
a map $\ev : \Set \to
\widetilde{G}$ such that
\begin{equation} \label{E:InjCrit}
\ev(\es \act \pi(x)) = \ev(\es) \opp \xx.
\end{equation}
holds for every~$\xx$ in some set that
generates~$\widetilde{G}$ and every~$\es$
in~$\Set$ such that $\es \act \pi(\xx)$ exists.
Then $\pi$ is an isomorphism.
\end{prop}

\begin{proof}
Let $\XX$ be the involved generating set
of~$\widetilde{G}$. First, for
every~$\ww$ in~$\WW\XX$, we have
\begin{equation} \label{E:Action}
\ev(\es \actt \pi(\ww)) = \ev(\es) \opp \cl\ww,
\end{equation}
where $\pi(\ww)$ is the word obtained by
replacing each letter~$\xx$ in~$\ww$
with~$\pi(\ww)$. We prove this using induction on
the length~$\ell$ of~$\ww$. For~$\ell = 1$ and
$\ww$ consisting of one letter in~$\XX$,
\eqref{E:Action} is true by hypothesis. Assume
that $\ww$ consists of one letter in~$\XX\inv$,
say $\ww = \xx\inv$. By Lemma~\ref{L:Actg}, $\ess
= \es
\act
\pi(\xx)\inv$ is equivalent to $\es = \ess \act
\pi(\xx)$. Hence, if $\es \act \pi(\xx)\inv$
exists, so does $(\es \act \pi(\xx)\inv) \act
\pi(\xx)$, and \eqref{E:InjCrit} gives
$$\ev(\es) = \ev((\es \act \pi(\xx)\inv) \act
\pi(\xx)) = \ev(\es \act \pi(\xx)\inv) \opp
\xx,$$ hence $\ev(\es \act \pi(\xx)\inv) =
\ev(\es) \opp \xx\inv$. Assume now $\ww = \wwi
\wwii$, with $\wwi, \wwii$ shorter than~$\ww$.
By definition, $\es \actt \pi(\ww)$ being
defined means that $\es \actt \pi(\wwi)$ and 
$(\es \actt \pi(\wwi)) \actt \pi(\wwii)$ are
defined, and, then, $\PAii$ and the
induction hypothesis give
$$\ev(\es \actt \pi(\ww)) 
= \ev((\es \actt \pi(\wwi)) \actt \pi(\wwii))
= \ev(\es \actt \pi(\wwi)) \opp \pi(\cl\wwii)
= (\ev(\es) \opp \pi(\cl\wwi)) \opp
\pi(\cl\wwii) = \ev(\es) \opp \pi(\cl\ww).$$

Now let $\gg$ be an element of~$\widetilde{G}$
satisfying $\pi(\gg) = 1$. Let $\ww$ be a word
in~$\WW\XX$ representing~$\gg$. By
Lemma~\ref{L:WordAction}, there
exists~$\es$ in~$\Set$ such that
$\es \actt \pi(\ww)$ is defined. Then,
\eqref{E:Action} gives 
$$\ev(\es) 
= \ev(\es \act \pi(\gg)) 
= \ev(\es \actt \pi(\ww))  
= \ev(\es) \opp \cl\ww
= \ev(\es) \opp \gg,$$
hence $\gg = 1$.
\end{proof}

\subsection{Group presentations}

If $G$ is a group and $\RR$ is a list of
relations satisfied in~$G$ by the elements of
some generating subset~$\XX$, there exists a
surjective homomorphism of the group~$\Gr(\XX;
\RR)$ onto~$G$. Proving that $(\XX; \RR)$
is a presentation of~$G$ amounts to proving that
the above morphism is injective, and this is
where Proposition~\ref{P:InjCrit} can be used. 

In the sequel, we shall consider partial actions
that satisfy strong freeness conditions. For $G$
acting on~$\Set$ and $\Disc \ince \Set$, we denote
by $\Disc \act G$ the set of all~$\esep \act \gg$
for $\esep$ in~$\Disc$ and~$\gg$ in~$G$.

\begin{defi}
Assume that $G$ has a partial action on~$\Set$. 
A subset~$\Disc$ of~$\Set$ is said to be
{\it discriminating} if, in $\PAiii$, we can
require $\es \in \Disc \act G$, no two elements
of~$\Disc$ lie in the same $G$-orbit, and each
element in~$\Disc$ has a trivial stabilizer.
\end{defi}

The first condition means that there
is an induced partial action on~$\Disc \act G$,
while the other ones guarantee that, for
each~$\es$ in~$\Disc \act G$, there exists a
unique~$\esep$ in~$\Disc$ and a unique~$\gg$
in~$G$ satisfying $\es = \esep \act \gg$. In
this case, we can select words describing the
connection between the elements of~$\Disc$ and the
elements of their orbits. When $\RR$ is a family
of relations for a group, we denote by~$\eq\RR$
the associated congruence. Our criterion takes the
following form. 

\begin{prop} \label{P:Criterion}
Let $G$ be a group with a partial action on
a set~$\Set$. Let $\XX$ be a subset
of~$G$ and $\RR$ be a collection of relations
satisfied in~$G$ by the elements of~$\XX$.
Assume that $\Disc$ is a discriminating
subset of~$\Set$ and that, for each~$\esep$
in~$\Disc$ and $\es$ in the $G$-orbit
of~$\esep$, a word~$\HH\es$ in~$\WW\XX$ is
chosen so that $\es = \esep 
\act \cl{\HH\es}$ holds.
Then a necessary and sufficient condition for
$(\XX; \RR)$ to be a presentation of~$G$ is
that, for all $\es, \ess$ in~$\Disc \act G$
and $\xx$ in~$\XX$,
\begin{equation} \label{E:Criterion}
\ess = \es \act
\xx \mbox{\qquad implies \qquad} 
\HH\ess \eq\RR \HH\es \opp \xx.
\end{equation}
\end{prop}

\begin{proof}
We begin with an auxiliary claim, namely that
$\ess = \es \act \gg$ implies
$\cl{\HH\ttt} = \cl{\HH\tt} \opp \gg$
for all $\es, \ess$ in~$\Disc \act G$ and~$\gg$
in~$G$. Indeed, assume $\ess = \es \act \gg$.
Let $\esep$ be an element of~$\Disc$ in the
orbit of~$\es$. Then $\esep$ also belongs to
the orbit of~$\ess$, and, by hypothesis, we
have $\es = \esep \act \cl{\HH\tt}$ and $\ess =
\esep \act \cl{\HH\ttt}$. On the other hand, we
also have
$\ess = \es \act \gg = (\esep \act \cl{\HH\tt})
\act \gg$, hence $\ess = \esep \act
(\cl{\HH\tt} \opp \gg)$. The hypothesis that
$\Disc$ is discriminating then implies
$\cl{\HH\ttt} = \cl{\HH\tt} \opp \gg$, as
expected. 

Let us show that \eqref{E:Criterion} is a necessary
condition. Assume $\ess = \es \act \xx$. By the
claim above, we deduce $\cl{\HH\ttt} =
\cl{\HH\tt} \opp \xx$, \ie, the words~$\HH\ttt$
and $\HH\tt \opp \xx$ represent the same
element of~$G$. If $(\XX; \RR)$ is a
presentation of~$G$, they must be
$\RR$-equivalent, and the condition is
necessary.

We turn to the converse. First,
let $\gg$ be an arbitrary element of~$G$. As
$\Disc$ is discriminating, there exists~$\es$
in~$\Disc \act G$ such that
$\es \act \gg$ exists. Let $\ess = \es \act
\gg$. By the claim above, we haves
$\cl{\HH\ttt} = \cl{\HH\tt} \opp \gg$. Now, by
hypothesis, the words~$\HH\ttt$ and $\HH\tt$
lie in~$\WW\XX$, so their classes belong to the
subgroup of~$G$ generated by~$\XX$, and so
does~$\gg$. Hence $\XX$ generates~$G$.
It remains to show that the relations of~$\RR$
make a presentation of~$G$. Let $\widetilde{G}$
be the presented group~$\Gr(\XX; \RR)$. The
set~$\XX$ generates~$G$, and the hypothesis is
that the relations of~$\RR$ are satisfied in~$G$.
Hence there exists a surjective homomorphism $\pi :
\widetilde{G}
\to G$ which is the identity on~$\XX$, and we
aim at proving that $\pi$ is injective. Now,
define $\ev :
\Disc\act G \to
\widetilde{G}$ so that $\ev(\es)$ is the
element of~$\widetilde{G}$ represented
by~$\HH\es$. If we assume
\eqref{E:Criterion}, then $\ess = \es \act \xx$
implies $\ev(\es \act \xx) = \ev(\es) \opp \xx$:
this is exactly Relation~\eqref{E:InjCrit} for
the partial action of~$G$ on~$\Disc \act G$, and
Proposition~\ref{P:InjCrit} then says that $\pi$
must be injective.
\end{proof}

The previous criterion will always be used as a
sufficient condition here. However knowing that
the condition is also necessary guarantees that
the presentations one obtains by introducing just
enough relations to witness for all equivalences
occurring in~\eqref{E:Criterion}  are in some
sense minimal. Also, adapting the criterion to the
context of monoids is easy, provided the
considered monoids admits left cancellation---but
we shall not use this version here.

\section{Thompson's group~$F$ as the
geometry group of associativity}
\label{S:PreF}

We describe now a realization of Thompson's
group~$F$ as the geometry group of the associativity
law. This is one way of formalizing the well-known
connection between~$F$ and the associativity
law, and it naturally leads to a presentation of~$F$ in
terms of a family of generators indexed by binary
addresses. Apart from more or less trivial geometric
relations, the only relations in this presentation
correspond to the well-known MacLane--Stasheff's
pentagons.

\subsection{Trees and associativity}

In the sequel, we consider finite, rooted binary 
trees---simply called trees. The number of leaves in a
tree is called its {\it size}. We denote by~$\et$ the
tree consisting of a single vertex and
by~$\tti \opp \ttii$, or simply $\tti \op
\ttii$, the tree with left subtree~$\tti$ and
right subtree~$\ttii$. Every tree has a unique
decomposition in terms of~$\et$ and the
product. 

\begin{figure} [htb]
\begin{picture}(100,17)(0, 0)
\put(0,4){\includegraphics{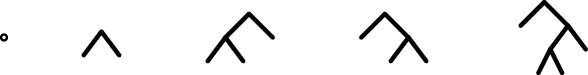}}
\put(0,0){$\et$}
\put(15,0){$\et \op \et$}
\put(36,0){$(\et \op \et) \op \et$}
\put(62,0){$\et \op (\et \op \et)$}
\put(88,0){$\et \op ((\et \op \et) \op \et)$}
\end{picture}
\caption{\smaller Typical trees with their
decomposition in terms of~$\et$}
\label{F:Tree}
\end{figure}

We also consider {\it $\LDS$-coloured}
trees, defined as trees in which the leaves wear
labels---or colours---taken from the set~$\LDS$.
We write $\et_\xx$ for~$\et$ with label~$\xx$,
and~$\Tree_\LDS$ for the set of all
$\LDS$-coloured trees. We use~$\To$ for the set of
all uncoloured trees, and see it as a subset
of~$\Ti$ by identifying an uncoloured tree with
the coloured tree where all leaves are
labelled~$1$.

The associativity law 
\begin{equation}
\tag{$\Ass$}
x(yz) = (xy)z
\end{equation}
gives rise to  an equivalence relation on
(coloured) trees: two trees~$\tt, \ttt$ are equivalent
up to associativity if we can transform~$\tt$
into~$\ttt$ by iteratively replacing one subtree of
the form~$\tti \op (\ttii
\op \ttiii)$ with the corresponding tree~$(\tti \op
\ttii) \op \ttiii$, or {\it vice versa}:
$$\begin{picture}(45,17)(0, 0)
\put(0,3){\includegraphics{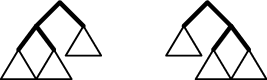}}
\put(2,0){$\tti$}
\put(8,0){$\ttii$}
\put(14,4){$\ttiii$}
\put(30,4){$\tti$}
\put(36,0){$\ttii$}
\put(42,0){$\ttiii$}
\put(21,12){$\leftrightarrow$}
\end{picture}$$
In order to describe this action precisely, we need an
indexation for the subtrees of a tree. One solution
is to describe the path from the root of the tree to
the root of the considered subtree using (for
instance) $\0$ for ``forking to the left'' and~$\1$
for ``forking to the right''.

\begin{defi}
A finite sequence of~$0$'s and~$1$'s is called an 
{\it address};  the empty address is
denoted~$\eadd$. For~$\tt$ a (coloured) tree and
$\a$ a short enough address, the
$\a$-{\it subtree} of~$t$ is the part of~$\tt$
that lies below~$\a$. The set of all~$\a$'s for
which the $\a$-subtree of~$\tt$ exists is called
the {\it skeleton} of~$\tt$.
\end{defi}

Formally, the $\a$-subtree is defined by the
following rules: the $\eadd$-subtree of~$\tt$
is~$\tt$, and, for $\a = 0\b$ (\resp
$1\b$), the
$\a$-subtree of~$\tt$ is the $\b$-subtree
of~$\tti$ (\resp $\ttii$) when $\tt$ is~$\tti
\op \ttii$, and it is undefined in other cases.
For instance, for $\tt = \et \op (( \et
\op \et) \op \et)$ (the rightmost example in
Figure~\ref{F:Tree}), the $\1\0$-subtree
of~$\tt$ is $\et \op \et$, while its
$\0\1$- and $\1\1\1$-subtrees are undefined. The
skeleton of~$\tt$ consists of~$\eadd$, $\0$,
$\1$, $\1\0$, $\1\0\0$, $\1\0\1$, $\1\1$.
 
Applying associativity to a tree~$\tt$ consists in
choosing an address~$\a$ in the skeleton of~$t$ and
either replacing the $\a$-subtree of~$\tt$, supposed to
have the form $\tti \op (\ttii \op \ttiii)$, by
the corresponding $(\tti \op \ttii) \op \ttiii$, or
performing the inverse substitution. We can see
this as applying an operator.

\begin{defi}
$(i)$ We denote by~$\AAp¥$ the partial
operator on~$\Ti$ that maps every tree of the
form $\tti \op (\ttii \op \ttiii)$ to the
corresponding tree
$(\tti \op \ttii) \op \ttiii$. 

$(ii)$ For $\a$ an address and $f$  a partial mapping
on trees, we define the $\a$-{\it shift}
of~$f$, denoted~$\DD\a f$, to be the partial
mapping consisting in applying~$f$ to the
$\a$-subtree of its argument (when the latter
exists). We write~$\dd¥$ for~$\DD\1¥$.

$(iii)$ For~$\a$ an address, we put $\AAp\a =
\DD\a\AAp¥$. We define~$\GMA$ to be the
monoid generated by all~$\AAp\a$'s and their
inverses using reversed composition.
\end{defi}

\begin{figure} [htb]
\begin{picture}(83,30)(0, 0)
\put(0,0){\includegraphics{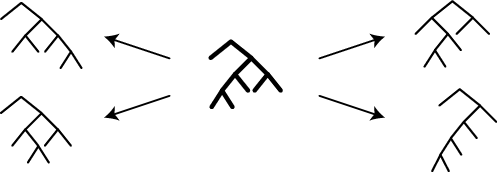}}
\put(23,7){$\AAm{\1\0}$}
\put(23,23){$\AAm\1$}
\put(57,7){$\AAp\1$}
\put(57,23){$\AAp¥$}
\put(40,23){$\eadd$}
\put(43,20){$\scriptstyle\1$}
\put(36,16.5){$\scriptstyle\1\0$}
\end{picture}
\caption{\smaller Two  operators~$\AAp\a$ and
two operators~$\AAm\a$ apply to the tree $\et
\op (((\et \op \et) \op \et) \op (\et \op \et))$}
\label{F:ExAA}
\end{figure}

\begin{exam} (Figure~\ref{F:ExAA})
Let $\tt =\et \op (((\et \op \et) \op \et) \op (\et \op
\et))$. Then $\tt$ lies in the domain of~$\AAp¥$,
as the $\eadd$-subtree of~$\tt$, \ie, $\tt$
itself, is $\tti \op (\ttii \op
\ttiii)$, with $\tti = \et$, $\ttii = (\et \op
\et) \op \et$, and $\ttiii  \et \op \et$. Then
the image of~$\tt$ under~$\AAp¥$ is $(\tti \op
\ttii)
\op \ttiii$, \ie,  $(\et \op ((\et \op \et) \op \et)) \op
(\et \op \et)$. Similarly, $\tt$ lies in the domain
of~$\AAp\1$, and in the images of~$\AAp\1$ and
of~$\AAp{\1\0}$, hence in the domains of
$\AAm\1$ and~$\AAm{\1\0}$. These are
the only operators~$\AApm\a$ applying to~$t$.
\end{exam}

We thus have a partial action of the
monoid~$\GMA$ on trees in the sense of
Section~\ref{S:Crit}; for~$f$ in~$\GMA$, we write
$\tt \act f$ for the image of~$\tt$ under~$f$, when
it exists. We use reversed composition in~$\GMA$ so
as to make our multiplication compatible with an
action on the right.

By construction, two trees~$\tt, \ttt$ are
equivalent up to associativity if and only if
some element of~$\GMA$ maps~$\tt$
to~$\ttt$. Thus the orbits for the
partial action of the monoid~$\GMA$ are the
equivalence classes with respect to
associativity. In particular, there is exactly
one orbit for each size inside~$\To$, and the
cardinal of the orbit of size~$n$ trees is
the $n$th Catalan number.

\subsection{Making $\GMA$ into a group}

Except the identity mapping, the elements
of~$\GMA$  are partial mappings, and the
monoid~$\GMA$ is not a group, but only an
inverse monoid, \ie, a monoid in which, for each
element~$g$, there exists~$g\inv$
satisfying $g g\inv g = g$ and $g\inv g g\inv =
g\inv$. For instance, the product~$\AAp¥ \AAm¥$ is
the identity of its domain, but the latter does not
contain~$\et$.

Every inverse monoid admits a maximal
quotient-group, called its universal
group~\cite{How, Pat}. In the general case, the
universal group may be much smaller than the
original monoid, typically when the latter
consists of partial mappings whose domains may
be disjoint. In the current case, no wild
collapsing occurs, and the induced action of the
universal group keeps the freeness properties
of the initial monoid action. As the same
construction will be used several times, we 
describe it in a general framework.

\begin{defi}
Two partial mappings~$\gg, \gg'$ are said
{\it near-equal}, denoted $\gg \approx \gg'$, if
there is at least one element~$\es$ such that
both $\es \act \gg$ and $\es \act \gg'$ are
defined, and $\es \act \gg = \es \act \gg'$
holds for every such~$\es$.
\end{defi}

\begin{lemm} \label{L:Mono}
Assume that $\mathcal{G}$ is a monoid consisting of
partial self-injections of a set~$\Set$ that is
closed under inverse, and there exists a
subset~$\Disc$ of~$\Set$ such that, for all
$\gg_1, \pp, \gg_n, \gg, \gg'$ in~$\mathcal{G}$,
\begin{gather}
\label{E:Con1}
\text{$\Dom(\gg_1) \cap \pp \cap \Dom(\gg_n)
\cap\Disc \act \mathcal{G}$ is nonempty,}\\
\label{E:Con2}
\text{$\gg \approx \gg'$ is true
whenever $\es \act \gg = \es \act \gg'$ holds
for some~$\es$ in~$\Disc \act \mathcal{G}$.}
\end{gather}
Then near-equality is a congruence on~$\mathcal{G}$,
the quotient-monoid is a group, the 
mappings of~$\mathcal{G}$ induce a partial action
of this group on~$\Set$, and the
set~$\Disc$ is discriminating for this partial
action.
\end{lemm}

\begin{proof}
Assume $\ggi \approx \ggii \approx \ggiii$.
By~\eqref{E:Con1}, there exists~$\es$
in~$\Disc \act \mathcal{G}$ such that $\es \act
\ggi, \es \act \ggii$, and $\es \act \ggiii$ are
defined. Then one necessarily has $\es \act \ggi
= \es \act \ggiii$, hence $\ggi \approx \ggiii$
by~\eqref{E:Con2}, and $\approx$ is an
equivalence relation. Next, $\ggi \approx \ggii$
implies $\gg \ggi \approx \gg\ggii$ and $\ggi
\gg \approx \ggii \gg$ for every~$\gg$, because \eqref{E:Con1} guarantees
that there exists~$\es$ in~$\Disc \act \mathcal{G}$
 for which $\es \act \gg$, $\es \act \gg \ggi$,
$\es \act \gg\ggii$, $\es \act \ggi$, $\es
\act \ggi \gg$, $\es \act \ggii$, and $\es \act
\ggii\gg$ are defined. So
$\approx$ is a congruence on~$\mathcal{G}$, and the
quotient-monoid~$\mathcal{G}/\!\approx$, henceforth
denoted~$G$, is well-defined. For each~$\gg$
in~$\mathcal{G}$, we have $\gg \gg\inv \approx \id$
because $\Dom(\gg)$ is nonempty, so $G$
is a group.

For~$\gg$ in~$\mathcal{G}$, let us denote by~$\cl\gg$
the class of~$\gg$ in~$G$. For~$\es$
in~$\Set$, and~$\ea$ in~$G$,  we define
$\es \act \ea$ to be~$\ess$ if $\es \act \gg =
\ess$ holds for some element~$\gg$
of~$\mathcal{G}$ satisfying $\cl\gg =
\ea$, if such an element exists. Then $\es \act
\ea$ is well-defined by definition of~$\approx$,
and we claim that one obtains in this way a partial
action of~$G$ on~$\Set$. Indeed,
Condition~$\PAi$ is trivial. As for~$\PAii$,
assume that $\es \act \ea$ and
$(\es \act \ea) \act \eb$ are defined. This means
that there exist~$\gg, \hh$ with $\ea = \cl\gg$
and $\eb = \cl\hh$ such that $\es \act \gg$ and
$(\es \act \gg) \act \hh$ are defined. But, then, 
$\es \act \gg\hh$ is defined, and, by
construction, we have $\cl{\gg\hh} = \cl\gg \,
\cl\hh$. Conversely, assume that $\es \act
\ea$ and $\es \act \ea\eb$ are defined, say $\es
\act \ea = \ess$ and $\es \act \ea\eb = \esss$. This
means that there exist~$\gg, \ggi$
in~$\mathcal{G}$ satisfying $\es \act \gg = \ess$,
$\es \act \ggi = \esss$, with $\cl\gg = \ea$ and
$\cl\ggi = \ea\eb$. Let
$\hh = \gg\inv \ggi$. Then $\hh$ belongs
to~$\mathcal{G}$, we have $\cl\hh = \ea\inv \ea \eb =
\eb$, and $\ess \act \hh = \esss$. This shows that
$(\es \act \ea) \act \eb$ is defined, and equal
to~$\esss$. So Condition~$\PAii$ is satisfied.
Then \eqref{E:Con1} implies~$\PAiii$  directly,
and we obtain a partial action of~$G$
on~$\Set$. Finally, the subset~$\Disc$ is
discriminating by~\eqref{E:Con2}.
\end{proof}

In order to apply the previous construction to
the monoid~$\GMA$ and its action on trees, we
describe the domain and the image of a
generic element of~$\GMA$ explicitly.

\begin{defi}
$(i)$ A mapping of~$\NNN$ to~$\Ti$ is
called a {\it substitution}. If $\tt$ is a tree
in~$\Ti$ and $\subst$ is a substitution, we
denote by~$\tt^\subst$ the tree obtained by
replacing each leaf~$\et_\xx$ in~$\tt$ with the
tree~$\subst(\xx)$.

$(ii)$ A coloured tree is said to be  {\it
injective} if its labels are pairwise
distinct.

$(iii)$ For~$\ff$ a partial mapping  of~$\Ti$
into itself, we say that a pair of trees $(\tt,
\tt')$ in~$\Ti$ is a {\it seed} for~$\ff$ if, as
a set of pairs, $\ff$ is the set of all
$(\tt^\subst, \tt'{}^\subst)$ with
$\subst$ a substitution.
\end{defi}

The pair $(\et_1 \op (\et_2 \op
\et_3), (\et_1 \op \et_2) \op \et_3)$ is a seed
for~$\AAp¥$: this is just saying that $\AAp¥$
consists of all pairs of the form $(\tt_1 \op (\tt_2
\op \tt_3), (\tt_1 \op \tt_2) \op \tt_3)$.
Then we have the following general result:

\begin{lemm}\label{L:DomF}
Each element of~$\GMA$ admits a seed consisting
of injective trees.
\end{lemm}

\begin{proof}
Let~$\ff$ be an element of~$\GMA$. We use
induction on the (minimal) length of a
decomposition of~$\ff$ in terms of the
operators~$\AAp\a$ and~$\AAm\a$. The pair~$(\et_1, \et_1)$  is a
seed for~$\gg = \id$, the pair $(\et_1 \op
(\et_2 \op \et_3), (\et_1 \op \et_2) \op \et_3)$ is
a seed for~$\gg = \AAp¥$, and it is easy to define
similarly a seed for~$\gg = \AApm\a$. Otherwise,
write $\ff = \ff_1 \ff_2$. By induction hypothesis,
$\ff_1$ and $\ff_2$ admit seeds, say $(\tt_1,
\tt'_1)$ and $(\tt_2, \tt'_2)$. If $\tt'_1$
happens to coincide with~$\tt_2$, then $(\tt_1,
\tt'_2)$ is a seed for~$\ff$. In the general case,
because $\tt'_1$ and~$\tt_2$ are injective, there
exist minimal substitutions~$\substi$
and~$\substii$ such that $\tt'_1{}^{\substi}$ and
$\tt_2^{\substii}$ coincide, and, then, the pair
$(\tt_1^{\substi}, \tt'_2{}^{\substii})$ is
a seed for~$\ff$.
\end{proof}

(Moreover, the seed is unique if the labels are
requested to make an initial segment of~$\NNN$.)

\begin{coro}
The monoid~$\GMA$ satisfies 
Conditions~\eqref{E:Con1} and~\eqref{E:Con2} of
Lem\-ma~\ref{L:Mono} with $\Set = \Ti$ and
$\Disc$ any subset of~$\Ti$ containing trees of
arbitrary large size.
\end{coro}

\begin{proof}
Let $\Disc$ be a subset of~$\Ti$ containing
trees of arbitrary large size, and
let~$\tt$ be an arbitrary tree. Then there
exists~$\esep$ in~$\Disc$ whose size is at
least that of~$\tt$. Using associativity, we
can transform~$\esep$ into a tree whose
skeleton includes that of~$\tt$, \ie, there
exists~$\ff$ in~$\GMA$ such that $\esep \act
\ff$ is defined and its skeleton includes that
of~$\tt$.

Let $\ff_1, \pp, \ff_n$ be elements
of~$\GMA$, and $(\tt_1, \tt'_1)$, \pp, $(\tt_n,
\tt'_n)$ be seeds for these elements. By the
above argument, there exists a tree~$\tt$
in~$\Disc \act \GMA$ whose skeleton includes
the skeletons of~$\tt_1, \pp, \tt_n$, hence there
exist substitutions~$\substi, \pp, \subst_n$
such that $\tt = \tt_i^{\subst_i}$ holds for
each~$i$, which implies that $\tt \act \ff_i$
is defined for each~$i$. So
Condition~\eqref{E:Con1} is satisfied.

Assume that $\ff_1, \ff_2$ belong to~$\GMA$, and
$\tt \act \ff_1 = \tt \act \ff_2$ holds for some
tree~$\tt$ in~$\Ti$. Let $(\tt_1, \tt'_1)$,
$(\tt_2, \tt'_2)$  be seeds for~$\ff_1$
and~$\ff_2$ respectively. As above, there exist
substitutions~$\substi, \substii$ such that the
trees $\tt_1^{\substi}$ and~$\tt_2^{\substii}$
coincide, they are injective, and their common
skeleton is the union of the skeletons of~$\tt_1$
and~$\tt_2$. The hypothesis that
$\tt \act \ff_1$ and $\tt \act \ff_2$ are defined
implies that the skeleton of~$\tt$ includes those
of~$\tt_1$ and~$\tt_2$, hence their union. Hences
there exists a substitution~$\subst$ satisfying
$\tt = (\tt_1^{\substi})^\subst  =
(\tt_2^{\substii})^\subst$. The hypothesis
that $\tt \act \ff_1$ and $\tt \act \ff_2$ are
equal then gives
$$ (\tt'_1{}^{\substi})^\subst 
= \tt \act \ff_1
= \tt \act \ff_2
= (\tt'_2{}^{\substii})^\subst.$$
This implies that the skeletons
of~$\tt'_1{}^{\substi}$ and $\tt'_2{}^{\substii}$
coincide. Moreover, the hypothesis
$\tt_1^{\substi} = \tt_2^{\substii}$ implies that
the sequence of labels in
$\tt_1^{\substi}$ and $\tt_2^{\substii}$
coincide. As associativity does not change the
order of the labels, the trees
$\tt'_1{}^{\substi}$ and $\tt'_2{}^{\substii}$
must coincide. This means that $\ff_1$ and
$\ff_2$ agree on every tree whose skeleton
includes that of $\tt_1^{\substi}$, \ie, on every
tree in the intersection of the domains of~$\ff_1$
and~$\ff_2$. In other words, $\ff_1 \approx
\ff_2$ holds, and Condition~\eqref{E:Con2} is
satisfied.
\end{proof}

By applying Lemma~\ref{L:Mono}, we obtain:

\begin{prop} \label{P:GroF}
Near-equality is a congruence on the
monoid~$\GMA$, and the quotient-monoid is a
group. The operators~$\AApm\a$ induce a partial
action of this group on~$\Ti$, and every
subset of~$\Ti$ containing trees of
unbounded sizes is discriminating for this partial
action. 
\end{prop}

\begin{defi} 
The {\it geometry group of associativity},
denoted~$\GGA$, is defined to be the
quotient-monoid~$\GMA / \! \approx$. 
\end{defi}

In the sequel, we still use~$\AA\a$ for the class
of~$\AAp\a$ in~$\GGA$. For~$\tt$ a tree and
$\gg$ an element of~$\GGA$, we denote by~$\tt
\act \gg$ the result of letting~$\gg$
act on~$\tt$. The elements of~$\GGA$ are
expressed by words on~$\GA$, and we also
use~$\act$ for the word action, \ie, we 
do not distinguish between~$\act$ and~$\actt$.
But we recall that $\tt
\act \ww$ exists only if $\tt \act \cl\wwo$
exists for each prefix~$\wwo$ of~$\ww$: for
instance, $(\et\op\et) \act \AA¥ \AA¥\inv$ is
not defined, since $(\et\op\et) \act \AA¥$ is
not.

It is straightforward to connect the
geometry group~$\GGA$ with Thompson's group~$F$:

\begin{prop}
The group~$\GGA$ is (isomorphic to) Thompson's
group~$F$, \ie, $F$ is the geometry group of
associativity.
\end{prop}

\begin{figure} [htb]
\begin{picture}(88, 42)(0, 0)
\put(0,0){\includegraphics{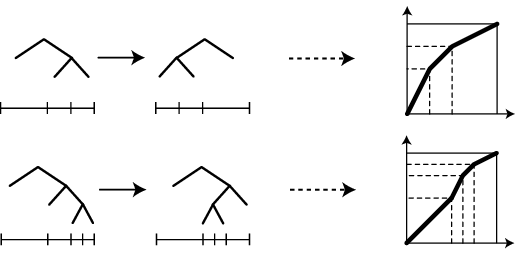}}
\put(19,34){$\AAp¥$}
\put(19,12){$\AAp\1$}
\end{picture}
\caption{\smaller From~$\GMA$ to~$F$: the
action of~$\AA¥$ and~$\AA\1$}
\label{F:FromGAToF}
\end{figure}

\begin{proof}
(Figure~\ref{F:FromGAToF}) We start with the
definition of~$F$ as a group of orientation
preserving piecewise linear homeomorphisms of the
unit interval, {\it cf.}~\cite{CFP}. Let
$\ff$ be an arbitrary element in~$\GMA$. We
map~$\ff$ to~$F$ as follows: let
$(\tt, \ttt)$ be a seed for~$\ff$; we
associate with~$\tt$ a dyadic decomposition~$0 =
r_0 < r_1 < \pp < r_n = 1$ of~$[0, 1]$, and,
similarly, let $0 = r'_0 < r'_1 < \pp < r'_n = 1$
be the dyadic decomposition associated
with~$\ttt$; then we map~$\ff$ to the unique
piecewise linear homeomorphism that maps~$r_i$
to~$r'_i$ and interpolates the values. We
obtain in this way a morphism~$\pi: \GMA \to F$.
The homeomorphisms associated with~$(\tt, \ttt)$
and~$(\tt^\subst, \ttt{}^\subst)$ coincide, and
this implies that $\pi$ factors
through~$\approx$. The injectivity of the
resulting morphism follows from the fact that
each element of~$F$ is determined  by its values
on a finite dyadic partition; its surjectivity
follows from the fact that the images of~$\AAp¥$
and~$\AAp\1$ generate~$F$.
\end{proof}

From now on, we identify~$F$ with~$\GGA$.

\subsection{Guessing relations in~$\GGA$}

Considering the group~$F$ as the geometry group of
associativity naturally leads to a
presentation of~$F$ in terms of the
generators~$\AA\a$. We proceed in two steps: 
first, we use the geometric definition
of the operators~$\AAp\a$ to guess a list of
relations; then, we prove
that these relations make a presentation using
the method of Section~\ref{S:Crit}.

\begin{figure} [htb]
\begin{picture}(108,45)(0, 0)
\put(0,0){\includegraphics{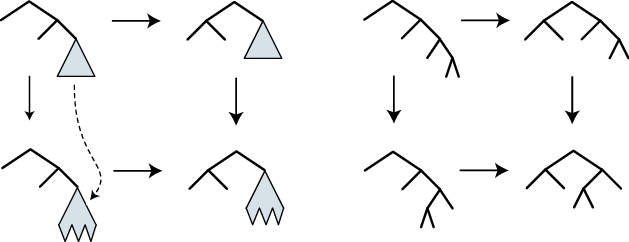}}
\put(6,42){$\tt$}
\put(41,42){$\ttt$}
\put(16,20){$f$}
\put(-1.5,24){$\ddd2 f$}
\put(42,24){$\dd f$}
\put(22,39){$\AAp¥$}
\put(22,14){$\AAp¥$}
\put(81,39){$\AAp¥$}
\put(81,14){$\AAp¥$}
\put(60.5,25){$\AAp{\1\1}$}
\put(98.5,25){$\AAp\1$}
\end{picture}
\caption{\smaller Geometric relations in~$\GMA$:
the general scheme and one example}
\label{F:Heir}
\end{figure}

Let us look for relations between the
operators~$\AAp\a$. We shall describe two types of
relations: the {\it geometric} relations,
and the {\it pentagon relations}. Geometric
relations arise when we consider inheritance
phenomena. Assume
$\ttt = \tt \act \AA¥$, \ie, assume that the
operator~$\AAp¥$ maps~$t$ to~$\ttt$. Then, by
definition, the $\1$-subtree of~$\ttt$ is a copy of
the
$\1\1$-subtree of~$\tt$. It follows that
performing any transformation in the latter
subtree and then applying~$\AAp¥$ has the same
result as applying~$\AAp¥$ first and performing
the considered transformation in the
$\1$-subtree of~$\ttt$. Therefore, the equality
\begin{equation}\label{E:Heir}
\dd^2 f \cdot \AAp¥ = \AAp¥ \cdot \dd f
\end{equation}
holds for every (partial) mapping~$f$
on trees (Figure~\ref{F:Heir}). In particular, for
$f = \AAp\a$, we
obtain
\begin{equation}
\AAp{\1\1\a} \cdot \AAp¥ = \AAp¥ \cdot
\AAp{\1\a},
\end{equation}
a typical example of what we shall call a {\it geometric
relation}.

We shall say that, under the action of~$\AAp¥$, the
address~$\1$ is the {\it heir} of the
address~$\1\1$, and, more generally, that
$\1\a$ is the heir of~$\1\1\a$. Inheritance
phenomena are quite general. Under the
action of~$\AAp\a$, for every~$\b$,
the address~$\a\0\0\b$ is the heir of~$\a\0\b$, the
address~$\a\0\1\b$ is the heir of~$\a\1\0\b$,
and~$\a\1\b$ is the heir of~$\a\1\1\b$
under~$\AAp\a$. Let us say that two
addresses~$\a, \b$ are incompatible, denoted $\a
\perp \b$, if neither is of prefix of the other, \ie, if
there exists~$\g$ such that
$\g\0$ is a prefix of~$\a$ and $\g\1$ is a prefix
of~$\b$, or {\it vice versa}. Then each
address~$\b$ with $\b \perp \a$ is its own heir
under the action of~$\AAp\a$. 

The argument leading to~\eqref{E:Heir} gives the
relation $\DD\g f \cdot \AAp\a = \AAp\a \cdot
\DD{\g'}f$ whenever $\g'$ is the heir of~$\g$
under~$\AAp\a$. In this way, we deduce a list of
geometric relations in~$\GMA$, namely
\begin{equation}\label{E:Geom}
\begin{cases}
\quad \AAp\b \opp \AAp\a 
= \AAp\a \opp \AAp\b
\mbox{\quad for $\b \perp \a$,}\\
\quad \AAp{\a\0\b} \opp \AAp\a
= \AAp\a \opp \AAp{\a\0\0\b},\\
\quad \AAp{\a\1\0\b} \opp \AAp\a
= \AAp\a \opp \AAp{\a\0\1\b},\\
\quad \AAp{\a\1\1\b} \opp \AAp\a
= \AAp\a \opp \AAp{\a\1\b}.
\end{cases}
\end{equation}

The geometric relations are rather trivial, and
we look for other, non-trivial relations
in~$\GMA$. As can be expected, MacLane's pentagon
enters the picture.

\begin{lemm}\label{L:Pent}
For each~$\a$, the following {\it pentagon
relation} holds in~$\GMA$:
\begin{equation} \label{E:Pent}
\AAp\a \opp \AAp\a = \AAp{\a\1} \opp  \AAp\a  \opp
\AAp{\a\0}.
\end{equation}
\end{lemm}

\begin{figure} [htb]
\begin{picture}(90, 26)(0, 0)
\put(0,0){\includegraphics{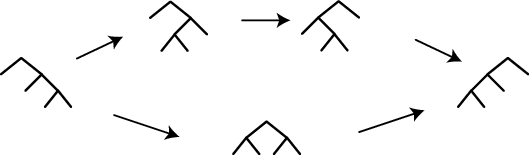}}
\put(12.5,20){$\AAp\1$}
\put(43.5,25){$\AAp¥$}
\put(74,20){$\AAp\0$}
\put(24,7){$\AAp¥$}
\put(64,7){$\AAp¥$}
\end{picture}
\caption{\smaller The pentagon relation}
\label{F:Pent}
\end{figure}

The verification (for $\a = \eadd$) is given in
Figure~\ref{F:Pent}. Keeping the same name for the
relations in~$\GMA$ and their counterparts
in~$\GGA$---hence in~$F$---we can summarize the
results as follows.

\begin{defi}
We denote by~$\GA$ the family of
all~$\AA\a$'s, and by~$\RA$ the family of all
geometry relations involving~$\AA¥$, namely the
translated copies of
\begin{gather}
\tag{$\RGperp$}
\AA{\0\a} \opp \AA{\1\b} = \AA{\1\b}
\opp \AA{\0\a},\\
\tag{$\RGA$}
\AA{\1\1\a} \AA¥ = \AA¥ \AA{\1\a},
\qquad
\AA{\1\0\a} \AA¥ = \AA¥ \AA{\0\1\a},
\qquad
\AA{\0\a} \AA¥ = \AA¥ \AA{\0\0\a};
\end{gather}
plus the pentagon relations, \ie, the translated
copies of
\begin{equation}
\tag{$\RP$}
\AA¥ \AA¥ =\AA\1 \AA¥ \AA\0.
\end{equation}
\end{defi}

\begin{prop} \label{P:RelF}
All relations of~$\RA$ are satisfied by
the elements~$\AA\a$ in~$\GGA$, \ie, in~$F$.
\end{prop}

\subsection{Constructing trees}

Our next aim is to prove that the relations
of Proposition~\ref{P:RelF} make a presentation
of~$F$. We apply the
method described in Section~\ref{S:Crit}, using the
partial action of~$\GGA$ on trees. We showed
that any family of trees containing trees of
arbitrary large size is discriminating, so,
according to Proposition~\ref{P:Criterion}, two
ingredients are needed, namely

- a family of trees~$\Disc$ containing trees
of unbounded size, and

- for every tree~$\tt$ in the orbit
of~$\Disc$, a distinguished word~$\HH\tt$
in~$\WW\GA$ connecting~$\tt$ with some
distinguished element of its orbit.

\begin{figure} [htb]
\begin{picture}(92, 29)(0,0)
\put(0,3){\includegraphics{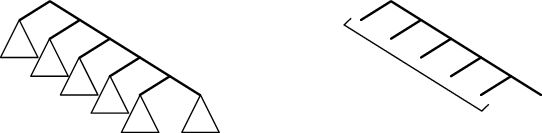}}
\put(48,12){$n-1$ vertices}
\put(2,13){$\tt_1$}
\put(7,10){$\tt_2$}
\put(11,7){$\pp$}
\put(21,0){$\tt_{n\!-\!1}$}
\put(33,0){$\tt_{n}$}
\put(9,27){$\vine{\tt_1, \pp, \tt_n}$}
\put(66,27){$\vine n$}
\end{picture}
\caption{\smaller The notation $\vine{\tt_1,
\pp, \tt_n}$ and the right vine $\vine n$}
\label{F:Vine}
\end{figure}

Both steps are easy: two trees are equal up to
associativity if and only if they have the same size, so
each family of trees containing exactly one
tree of size~$n$ for each~$n$ is convenient. In the
current case, we shall use the {\it right vines}
(or combs) of Figure~\ref{F:Vine}.

\begin{defi}
Let $\tt_1, \pp, \tt_n$ be trees. We put
$$\vine{\tt_1, \pp, \tt_n} = \tt_1 \op (\tt_2 \op
\pp \op (\tt_{n-1} \op \tt_n) \pp);$$
we define the {\it right vine} $\vine n$ to be
$\vine{\et, \pp, \et}$ with $n$~times~$\et$.
\end{defi}

With this notation, applying the
operator~$\AAp¥$ means replacing
$\vine{\tti, \ttii, \pp}$ with $\vine{\tti \op
\ttii, \pp}$. As there are vines of each size,
we immediately get:

\begin{lemm}
Vines form a discriminating subset of~$\Ti$ 
for the action of~$\GGA$.
\end{lemm}

If $\tt$ is a size~$n$ tree, there exists a
(unique) element of~$\GGA$ mapping the right
vine~$\vine n$ to~$\tt$: in order to
obtain~\eqref{E:Criterion} and possibly apply
Proposition~\ref{P:Criterion}, it suffices to
select a distinguished word~$\HH\tt$ representing
that element, \ie, to describe how
$\tt$ can be {\it constructed} from~$\vine n$
using associativity. Several solutions exist.
We give now an inductive definition that leads
to short computations, but requires that we
introduce two words~$\HH\tt, \HHH\tt$ for each
tree~$\tt$.

\begin{defi}
$(i)$ For $\ww$ a word involving letters indexed
by addresses, we denote by $\DD\a\ww$ the
word obtained by appending~$\a$ at the
beginning of each index; we use $\dd\ww$
for~$\DD\1\ww$.

$(ii)$ For each tree~$\tt$, we define two
words~$\HH\tt, \HHH\tt$ using the
inductive rules: 
\begin{align*}
\HH\tt 
&= \HHH\tt = \e 
&&\text{for $\tt$ of size~$1$,}\\
\HH\tt = \HHH\tti \opp \dd{\HH\ttii}, 
&\qquad 
\HHH\tt = \HHH\tti \opp \dd{\HHH\ttii} \opp \AA¥
&&\text{for $\tt = \tti \op \ttii$.}
\end{align*}
\end{defi}

The following characterization of the words~$\HH\tt$
and~$\HHH\tt$ is not needed in the sequel, but it
should make the construction concrete. Each
tree~$\tt$ admits a unique decomposition in terms
of the basic tree~$\et$. Besides the algebraic
notation $\tti \opp \ttii$ for the product
of~$\tti$ and~$\ttii$, we can also use the right
Polish notation in which this product is denoted
$\tti \ttii \oppol$. For instance, the
Polish expression of $\et \op ((\et
\op \et) \op \et)$ is $\et \et \et \oppol \et
\oppol \oppol$. In the next proposition, a
length~$\ell$ word~$\ww$ is considered as a
sequence of symbols indexed by $\{1, \pp,
\ell\}$, and $\ww(p)$ denotes the $p$th
symbol in~$\ww$.

\begin{prop} \label{P:Polo}
For~$\ww$ a length~$\ell$ word and $0 \le p
\le \ell$, define the defect~$\defect\ww p$
of~$p$ in~$\ww$  by the rules: $\defect\ww0
= -1$, $\defect\ww{p} = \defect\ww{p-1}  -1$
for $\ww(p) = \oppol$, and $\defect\ww{p}
= \defect\ww{p-1} +1$ otherwise. Then, for each
tree~$\tt$, the word~$\HHH\tt$ is
obtained from the Polish expression
of~$\tt$ by deleting the symbols~$\et$, and
replacing each defect~$i$ symbol~$\oppol$
with~$\AA{\1^i}$. The word~$\HH\tt$ is
obtained similarly, except that the final
symbols~$\oppol$, \ie, those followed by
no~$\et$, do not contribute.
\end{prop}

\begin{proof}
It is standard that a word~$\ww$ is the Polish
expression of a tree if and only if
the defect of each symbol is nonnegative, and the
defect of the last symbol is~$0$. For~$\tt$ a tree,
define the {\it enhanced decomposition} of~$\tt$
to be the Polish expression with the defect of
each symbol appended. Then the
enhanced decomposition of~$\tti \op
\ttii$ is the enhanced decomposition
of~$\tti$, followed by the enhanced decomposition
of~$\ttii$ with all defects shifted by~$1$, followed by
the symbol~$\oppol$ with $0$ defect. So, the
enhanced decomposition and the word~$\HHH\tt$ obey
parallel inductive rules. Hence, as the
correspondence of the proposition clearly holds for the
basic tree~$\et$, it inductively holds for every tree. A
similar argument gives the connection
between~$\HHH\tt$ and~$\HH\tt$.
\end{proof}

For instance, for $\tt = \et \op ((\et
\op \et) \op \et)$, the enhanced decomposition
of~$\tt$ is $\et 0 \et 1 \et 2 \oppol 1 \et 2
\oppol 1 \oppol 0$, and a direct translation
yields $\HHH\tt = \AA\1 \AA\1 \AA¥$, and
$\HH\tt = \AA\1$ (the last two symbols~$\oppol$
are dismissed). A consequence of
Proposition~\ref{P:Polo} is that, for each
tree~$\tt$, we have
\begin{equation}
\HHH\tt = \HH\tt \opp \AA{1^{\hh-1}} \pp
\AA\1 \AA¥,
\end{equation}
where $\hh$ is the length of the rightmost branch
in~$\tt$.

We aim at proving that the trees~$\vine n$ and the
words$~\HH\tt$ satisfy the requirements
of~Proposition~\ref{P:Criterion} and therefore
lead to a presentation of~$\GGA$, \ie, of~$F$. In
the sequel, we use mixed expressions like
$\vine{p, \tt, q, \pp}$ where
$p, q$ are numbers and~$\tt$ is a tree to mean
$\vine{\et ,
\pp, \et, \tt, \et, \pp, \et, \pp}$ with
$p$~$\et$ in the first block and $q$ in the
second. 

\begin{lemm} \label{L:ConF}
For each size~$\nn$ tree~$\tt$ and each
tree~$\ttt$, we have
\begin{equation}
\begin{CD}
\vine n
@>\HH\tt>>
\tt
\end{CD}
\text{\qquad and \qquad}
\begin{CD}
\vine{n, \ttt}
@>\HHH\tt>>
\vine{\tt, \ttt}
\end{CD},
\end{equation}
\ie, $\HH\tt$ constructs~$\tt$
from~$\vine\nn$, and $\HHH\tt$
constructs~$\tt \op \ttt$ from~$\vine{\nn,
\ttt}$.
\end{lemm}

\begin{figure} [htb]
\begin{picture}(117,45)(0, 0)
\put(0,0){\includegraphics{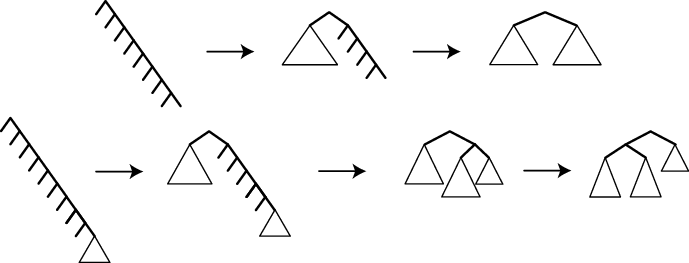}}
\put(20,45){$\vine\nn$}
\put(59,42){$\vine\nnii$}
\put(52,30){$\tti$}
\put(85,30){$\tti$}
\put(95,30){$\ttii$}
\put(36,39){$\HHH\tti$}
\put(70,39){$\dd\HH\ttii$}
\put(20,0){$\ttt$}
\put(46,1){$\ttt$}
\put(82,9.5){$\ttt$}
\put(113,12){$\ttt$}
\put(31,10){$\tti$}
\put(71,10){$\tti$}
\put(101,8){$\tti$}
\put(77,8){$\ttii$}
\put(108,8){$\ttii$}
\put(17,19){$\HHH\tti$}
\put(55,19){$\dd\HHH\ttii$}
\put(91,19){$\AA¥$}
\end{picture}
\caption{\smaller For~$\tt$ a tree of size~$n$, 
the word~$\HH\tt$ describes how to construct~$\tt$
from~$\vine n$, and $\HHH\tt$ describes
how to construct $\vine{\tt, \ttt}$
from~$\vine{\nn, \ttt}$; the figure illustrates
the inductive argument for $\tt = \tti \op
\ttii$}
\label{F:Blue}
\end{figure}

\begin{proof}
We use induction on~$n$. For $n = 1$, the result is
obvious. Otherwise, assume $\tt = \tti \op \ttii$, and let
$\nni$ and~$\nnii$ be the respective sizes
of~$\tti$ and~$\ttii$. Then we have
$\HH\tt =
\HHH\tti \opp
\dd{\HHH\ttii}$. By induction hypothesis,  $\HHH\tti$
maps~$\vine n$, \ie, $\vine{\nni, \nnii}$,
to~$\tti \op \vine\nnii$, \ie, $\vine{\tti,
\nnii}$. Then, by induction hypothesis again,
$\HH\ttii$ maps~$\vine\nnii$ to~$\ttii$, hence
$\dd{\HH\ttii}$ maps $\tti \op \vine\nnii$ to
$\tti
\op \ttii$. So $\HH\tt$ maps~$\vine n$ to~$\tti \op \ttii$,
\ie, to~$\tt$ (Figure~\ref{F:Blue} top):
$$\begin{CD}
\vine n = \vine{\nni, \nnii}
@>\HHH\tti>>
\vine{\tti, \nnii}
@>\dd{\HH\ttii}>>
\vine{\tti, \ttii} = \tt.
\end{CD}$$
Similarly, we have $\HHH\tt = \HHH\tti \opp
\dd{\HHH\ttii} \opp \AA¥$, and the diagram is
now:
$$\begin{CD}
\vine{\nn, \ttt} = \vine{\nni, \nnii, \ttt}
@>\HHH\tti>>
\vine{\tti, \nnii, \ttt}
@>\dd{\HHH\ttii}>>
\vine{\tti, \ttii, \ttt} 
@>{\AA¥}>>
\vine{\tti \op \ttii, \ttt} = \vine{\tt, \ttt}
\end{CD},$$
as is easily checked on Figure~\ref{F:Blue}
bottom.
\end{proof}

So Condition~\eqref{E:Criterion} is satisfied. Then
Proposition~\ref{P:Criterion} tells us that
a family of relations involving the
generators~$\AA\a$ makes a presentation
of~$\GGA$ if and only if it contains enough
relations to make the words $\HH\ttt$ and
$\HH\tt \opp \AA\a$ equivalent whenever
$\ttt$ is the image of~$t$ under~$\AA\a$. As we
will show now, this is the case for the
relations~$\RA$ of Proposition~\ref{P:RelF}. Due
to our inductive construction, it is
convenient to prove two results simultaneously,
namely one for~$\HH\tt$ and one
for~$\HHH\tt$. Note that the argument proving
$\cl{\HH\ttt} = \cl{\HH\tt} \opp
\AA\a$ when $\AA\a$ maps~$\tt$ to~$\ttt$
similarly proves $\cl{\HHH\ttt} = \cl{\HHH\tt}
\opp \AA{\0\a}$, as, writing $n$ for the common size
of~$\tt$ and~$\ttt$, both operators map~$\vine
N$ to~$\vine{\ttt, N-n}$ for
$N > n$. In the sequel, we use
$\eqG$ and~$\eqP$ to indicate that we
specifically use a  geometric (\ie, a twisted
commutation) or a pentagon relation.

\begin{lemm}\label{L:BluF} 
Assume $\ttt = \tt \act \AA\a$. Then we have
\begin{equation} \label{E:BluF}
\HH\ttt \eqA \HH\tt \opp \AA\a
\mbox{\quad and \quad}
\HHH\ttt \eqA \HHH\tt \opp \AA{\0\a}.
\end{equation}
\end{lemm}

\begin{proof}
We use induction on the length of~$\a$ as a sequence
of~$0$'s and~$1$'s. Assume first that $\a$ is the
empty address. The hypothesis that $\ttt  =
\tt \act \AA¥$ holds, \ie, that the
operator~$\AAp¥$ maps~$\tt$ to~$\ttt$, means
that there exist
$\tti, \ttii, \ttiii$ such that $\tt$ is $\tti \op (\ttii
\op \ttiii)$ and
$\ttt$ is $(\tti \op \ttii) \op \ttiii$. Then we find
\begin{gather*}
\HH\ttt
= \HHH\tti \opp \dd{\HHH\ttii} \opp \AA¥ \opp
\dd{\HH\ttiii} 
\eqG
 \HHH\tti \opp \dd{\HHH\ttii} \opp
\ddd 2{\HH\ttiii} \opp \AA¥
= \HH\tt \opp \AA¥,\\
\HHH\ttt
= \HHH\tti \opp \dd{\HHH\ttii} \opp \AA¥ \opp
\dd{\HHH\ttiii} \opp \AA¥ 
\eqG \HHH\tti \opp \dd{\HHH\ttii} \opp
\ddd 2{\HHH\ttiii}  \opp \AA¥ \AA¥ 
\eqP \HHH\tti \opp \dd{\HHH\ttii} \opp
\ddd 2{\HHH\ttiii}  \opp \AA\1 \AA¥ \AA\0
= \HHH\tt \opp \AA\0.
\end{gather*}
Assume now $\a = \0\b$. The hypothesis that
$\AAp\a$ maps~$\tt$ to~$\ttt$ means that there
exist $\tti, \ttii, \ttti$ such that $\tt$
is~$\tti \op
\ttii$, $\ttt$ is $\ttti \op \ttii$, and $\AAp\b$
maps~$\tti$ to~$\ttti$. Using the induction
hypothesis~$(IH)$, we find
\begin{align*}
\HH\ttt
= \HHH\ttti \opp \dd{\HH\ttii}
&\eqIH \HHH\tti \opp \AA{\0\b}
\opp \dd{\HH\ttii}
\eqG \HHH\tti \opp \dd{\HH\ttii}
 \opp \AA{\0\b}
= \HH\tt \opp \AA\a,\\
\HHH\ttt
= \HHH\ttti \opp \dd{\HHH\ttii} \opp \AA¥ 
&\eqIH \HHH\tti \opp \AA{\0\b}
 \opp \dd{\HHH\ttii} \opp \AA¥ \\  
&\eqG \HHH\tti \opp \dd{\HHH\ttii}
 \opp \AA{\0\b}\AA¥
\eqG \HHH\tti \opp \dd{\HHH\ttii}
 \opp \AA¥ \AA{\0\0\b}
= \HHH\tt \opp \AA{\0\a}.
\end{align*}
Finally, assume $\a = \1\b$. With similar notation, we
have $\tt = \tti \op \ttii$ and $\ttt = \tti \op \tttii$
with $\AAp\b$ mapping~$\ttii$ to~$\tttii$, and we
find now
\begin{gather*}
\HH\ttt
= \HHH\tti \opp \dd{\HH\tttii}
\eqIH 
\HHH\tti \opp \dd{\HH\ttii} \opp  \AA{\1\b}
= \HH\tt \opp \AA\a,\\
\HHH\ttt
= \HHH\tti \opp \dd{\HHH\tttii} \opp \AA¥ 
\eqIH 
\HHH\tti \opp \dd{\HHH\ttii} \opp
\AA{\1\0\b}\AA¥ 
\eqG 
\HHH\tti \opp \dd{\HHH\ttii} \opp
\AA¥ \AA{\0\1\b} 
 = \HHH\tt \opp \AA{\0\a},
\end{gather*}
which completes the proof.
\end{proof}

Applying Proposition~\ref{P:Criterion}, we deduce:

\begin{prop} \label{P:PreF}
The relations~$\RA$, \ie, the geometric relations
for~$\AA¥$ plus the pentagon relations, make a
presentation of the group~$\GGA$, \ie,~$F$, in
terms of the generators~$\AA\a$.
\end{prop}

\subsection{The standard presentation}

There is a well-known presentation
of~$F$ in terms of an infinite sequence of 
generators, usually denoted~$x_i$,  indexed by
nonnegative integers~\cite{CFP}. It is easy to
establish the connection between these
generators and our current generators~$\AA\a$
and, using Proposition~\ref{P:Criterion} again, to
re-obtain the standard presentation of~$F$ as a
direct corollary. 

\begin{defi}
$(i)$ For $i \ge 1$, we put $\aa i =
\AA{\1^{i-1}}$, and we denote by~$\Ga$ the
family of all~$\aa i$'s.

$(ii)$ We denote by~$\Ra$ the subfamily
of~$\RA$ consisting of those relations in~$\RA$
that involve the generators of~$\Ga$ exclusively,
namely the relations $\aa i \aa{j-1} = \aa j
\aa i$ for $j \ge i+2$.
\end{defi}

\begin{prop} \label{P:Pref}
The set~$\Ga$ generates~$\GGA$, \ie,~$F$,
and the relations~$\Ra$ make a presentation
of~$\GGA$ in terms of the generators~$\aa i$.
\end{prop}

\begin{proof}
By construction, the words~$\HH\tt$ belong
to~$\WW{\Ga}$, and we can apply
Proposition~\ref{P:Criterion} to the family~$\Ga$.
So, in order to prove that $\Ra$ makes a
presentation,  it suffices to check that the
relations of~$\Ra$ are sufficient to establish the
equivalence of $\HH\ttt$ and~$\HH\tt \opp \aa i$
when $\aa i$ maps~$\tt$ to~$\ttt$. Looking at the
proof of Lemma~\ref{L:BluF} immediately shows that
this is true.
\end{proof}

As the $\aa i$'s generates~$F$, each~$\AA\a$
can be expressed in terms of the~$\aa
i$'s. For $\a$ an address containing at least
one~$\0$, say
$\a = \1^p\0\0^{e_0} \1\0^{e_1} \1 \pp \1 \0^{e_q}$ with
$p, q, e_0, \pp, e_q \ge 0$, one can check
$$\AA\a = (\aa{p+1}^{e_0+1} \aa{p+2}^{e_1+1}
\pp
\aa{p+q+1}^{e_q + 1}) (\aa{p+q+1}
\aa{p+q+2}\inv) 
 (\aa{p+1}^{e_0+1} \aa{p+2}^{e_1+1} \pp
\aa{p+q+1}^{e_q + 1})\inv.
$$
For instance, for $\a = \0\1\1\0\0$, we
have $\AA{\0\1\1\0\0} = \aa1 \aa2 \aa3^4 \aa4\inv
\aa3^{-3} \aa2\inv \aa1\inv$. As was noted in the
introduction, the current~$\aa i$ corresponds
to~$x_{i-1}\inv$ in literature about~$F$.

\subsection{The lattice structure of~$F$}

It is known that $F$ is a finitely presented
group, generated by the two elements~$\aa1,
\aa2$. Using infinite presentations has
disadvantages, and it may seem strange to replace
the infinite family~$\Ga$, which requires a very
simple set of relations, with the still larger
family~$\GA$ that involves a seemingly more
complicated set of relations. However, one of the
interests of the presentation~$(\GA; \RA)$ of~$F$
is that it is more symmetric, giving the same role
to the left and right directions, contrary
to~$\Ga$ that priviledges the right one.

In particular, considering the
generators~$\AA\a$ makes it natural to
introduce the submonoid~$\Fp$ of~$F$ generated
by these elements. Using a monoid version
of Proposition~\ref{P:Criterion} and a convenient
combinatorial methods, one can show that
$\Fp$ admits, as a monoid, the
presentation~$(\GA, \RA)$ and that it is
isomorphic to the geometry monoid of oriented
associativity~$\GGp\Ass$ defined as~$\GG\Ass$ but
considering the positive operators~$\AAp\a$
only~\cite{Dfg}. Contrary to the submonoid of~$F$
generated by the elements~$\aa i$, the
monoid~$\Fp$ admits both left and right least
common multiples, and one obtains in this way a
double lattice structure on~$F$. 

Another interest of considering the
generators~$\AA\a$ is that the associated Cayley
graph is closely connected with Stasheff's
associahedra: essentially, the graph is a direct
limit of the associahedra, which appear as the
orbits of the (partial) action of~$F$ on binary
trees. These aspects, first described by
P.~Greenberg in~\cite{Gre}, will be
further investigated in~\cite{Dhz}.

\section{The geometric presentation of
Thompson's group~$V$}
\label{S:PreV}

Our approach to Thompson's group~$F$ was based on
its connection with the associativity. We now
develop a similar approach for Thompson's
group~$V$. The latter appears
when the commutativity law $xy = yx$ is
added. As in Section~\ref{S:PreF}, the geometry of
the commutativity operators leads to a natural
presentation: in addition to the
geometric and pentagon relations, the only new
relations are the MacLane--Stasheff hexagon
relations, plus some torsion relations.

\subsection{The geometry monoid of a family
of algebraic laws} \label{S:Iden}

The approach developed in Section~\ref{S:PreF} for
the special case of associativity extends to
arbitrary algebraic laws. The general form of an
identity~$\IId$ is $\tau_\smm = \tau_\smp$,
where $\tau_\smm, \tau_\smp$ are formal 
combinations of variables, or, equivalently,
coloured trees. For each such~$\IId$, we can
consider the partial operator~$\IIdp¥$ on~$\Ti$
such that a tree~$\tt$ belongs to the domain
of~$\IIdp¥$ if it can be written as
$\tau_\smm^\subst$ for some substitution~$\subst$,
and, then, define
$\tt \act \IIdp¥$ to
be~$\tau_\smp^\subst$. The operator~$\IIdm¥$ is
defined symmetrically, and, as above, we denote
by~$\IIdpm\a$ the translated
copy~$\DD\a{\IIdpm¥}$, \ie, the result of
letting~$\IIdpm¥$ act on the $\a$-subtree.

\begin{defi}
For $\IId, \JJd, \pp$ algebraic laws, we
define the {\it geometry monoid} of~$\IId, \JJd,
\pp$, denoted $\GM{\IId, \JJd, \pp}$, to
be the monoid generated by all partial
operators~$\IIdpm\a, \JJdpm\a, \pp$ 
\end{defi}

Thus, the monoid~$\GMA$ of
Section~\ref{S:PreF} is the geometry monoid of
the associativity law.
Formally, the definition of the
operators~$\IIdp\a$ and, therefore, of the
geometry monoid, depends on the considered
family of (coloured) trees. We shall
forget about this here, which amounts to
assuming that we work once for all inside a
sufficiently large family of coloured
trees~$\Ti$. 

In this framework, it is obvious that two
trees~$\tt, \ttt$ are $\{\IId, \JJd,
\pp\}$-equal if and only if some element
of~$\GM{\IId, \JJd, \pp}$ maps~$\tt$ to~$\ttt$.
At this degree of generality, we cannot expect
really deep results, and going further
requires to restrict the considered laws.
An unpleasant phenomenon is that,
in general, the geometry monoid~$\GM{\IId, \JJd,
\pp}$ contains the empty mapping, \ie, there exist
products of operators~$\IIdpm\a$, $\JJdpm\b$, \pp
applying to no tree, typically because the
labels cannot be compatible. This however is
excluded when the laws are simple
enough.

\begin{defi}
A law $\tau_\smm =  \tau_\smp$ is said to be
{\it linear} if the same variables occur
in~$\tau_\smm$ or in~$\tau_\smp$ and each of
them occurs exactly once.
\end{defi}

The associativity law
$x(yz) = (xy)z$ is linear, as $x$, $y$, and~$z$
occur only once on each side,
while the self-distributivity law $x(yz) =
(xy)(xz)$ is not, as $x$ is repeated twice on the
right.

\begin{lemm} \label{L:Seed}
Assume that $\IId, \JJd, \pp$ are linear laws.
Then each operator~$f$ in~$\GMIJ$ admits a seed
consisting of injective trees, \ie, there
exists a pair of injective trees~$(\tt, \ttt)$
such that, as a pair of trees,  $f$~is the set
of all substitutes of~$(\tt, \ttt)$.
\end{lemm}

\begin{proof}
The point is that, if $\tt_1, \tt_2$ are injective
trees, then there always exists
substitutions~$\substi, \substii$ such that
$\tt_1^{\substi}$ and $\tt_2^{\substii}$ are
equal, which need not be the case when
some labels in $\tt_1$ or~$\tt_2$ occur twice.
Then the substitutions may be chosen so that the
common skeleton of $\tt_1^{\substi}$ and
$\tt_2^{\substii}$ is the union of the skeletons
of~$\tt_1$ and~$\tt_2$, and the proof is the
same as for Lemma~\ref{L:DomF}.
\end{proof} 

In the previous case, Lemma~\ref{L:Mono} applies,
and, as in the case of~$\GMA$, it leads to a group.

\begin{prop} \label{P:GeoG}
Let $\IId, \JJd, \pp$ be linear algebraic laws.
Then near-equality is a congruence on~$\GMIJ$,
and the quotient-monoid is a group. The
operators~$\IIdpm\a$, $\JJdpm\a$, \pp
induce a partial action of this  group
on trees. Injective trees form a discriminating family
for this action.
\end{prop}

\begin{defi}
Under the above hypothesis, the
group~$\GMIJ/\!\approx$ is called the {\it
geometry group} of the laws~$\IId, \JJd,
\pp$, and it is denoted~$\GGIJ$.
\end{defi}

(Here, we restrict to laws that involve a
single binary operation. A similar approach is
of course possible for laws involving more than
one operation, at the expense of considering
trees in which the internal nodes are
marked with operation symbols and their
degree is adjusted to the arity of the
operation.)

\subsection{Commutativity operators
and Thompsons's group~$V$}

The commutativity law
\begin{equation}
\tag{$\Com$}
xy=yx
\end{equation}
is eligible for the previous approach. Here the
basic operator is the operator exchanging the
left and the right subtrees of a tree:

\begin{defi}
We denote by~$\CCp¥$ the (partial) operator that
maps every tree of the form $\tti \op \ttii$ to the
corresponding tree $\ttii \op \tti$. For each
address~$\a$, we put~$\CCp\a = \DD\a{\CCp¥}$.
We define~$\GMAC$ to be the monoid generated
by all operators $\AAp\a$ and~$\CCp\a$ and
their inverses.
\end{defi}

Associativity and commutativity are linear
laws, hence Lemma~\ref{L:Seed} and,
therefore, Proposition~\ref{P:GeoG} apply. So,
near-equality is a congruence on the
monoid~$\GMAC$, and we obtain a group,
denoted $\GGAC$ by identifying
near-equal operators. As in Section~\ref{S:PreF},
we shall use~$\AA\a$ for the class of~$\AAp\a$
in~$\GGAC$, and, similarly, $\CC\a$ for the
class of~$\CCp\a$. We still denote by~$\GA$ the
family of all~$\AA\a$'s, and, similarly, we
use~$\GC$ for the family of all~$\CC\a$'s.

\begin{prop} \label{P:GroV}
The group~$\GGAC$ is (isomorphic to) Thompson's
group~$V$, \ie, $V$ is the geometry group of
associativity and commutativity.
\end{prop}

\begin{figure} [htb]
\begin{picture}(100, 20)(0, 0)
\put(0,0){\includegraphics{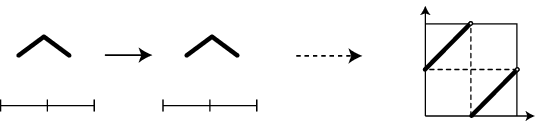}}
\put(20,14){$\CCp¥$}
\put(2,7){$1$}
\put(10.5,7){$2$}
\put(31,7){$2$}
\put(39.5,7){$1$}
\end{picture}
\caption{\smaller From~$\GGAC$ to~$V$: the
action of~$\CCp¥$}
\label{F:FromGCToV}
\end{figure}

\begin{proof}
(Figure~\ref{F:FromGCToV}) We associate with each
element of~$\GGAC$ an element of~$V$, \ie, a
piecewise linear mapping of~$[0, 1]$ into
itself as we did for~$\GGA$ and~$F$ in
Section~\ref{S:PreF}: we associate to each tree a
dyadic partition of~$[0, 1]$, and we map~$f$ to
the piecewise linear function that maps the
partition associated to~$\ttt$ to the partition
associated to~$\tt$, where $(\tt, \ttt)$ is a seed
for~$f$---we again reverse the orientation to
obtain a homomorphism  with composition---and
interpolates the values. The latter homomorphism
is surjective since, as was shown in
Section~\ref{S:PreF}, its image includes~$F$, and
it contains the mappings denoted~$C$ and~$\pi_0$
in~\cite{CFP}, which correspond to~$\AA¥ \CC\0
\AA¥\inv$ and
$\AA¥ \CC\0 \AA¥\inv \CC\1$ respectively.
\end{proof}

In the sequel, we identify $\GGAC$ and~$V$. 

\subsection{Guessing relations in~$\GMAC$}

As in the case of~$\GMA$, geometric inheritance
provides twisted commutation relations
in~$\GMAC$, and therefore in the
group~$\GGAC$. First, we observed that
$\dd^2 f \cdot \AAp¥ =
\AAp¥ \cdot \dd f$ holds for every mapping~$f$,
and used it for~$f = \AAp\a$ to obtain
$\AAp{\1\1\a} \opp \AAp¥ =
\AAp¥ \opp \AAp{\1\a}$. Applying it now to $f =
\CCp\a$, we deduce 
$\CCp{\1\1\a} \opp \AAp¥ = \AAp¥ \opp
\CCp{\1\a}$ similarly. In this way, using~$\BB\a$
to represent either~$\AAp\a$ or~$\CCp\a$, we
obtain the following relations:
\begin{equation}\label{E:Geo1}
\begin{cases}
\quad \BB\b \AAp\a 
= \AAp\a \BB\b
\mbox{\quad whenever $\b \perp \a$ holds,}\\
\quad \BB{\a\1\1\b} \AAp\a
= \AAp\a \BB{\a\1\b},\\
\quad \BB{\a\1\0\b} \AAp\a
= \AAp\a \BB{\a\1\0\b},\\
\quad \BB{\a\0\b} \AAp\a
= \AAp\a \BB{\a\0\0\b}.
\end{cases}
\end{equation}
New inheritance phenomena appear
with~$\CCp\a$: its action on a
tree~$t$ exchanges the $\a\0$- and
$\a\1$-subtrees of~$\tt$, and we deduce the
following relations, where $\BB\a$ still stands
for~$\AAp\a$ or~$\CCp\a$:
\begin{equation}\label{E:Geo2}
\begin{cases}
\quad \BB\b \CCp\a 
= \CCp\a \BB\b
\mbox{\quad whenever $\b \perp \a$ holds,}\\
\quad \BB{\a\0\b} \CCp\a
= \CCp\a \BB{\a\1\b},\\
\quad \BB{\a\1\b} \CCp\a
= \CCp\a \BB{\a\0\b}.
\end{cases}
\end{equation}
The relations mentioned in~\eqref{E:Geo1}
and~\eqref{E:Geo2} will be called the $\AA¥$-
and $\CC¥$-{\it geometric} relations, respectively.
Apart from the geometric relations, we
know that the pentagon relations, \ie,
\begin{equation}
\AAp\0 \AAp¥ \AAp\1 = \AAp¥{}^2
\end{equation}
and its shifted copies, are satisfied in~$\GMA$, hence
in~$\GMAC$. Two more types arise now.

\begin{lemm}
The following relations and their translated copies
hold in~$\GMAC$:
\begin{gather}
\label{E:Hexa}
\AAp¥ \CCp¥ \AAp¥ = \CCp\0 \AAp¥ \CCp\1,\\
\label{E:Tors}
\CCp¥{}^2 \approx \id.
\end{gather}
\end{lemm}

\begin{figure}[htb]
\begin{picture}(108,34)(0, 0)
\put(0,3){\includegraphics{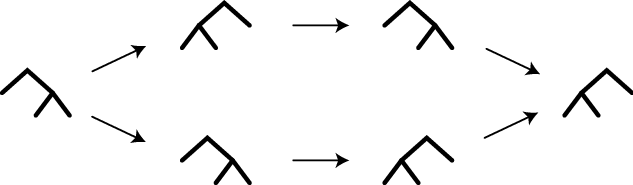}}
\put(18,26){$\AAp¥$}
\put(53,32){$\CCp¥$}
\put(87,26){$\AAp¥$}
\put(17,9){$\CCp\1¥$}
\put(53,3){$\AAp¥$}
\put(86,9){$\CCp\0¥$}
\put(-1,16){$\scriptstyle1$}
\put(5,12.5){$\scriptstyle2$}
\put(11,12.5){$\scriptstyle3$}
\put(30.5,4){$\scriptstyle1$}
\put(36,0.5){$\scriptstyle3$}
\put(41.5,0.5){$\scriptstyle2$}
\put(65,0.5){$\scriptstyle1$}
\put(70,0.5){$\scriptstyle3$}
\put(76,4){$\scriptstyle2$}
\put(95,12){$\scriptstyle3$}
\put(101,12){$\scriptstyle1$}
\put(107,16){$\scriptstyle2$}
\put(30.5,23){$\scriptstyle1$}
\put(35,23){$\scriptstyle2$}
\put(42,27){$\scriptstyle3$}
\put(65,27){$\scriptstyle3$}
\put(70,23){$\scriptstyle1$}
\put(76,23){$\scriptstyle2$}
\end{picture}
\caption{\smaller The hexagon relation}
\label{F:Hexa}
\end{figure}

The verification for~\eqref{E:Hexa}, which 
corresponds to two ways of going from~$(\tti \op
\ttii) \op
\ttiii$ to~$(\ttii \op \ttiii) \op \tti$, is given
in Figure~\ref{F:Hexa}. The involutivity
of~$\CCp¥$ is obvious---but, as
$\CCp¥$ is defined only on those trees that are
not~$\et$, we obtain a $\approx$-relation, not an
equality.

\begin{defi} 
Let~$\RAC$ consist of all
$\AA¥$- and $\CC¥$-geometric
relations, \ie, the translated
copies of
\begin{gather}
\tag{$\RGperp$}
\BB{\0\a} \opp \BBB{\1\b} = \BBB{\1\b}
\opp \BB{\0\a},\\
\tag{$\RGA$}
\BB{\1\1\a} \opp \AA¥ = \AA¥ \opp \BB{\1\a},
\qquad
\BB{\1\0\a} \opp \AA¥ = \AA¥ \opp \BB{\0\1\a},
\qquad
\BB{\0\a} \opp \AA¥ = \AA¥ \opp \BB{\0\0\a},\\
\tag{$\RGC$}
\BB{\0\a} \opp \CC¥ = \CC¥ \opp \BB{\1\a},
\qquad 
\BB{\1\a} \opp \CC¥ = \CC¥ \opp \BB{\0\a},
\end{gather}
with $\BB¥, \BBB¥ = \AA¥$
or~$\CC¥$, plus the pentagon relations,
\ie, the translated copies of
\begin{equation}
\tag{$\RP$}
\AA¥ \AA¥ =\AA\1 \AA¥ \AA\0,
\end{equation}
plus the hexagon relations, defined to be the
translated copies of
\begin{equation}
\tag{$\RH$}
\AA¥ \CC¥ \AA¥ =\CC\1 \AA¥ \CC\0
\text{\quad and \quad} 
\AA¥\inv \CC¥ \AA¥\inv =\CC\0 \AA¥\inv
\CC\1.
\end{equation}
\end{defi}

The relations in~$\GMAC$ induce
similar relations in~$\GGAC$, \ie, in~$V$. So,
we may state:

\begin{prop} \label{P:RelV}
All relations in~$\RAC$ plus the torsion relations
$\CC\a^2 = 1$ are satisfied by the elements
of~$\GA$ and~$\GC$ in the group~$\GGAC$, \ie,
in~$V$.
\end{prop}

Distinguishing two hexagon relations, which are
equivalent when the torsion relations~$\CC\a^2 =
1$ are present, may seem strange. The reason is
that we consider  a torsion-free version
of~$V$ in Section~\ref{S:PreB}, and it is
convenient to keep track of the torsion relations
from now on.

\subsection{Restricting the family of
generators}

As in the case of~$F$, we shall consider two
families of generators for the group~$V$:
besides the families~$\GA$ and~$\GC$
comprising all~$\AA\a$'s and~$\CC\a$'s, we
shall also consider the proper subfamilies
corresponding to right branch addresses.

\begin{defi}
For $i \ge 1$, we put $\cc i =
\CC{\1^{i-1}}$. We denote by~$\Gc$ the family
of all~$\cc i$'s.
\end{defi}

Thus $\cc i$ is an exact counterpart
to~$\aa i$. We now list some relations
satisfied by the elements of~$\Ga$
and~$\Gc$ in~$\GGAC$. A disadvantage of 
restricting the families of generators is that
expressing the geometric phenomena is less simple
than with the whole families~$\GA$ and~$\GC$.

\begin{defi}
We define $\Rac$ to consist of the
following relations:
\begin{gather}
\label{E:Rac1}
\aa i \bb{j-1}
= \bb j \aa i
\qquad\text{for $j \ge i+2$ and
$\bb¥ = \aa¥$ or~$\cc¥$},\\
\label{E:Rac2}
\cc i  \aa i\inv \cc{i+1}\inv\bb j
= \bb j  \cc i \aa i\inv \cc{i+1}\inv
\qquad\text{for $j \ge i+2$ and
$\bb¥ = \aa¥$ or~$\cc¥$},\\
\label{E:Rac3}
\aa{i+1} \aa i \cc i^e \aa{i+1}
= \aa i^2 \cc i^e
\qquad\text{for $e = \pm 1$},\\
\label{E:Rac4}
\aa i \cc i \cc{i+1} \aa i
= \cc{i+1} \cc i,\\
\label{E:Rac5}
\cc{i+1} \cc i \aa i\inv \cc{i+1}
= \cc i \aa i\inv \cc i \aa i\inv.
\end{gather}
\end{defi}

\begin{lemm}
All relations in~$\Rac$ follow from~$\RAC$ (and
the definitions $\aa i = \AA{\1^{i-1}}$,
$\cc i = \CC{\1^{i-1}}$).
\end{lemm}

\begin{proof}
It is sufficient to establish the relations
for $i = 1$ and then use~$\dd^{i-1}$
to deduce the general version.
Relations~\eqref{E:Rac1} 
and~\eqref{E:Rac2} are of purely
geometric nature: \eqref{E:Rac1} is a
$\AA¥$-geometric
relation, and \eqref{E:Rac2} follows from
$$\CC¥  \AA¥\inv \CC\1\inv\BB{\1\1\a}
\eqG
\CC¥ \AA¥\inv \BB{\1\0\a} \CC\1\inv
\eqG
\CC¥ \BB{\0\1\a} \AA¥\inv \CC\1\inv 
\eqG
\BB{\1\1\a} \CC¥ \AA¥\inv \CC\1\inv,$$ 
which is valid both for $\BB¥ = \AA¥$ or~$\CC¥$.
Relations~\eqref{E:Rac3} use the
pentagon relations:
$$\AA\1 \AA¥ \CC¥^e \AA\1
\eqG
\AA\1 \AA¥ \AA\0 \CC¥^e 
\eqP
\AA¥^2 \CC¥^e.$$
Finally, appealing to the hexagon relations, we
find
\begin{gather*}
\AA¥ \CC¥ \CC\1
\eq¥
\AA¥ \CC¥ \AA¥ \AA¥\inv \CC\1
\eqH
\CC\1 \AA¥ \CC\0  \AA¥\inv \CC\1
\eqH
\CC\1 \AA¥ \AA¥\inv \CC¥ \AA¥\inv
\eq¥
\CC\1 \CC¥ \AA¥\inv,\\
\CC\1 \CC¥ \AA¥\inv \CC\1
\eqG
\CC¥ \CC\0 \AA¥\inv \CC\1
\eqH
\CC¥ \AA¥\inv \CC¥ \AA¥\inv,
\end{gather*}
which gives \eqref{E:Rac4} and \eqref{E:Rac5}.
\end{proof}

\subsection{Constructing trees}

We aim at proving that the
relations~$\RAC$ and~$\Rac$ make presentations
of the group~$V$. As in the case of~$F$, we
shall use the criterion of
Proposition~\ref{P:Criterion}. So, as in
Section~\ref{S:PreF}, the point is to
introduce for each tree~$t$ a distinguished
word~$\HH\tt$ that describes the construction
of~$t$ from some distinguished tree in its
$V$-orbit.

In contrast to the case of associativity,
considering commutativity requires that we take
labels into account. Indeed, uncoloured trees are
not discriminating: for instance,
the operators~$\id$ and~$\CC¥$ are not
(near)-equal, but both fix~$\et \op \et$.
We use coloured versions of the right
vines~$\vine n$.

\begin{figure}[htb]
\begin{picture}(0, 19)(33, 0)
\put(0.5,3){\includegraphics{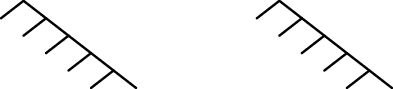}}
\put(0,13){$\scriptstyle2$}
\put(3.7,10){$\scriptstyle5$}
\put(7.5,7){$\scriptstyle6$}
\put(11.2,4){$\scriptstyle1$}
\put(15,1){$\scriptstyle3$}
\put(22.5,1){$\scriptstyle4$}
\put(43,13){$\scriptstyle1$}
\put(46.7,10){$\scriptstyle2$}
\put(50.5,7){$\scriptstyle3$}
\put(54.2,4){$\scriptstyle4$}
\put(58.5,1){$\scriptstyle5$}
\put(66,1){$\scriptstyle6$}
\end{picture}
\caption{\smaller The coloured
right vines $\vine{\{2, 5, 6\}, \{1, 3, 4\}}$
and $\vine{\{2, 5, 6, 1, 3, 4\}}$; the latter
is also $\vine{\{1, 2, 3, 4, 5, 6\}}$}
\label{F:DecVine}
\end{figure}

\begin{defi}
For $\II_1, \pp, \II_k$ finite subsets
of~$\NNN$, we define the {\it coloured right
vine} $\vine{\II_1,\!...,\!\II_k}$ by
$$\vine{\II_1, \pp, \II_k} =
\et_{\ell_1} \op ( \et_{\ell_2} \op ( \pp
 \op ( \et_{\ell_{n-1}} \op \et_{\ell_n}) \pp
)),$$ where $(\ell_1, \pp, \ell_n)$ is the
increasing enumeration of~$\II_1$, followed by
the increasing enumeration of~$\II_2$, etc. 
(Figure~\ref{F:DecVine}).
\end{defi}

In particular $\vine\II$ is the vine  in
which the labels are the elements
of~$\II$ enumerated in increasing order. By
construction, all coloured vines~$\vine{\II}$
are injective trees, so, clearly, we have:

\begin{lemm}
Coloured vines make a discriminating family
for the action of~$\GGAC$ on~$\Ti$.
\end{lemm}

The scheme is now the same as in
Section~\ref{S:PreF}: in order to apply
Proposition~\ref{P:Criterion}, we
select, for each tree~$\tt$ with labels~$\II$, a
word~$\HH\tt$ that describes how~$\tt$ can be
constructed from the vine~$\vine{\II}$. For the
skeleton, we can use associativity as
in Section~\ref{S:PreF}. For the labels, we use
commutativity, \ie, operators~$\CCp\a$.
The first step for an inductive construction is to
define an operator that maps
$\vine{\II \cup \JJ}$ to~$\vine{\II, \JJ}$ for
disjoint~$\II, \JJ$. To this end,
we introduce new elements of~$\GGAC$.

\begin{figure} [htb]
\begin{picture}(92, 21)(0, 0)
\put(0,0){\includegraphics{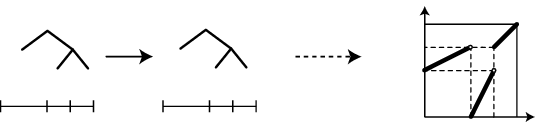}}
\put(20,14){$\SSp¥$}
\put(3,10){$\scriptstyle1$}
\put(9,7){$\scriptstyle2$}
\put(14,7){$\scriptstyle3$}
\put(29,10){$\scriptstyle2$}
\put(35.5,7){$\scriptstyle1$}
\put(40.5,7){$\scriptstyle3$}
\end{picture}
\caption{\smaller The action of~$\SSp¥$; the
operator~$\ss i$, \ie, $\SS{\1^{i-1}}$, switches
the $i$th and the $(i+1)$th factors of
$\vine{\tt_1, \pp, \tt_n}$, as it maps the latter
to
$\vine{\tt_1, \pp, \tt_{i+1}, \tt_i, \pp \tt_n}$.}
\label{F:ActS}
\end{figure}

\begin{defi}
(Figure~\ref{F:ActS})
For each address~$\a$, we put $\SS\a = \CC\a
\AA\a\inv \CC{\a\1}\inv$. We denote by~$\GS$
the family of all~$\SS\a$'s, and by~$\RACS$ the
family obtained by adding the definition
of~$\SS\a$ to~$\RAC$. For $i \ge 1$, we put $\ss
i = \SS{\1^{i-1}}$, \ie, $\ss i =  \cc i \aa i\inv
\cc{i+1}\inv$, and we denote by~$\Gs$ the
family of all~$\ss i$'s.
\end{defi}

\begin{defi}\label{D:PPSS}
For $\II, \JJ$ finite disjoint subsets of~$\NNN$, the
word~$\PP\II\JJ$ is inductively determined by
$\PP\varnothing\varnothing  = \e$ and the
rules: for~$\ell$ smaller than all
elements of~$\II$ and~$\JJ$, 
\begin{equation*}
\PP{\{\ell\} \cup \II}\JJ
=  \dd{\PP\II\JJ}, \quad
\PP\II{\{\ell\} \cup \JJ} =
\begin{cases}
\ss1 \ss2 \pp \ss{p-1} \cc p
&\text{if $\II$ has $p$~elements and $\JJ$ is
empty,}\\
\dd{\PP\II\JJ} \opp \ss1 \ss2 \pp \ss p
&\text{if $\II$ has $p$~elements and $\JJ$ is
nonempty.}
\end{cases}
\end{equation*}
The word~$\PPP\II\JJ$ is defined similarly, 
except that $\PP\II{\{\ell\}}$  is defined to be 
$\ss1 \ss2 \pp \ss p$. 
\end{defi}

\begin{exam}
Let $\II = \{2, 5, 6\}$ and $\JJ = \{1, 3, 4\}$. By
considering the elements of~$\II \cup \JJ$ in
decreasing order, we find successively
$\PP{\varnothing}{\varnothing}  = \e$, 
$\PP{\{6\}}{\varnothing}  = \e$, 
$\PP{\{5, 6\}}{\varnothing}  = \e$, 
$\PP{\{5, 6\}}{\{4\}}  = \ss1 \cc2$, \\
$\PP{\{5, 6\}}{\{3, 4\}} = \ddp{\ss1 \cc2}
\opp \ss1 \ss2 = \ss2 \cc3 \ss1 \ss2$, 
$\PP{\{2, 5, 6\}}{\{3, 4\}} =
\ddp{\ss2 \cc3 \ss1 \ss2} =
\ss3 \cc4 \ss2 \ss3$, \\
$\PP{\{2, 5, 6\}}{\{1, 3, 4\}} =
\ddp{\ss3 \cc4 \ss2 \ss3} \opp \ss1 \ss2 \ss3 
= \ss4 \cc5 \ss3 \ss4 \ss1 \ss2 \ss3$, and
$\PPP{\{2, 5, 6\}}{\{1, 3,
4\}} = \ss4 \ss5 \ss3 \ss4 \ss1 \ss2 \ss3$.
\end{exam}

\begin{lemm}\label{L:ActP}
For all sets~$\II, \JJ$, and for every tree~$\tt$, we have
$$\begin{CD}
\vine{\II \cup \JJ} 
@>\PP\II\JJ>>
\vine{\II, \JJ}
\end{CD}
\text{\qquad and \qquad}
\begin{CD}
\vine{\II \cup \JJ, \tt} 
@>\PPP\II\JJ>>
\vine{\II, \JJ, \tt}
\end{CD}.$$
\end{lemm}

\begin{proof}
We use induction on the cardinality of~$\II \cup
\JJ$. The result is clear if $\II$ and $\JJ$ are
empty. Assume that $\ell$ is smaller than
all elements in~$\II$ and~$\JJ$. The induction
hypothesis asserts that $\PP\II\JJ$ maps
$\vine{\II \cup \JJ}$ to $\vine{\II,\JJ}$, hence
$\dd{\PP\II\JJ}$ maps $\et_\ell \op \vine{\II\cup\JJ}$,
which is $\vine{\{\ell\} \cup \II \cup \JJ }$, to
$\et_\ell \op \vine{\II, \JJ}$, \ie, to
$\vine{\{\ell\} \cup \II, \JJ}$, as expected
for~$\PP{\{\ell\} \cup \II}\JJ$.

Let us consider $\PP\II {\{\ell\} \cup \JJ}$. Let $p$
be the cardinal of~$\II$. Assume first $\JJ \not=
\emptyset$.  We have seen that 
$\dd{\PP\II\JJ}$ maps $\vine{\{\ell\} \cup \II
\cup \JJ }$ to
$\vine{\{\ell\} \cup \II, \JJ}$. Then the iterated
transposition~$\ss1 \ss2 \pp \ss p$ carries the
leftmost leaf of~$\vine{\{\ell\} \cup \II, \JJ}$,
\ie,~$\et_\ell$, through $p$~leaves to the
right, \ie, we obtain $\vine{\II, \{\ell\}, \JJ}$,
which is also $\vine{\II, \{\ell\} \cup \JJ}$. 
Finally, if $\JJ$ is empty, then $\PP\II\JJ$
is~$\e$, as an induction shows, and $\ss1 \ss2 \pp
\ss{p-1} \cc p$ maps $\vine{\{\ell\}, \II}$ to
$\vine{\II, \{\ell\}}$. So, in each case,
$\PP{\II}{\{\ell\} \cup \JJ}$ maps $\vine{\{\ell\}
\cup \II \cup \JJ }$ to $\vine{\II, \{\ell\} \cup
\JJ}$.

The argument is similar for~$\PPP\II\JJ$.
\end{proof}

We are now ready for defining a word~$\HH\tt$ that
describes how to construct a coloured tree~$\tt$
with labels~$\II$ from the right
vine~$\vine{\II}$. The current construction is
similar to that of Section~\ref{S:PreF}. The
only change is that, in the induction step, we
first sort the labels in order to push to the
initial positions the labels that correspond to
the left subtree. This is exactly  what (the
operators associated with) $\PP\II\JJ$
and~$\PPP\II\JJ$ do. So the following definition
should be natural. 

\begin{defi}
For each injective tree~$\tt$, the words~$\HH\tt,
\HHH\tt$ are defined by the rules: 
\begin{align*}
\HH\tt 
&= \HHH\tt = \e 
&&\text{for $\tt$ of size~$1$,}\\
\HH\tt 
= \PP{\IIi}{\IIii} \opp \HHH\tti \opp
\dd{\HH\ttii},
&\quad \HHH\tt 
= \PPP{\IIi}{\IIii} \opp \HHH\tti \opp
\dd{\HHH\ttii} \opp \AA¥ 
&&\text{for $\tt = \tti \op \ttii$ and $\II_k$ the
labels in~$\tt_k$.}
\end{align*}
\end{defi}

The following result is the exact counterpart to
Lemma~\ref{L:ConF}. For $\II = \{\ell_1, \pp,
\ell_n\}$ and $\tt$ a tree, we use
$\vine{\II, \tt}$ for $\vine{\et_{\ell_1}, \pp,
\et_{\ell_n}, \tt}$.

\begin{lemm} \label{L:ConV}
For each injective tree~$\tt$ with labels~$\II$,
and each tree~$\ttt$, we have
\begin{equation}
\begin{CD}
\vine\II
@>\HH\tt>>
\tt
\end{CD}
\text{\qquad and \qquad}
\begin{CD}
\vine{\II, \ttt}
@>\HHH\tt>>
\vine{\tt, \ttt}
\end{CD},
\end{equation}
\ie, $\HH\tt$ constructs~$\tt$ from~$\vine\II$,
and $\HHH\tt$ constructs~$\vine{\tt, \ttt}$
from~$\vine{\II, \ttt}$.
\end{lemm}

\begin{proof}
The inductive verification is the same as for
Lemma~\ref{L:ConF}. The diagrams are now:
\begin{gather*}
\begin{CD}
\vine \II = \vine{\IIi \cup \IIii}
@>\PP\IIi\IIii>>
\vine{\IIi, \IIii}
@>\HHH\tti>>
\vine{\tti, \IIii}
@>\dd{\HH\ttii}>>
\vine{\tti, \ttii} = \tt,
\end{CD}\\
\begin{CD}
\vine{\II, \ttt} = \vine{\IIi \cup \IIii, \ttt}
@>\PPP\IIi\IIii>>
\vine{\IIi, \IIii, \ttt}
@>\HHH\tti>>
\vine{\tti, \IIii, \ttt}
\end{CD}\hspace{3cm}\\
\hspace{5cm}\begin{CD}
@>\dd{\HHH\ttii}>>
\vine{\tti, \ttii, \ttt} 
@>{\AA¥}>>
\vine{\tti \op \ttii, \ttt}
= \vine{\tt, \ttt}
\end{CD}
\end{gather*}
for $\tt = \tti \op \ttii$ and $\IIi, \IIii$ the sets
of labels in~$\tti$ and~$\ttii$
respectively.
\end{proof}

\subsection{Derived relations}

In order to apply Proposition~\ref{P:Criterion} and
prove that the relations of
Proposition~\ref{P:RelV} make a presentation of
the group~$\GGAC$, \ie,  of~$V$, we have
to check that there are enough relations to
establish the equivalence of~$\HH\ttt$ and
$\HH\tt \opp \BB\a$ whenever $\BB\a$ maps~$\tt$
to~$\ttt$, where $\BB¥$ is either~$\AA¥$
or~$\CC¥$. The needed verifications are easy,
but longer than in the case of~$\GGA$, and we 
begin with some technical, but easy
preparatory results asserting that certain
relations involving the letters~$\AA\a$,
$\CC\a$, and~$\SS\a$ follow from~$\RACS$.

\begin{lemm} The following relations follow
from~$\RACS$:

(i) The  $\AA¥$- and~$\CC¥$-geometric relations
of~$\RAC$ in which $\BB¥$ or~$\BBB¥$ is
replaced with~$\SS¥$;

(ii) The $\SS¥$-geometric relations, defined to be
the translated copies of
\begin{equation}
\tag{$\RGS$}
\BB{\1\1\a} \opp \SS¥  = \SS¥ \opp \BB{\1\1\a},
\qquad
\BB{\1\0\a} \opp \SS¥ = \SS¥ \opp \BB{\0\a},
\qquad
\BB{\0\a} \opp \SS¥ = \SS¥ \opp \BB{\1\0\a},
\end{equation}
in which $\BB¥$ stands for~$\AA¥, \CC¥$
or~$\SS¥$, 

(iii) The translated copies of the relations
\begin{gather}
\label{E:RACS1}
\SS¥ \AA¥ =  \AA¥ \CC\0,
\text{\qquad}
\SS¥ \AA\1 \AA¥ 
= \AA\1 \AA¥ \SS\0,\\
\label{E:RACS2}
\SS\1 \SS¥ \AA\1 = \AA¥ \SS¥, 
\quad
\SS¥ \SS\1 \AA¥ = \AA\1 \SS¥,
\quad 
\SS¥ \SS\1 \SS¥ = \SS\1 \SS¥ \SS\1.
\end{gather}
\end{lemm}

\begin{proof}
The extension of the $\AA¥$- and
$\CC¥$-geometric relations to~$\SS\a$ is
obvious, as $\SS\a$ is defined from~$\CC\a$,
$\AA\a$, and~$\CC{\1\a}$.  The
$\SS¥$-geo\-metric relations follow from the
other geometric relations. For instance, we find
\begin{align*}
\SS¥ \BB{\1\1\a} 
= \AA¥ \CC\0 \AA¥\inv \BB{\1\1\a}
&\eq{\RGA}
\AA¥ \CC\0  \BB{\1\a} \AA¥\inv
\eq{\RGperp} 
\AA¥ \BB{\1\a} \CC\0   \AA¥\inv
\eq{\RGA}  
\BB{\1\1\a} \AA¥  \CC\0   \AA¥\inv
=\BB{\1\1\a} \SS¥.
\end{align*}
The first relation in~\eqref{E:RACS1}
follows from the definition and an hexagon
relation:
$$\SS¥ \AA\0 = \CC¥ \AA¥\inv \CC\1\inv \AA¥
\eqH \AA¥ \CC\0.$$
The second relation
comes by cancelling~$\AA\0$ on the right in
\begin{equation*}
\SS¥ \AA\1 \AA¥ \AA\0
\eqP \SS¥ \AA¥ \AA¥
\eq{\eqref{E:RACS1}} \AA¥ \CC\0 \AA¥ 
\eqG \AA¥ \AA¥ \CC{\0\0}
\eqP \AA\1 \AA¥ \AA\0 \CC{\0\0}
\eq{\eqref{E:RACS1}} \AA\1 \AA¥ \SS\0 \AA\0.
\end{equation*}
Then we observe that the hexagon relation implies
\begin{equation} \label{E:Pre1}
\CC\1 \SS¥
\eq¥
\CC\1 \SS¥ \AA¥ \AA¥\inv
\eq{\eqref{E:RACS1}}
\CC\1 \AA¥ \CC\0 \AA¥\inv
\eqH
\AA¥ \CC¥ \AA¥ \AA¥\inv
\eq¥ \AA¥ \CC¥.
\end{equation}
Next, the first two relations in~\eqref{E:RACS2}
are obtained by  cancelling~$\AA¥$ on the right
in
\begin{align*}
\SS\1 \SS¥ \AA\1 \AA¥
&\eq{\eqref{E:RACS1}} 
\SS\1 \AA1 \AA¥ \SS\0 
\eq{\eqref{E:RACS1}} 
\AA\1 \CC{\1\0} \AA¥ \SS\0 \\
&\eqG 
\AA\1¥ \AA¥ \CC{\0\1} \SS\0 
\eq{\eqref{E:Pre1}} 
\AA\1¥ \AA¥ \AA\0 \CC\0 
\eqP 
\AA¥ \AA¥ \CC\0 
\eq{\eqref{E:RACS1}} 
\AA¥ \SS¥ \AA¥,\\
\SS¥ \SS\1 \AA¥ \AA¥
&\eqH \SS¥ \SS\1 \AA\1 \AA¥ \AA\0
\eq{\eqref{E:RACS1}} 
\SS¥ \AA\1 \CC{\1\0} \AA¥ \AA\0
\eqG \SS¥ \AA\1 \AA¥ \CC{\0\1}  \AA\0\\
&\eq{\eqref{E:RACS1}} \AA\1 \AA¥ \SS\0 \CC{\0\1}  \AA\0
\eq¥ 
\AA\1 \AA¥ \CC\0
\eq{\eqref{E:RACS1}} 
\AA\1 \SS¥ \AA¥.
\end{align*}
Finally, we have
$$\SS¥ \CC\1 \SS¥
\eq{\eqref{E:Pre1}} 
\SS¥ \AA¥ \CC¥
\eq{\eqref{E:RACS1}} 
\AA¥ \CC\0 \CC¥
\eqG \AA¥ \CC¥ \CC\1
\eq{\eqref{E:Pre1}} 
\CC\1 \SS¥ \CC\1,$$
so, using $\SS¥ \AA\1 \AA¥
\eq{\eqref{E:RACS1}} 
\AA\1 \AA¥ \SS\0$, and 
$\SS\1 \AA\1 \AA¥
\eq{\eqref{E:RACS1}} 
\AA\1 \CC{\1\0} \AA¥
\eqG 
\AA\1 \AA¥ \CC{\0\1}$, we deduce
$$ \AA\1 \AA¥ \SS¥ \SS\1 \SS¥
\eq¥
\DD\0{(\SS¥ \CC\1 \SS¥)} \opp  \AA\1 \AA¥
\eq¥
\DD\0{(\CC\1 \SS¥ \CC\1)} \opp \AA\1 \AA¥
\eq¥
\AA\1 \AA¥ \SS\1 \SS¥ \SS\1,$$
which implies the third relation
in~\eqref{E:RACS2} by cancelling~$\AA\1
\AA¥$ on the left.
\end{proof}

On the other hand, we observe that, by
construction, the words~$\HH\tt$
and~$\HHH\tt$ involve the letters~$\aa i$, $\cc
i$, and~$\ss i$ only. So it will be convenient to
work with the following restricted list.

\begin{defi}
We define~$\Racs$ to consist of  the
following relations:
\begin{gather}
\label{E:Racs1}
\aa i \bb{j-1}
= \bb j \aa i,
\text{\qquad with $j \ge i + 2$ and $\bb¥ =
\aa¥, \cc¥$ or~$\ss¥$},\\
\label{E:Racs2}
\ss i \bb j 
= \bb j \ss i,
\text{\qquad with $j \ge i + 2$ and 
$\bb¥ = \aa¥,  \cc¥$ or~$\ss¥$},\\
\label{E:Racs3}
\ss i \ss{i+1} \aa i 
= \aa{i+1} \ss i,
\text{\qquad and\qquad }
\ss{i+1} \ss i \aa{i+1}
= \aa i \ss i,\\
\label{E:Racs4}
\ss i \bb{i+1} \ss i 
= \bb{i+1} \ss i \bb{i+1},
\text{\qquad with $\bb¥ = \ss¥$ 
or~$\cc¥$}.
\end{gather}
\end{defi}

\begin{lemm} \label{L:Coxe}
All relations in~$\Racs$ are consequences
of~$\Rac$, hence of~$\RAC$ (plus the
definitions of~$\aa i, \cc i$ and~$\ss i$).
Furthermore, $\ss i^2 = 1$ follows from~$\Rac$
completed with the relations~$\cc i^2 = 1$.
\end{lemm}

\begin{proof}
When $\bb¥$ is~$\aa¥$ or~$\cc¥$,
\eqref{E:Racs1} coincides
with~\eqref{E:Rac1}; for $\bb j = \ss j$, we apply 
\eqref{E:Rac1} to~$\cc j$, $\aa j\inv$, and
$\cc{j+1}\inv$ successively. Similarly,
\eqref{E:Racs2} for $\bb j = \aa j$
or~$\cc j$ directly follows
from~\eqref{E:Rac2} owing to the
definition of~$\ss i$; the relation for
$\bb j = \ss j$ then follows by
replacing~$\ss j$ with its definition.
As for~\eqref{E:Racs3}, we find
\begin{align*}
\ss i \ss{i+1} \aa i
\eq¥ \cc i \aa i\inv \aa{i+1}\inv
\cc{i+2}\inv \aa i
\eq{\eqref{E:Rac1}}
\cc i \aa i\inv \aa{i+1}\inv
\aa i \cc{i+1}\inv 
\eq{\eqref{E:Rac3}}
\aa{i+1} \cc i \aa i\inv \cc{i+1}\inv
\eq{\eqref{E:Racs5}}
\aa{i+1} \ss i.
\end{align*}
Next, we observe that the relation~\eqref{E:Rac4}
of~$\Rac$ implies
\begin{equation} \label{E:Racs5}
\ss i = \cc i \aa i\inv \cc{i+1}\inv
\eq\Rac \cc{i+1}\inv \aa i \cc i,
\end{equation}
and we deduce symmetrically
\begin{align*}
\ss{i+1} \ss i \aa{i+1}
\eq¥ \cc{i+2}\inv \aa{i+1} \aa i \cc i 
\aa{i+1}
\eq{\eqref{E:Rac3}}
\cc{i+2}\inv \aa i^2 \cc i  \aa{i+1}
\eq{\eqref{E:Rac1}}
 \aa i \cc{i+1}\inv \aa i \cc i  \aa{i+1}
= \aa i \ss i.
\end{align*}
For~\eqref{E:Racs4} with $\bb 2 = \cc 2$, we have
$$\ss1 \cc2 \ss1
\eq¥ 
\cc1 \aa1\inv \ss1
= \cc1 \aa1\inv \cc1 \aa1\inv \cc2\inv 
\eq{\eqref{E:Rac5}}
\cc2 \cc1 \aa1\inv 
\eq¥ 
\cc2 \ss1 \cc2 .$$
As for~\eqref{E:Racs4} with $\bb 2 = \ss 2$,
we have
\begin{align*}
\ss1 \ss2 \ss1
= \ss1 \cc2 \aa2\inv \cc3\inv \ss1
&\eq{\eqref{E:Racs2}}
\ss1 \cc2 \aa2\inv \ss1 \cc3\inv 
\eq{\eqref{E:Racs3}}
\ss1 \cc2 \ss1 \aa1\inv \ss2\inv \cc3\inv,\\
\ss2 \ss1 \ss2 
= \cc2 \aa2\inv \cc3\inv \ss1 \ss2
&\eq{\eqref{E:Racs2}}
\cc2 \aa2\inv \ss1 \cc3\inv \ss2
\eq{\eqref{E:Racs3}}
\cc2 \ss1 \aa1\inv  \ss2\inv \cc3\inv \ss2 \\
& \eq¥
\cc2 \ss1 \cc2 \cc2\inv \aa1\inv  \ss2\inv
\cc3\inv \ss2
\eq{\eqref{E:Racs1}}
\cc2 \ss1 \cc2 \aa1\inv  \cc3\inv  \ss2\inv
\cc3\inv \ss2.
\end{align*}
Applying the relations $\ss1 \cc2 \ss1 \eq¥ \cc2
\ss1 \cc2$ and $\ss2 \cc3 \ss2 \eq¥ \cc3 \ss2
\cc3$---hence $\ss2\inv \cc3\inv \eq¥ \cc3\inv  \ss2\inv
\cc3\inv \ss2$---which were established above, we
deduce $\ss1 \ss2 \ss1 \eq¥ \ss2 \ss1
\ss2$.

Finally, we have seen that $\Rac$ implies
$\ss1 = \cc1 \aa1\inv \cc2\inv \eq¥ \cc2\inv
\aa1 \cc1$, hence $\ss1^2 \eq¥ \cc1 \cc2^{-2}
\cc1$: so $\cc1^2 \eq¥  \cc2^2 \eq¥  1$ implies
$\ss1^2 \eq¥  1$.
\end{proof}

For furure inductive arguments, we need
some results about the auxiliary
words~$\PP\II\JJ$ and~$\PPP\II\JJ$.

\begin{lemm} \label{L:PPPP}
For $\II, \JJ, \KK$ disjoint with $p = \card\II
\ge 1$, we have
\begin{equation} 
\label{E:PPPP}
\bb{\II \cup \JJ, \KK} \opp \PPP\II\JJ
\eqacs \bb{\II, \JJ \cup\KK}  \opp
\ddd{p}{\bb{\JJ,\KK}} 
\text{\qquad for $\bb¥ = \ss¥$ and $\bb¥ =
\cc¥$.}
\end{equation}
\end{lemm}

\begin{proof}
We begin with the auxiliary formulas
\begin{gather}
\label{E:Cox3*}
\ss1 \ss2 \pp \cc{k+1} \ss i
\eqacs \cc{i+1} \ss1 \ss2 \pp \cc{k+1}
\mbox{,\quad for $1 \le i \le k$},\\
\label{E:Cox4*}
\ss1 \ss2 \pp \ss k \cc{k+1} \ss k
\eqacs \cc{k+1} \ss1 \ss2 \pp \ss k \cc{k+1}
\mbox{,\quad for $1 \le k$}.
\end{gather}
A direct inductive verification is possible; we can also 
observe that Lemma~\ref{L:Coxe} shows that
$(\ss1, \pp, \ss{k+1})$ and $(\ss1, \pp, \ss{k},
\cc{k+1})$ satisfy the relations of Artin's
presentation of the braid group~$B_{k+2}$:
therefore, every braid relation between the
standard generators~$\sigma_i$ of~$B_{k+2}$
must hold between the~$\ss i$'s, which is the
case for the counterpart of~\eqref{E:Cox3*}
and~\eqref{E:Cox4*}.

Next, we claim that the following relations are true,
where $q$ denotes~$\card\JJ$:
\begin{gather}
\label{E:Cox5}
\ss1 \ss2 \pp \ss{q + r} \PPP\JJ\KK
\eqacs \dd{\PPP\JJ\KK} \opp \ss1 \ss2
\pp \ss{q + r},\\
\label{E:Cox5*}
\ss1 \ss2 \pp \ss{q + r-1} \cc{q+r}
\PPP\JJ\KK
\eqacs \dd{\PP\JJ\KK}
\ss1 \ss2 \pp \ss{q + r-1} \cc{q+r}.
\end{gather}
Indeed, an easy induction shows that the
word~$\PP\JJ\KK$ is a product of~$\ss i$'s
with $1 \le i \le q + r -2$, and, if it not empty,
of~$\cc{q+r-1}$ occurring once, and that
$\PPP\JJ\KK$ is obtained from~$\PP\II\JJ$ by
replacing the possible $\cc{q+r-1}$
with~$\ss{q+r-1}$.
Then \eqref{E:Cox5} comes by
applying~\eqref{E:Cox3*} with $k = q + r - 1$ to the
letters~$\ss i$ in~$\PP\JJ\KK$, and so
does \eqref{E:Cox5*} using \eqref{E:Cox4*} for the
possible letter~$\cc{q+r-1}$ of~$\PPP\JJ\KK$.

We turn to the first formula
in~\eqref{E:PPPP}. The result is trivial for $\II = \JJ = \KK
= \emptyset$. For an induction, it is sufficient to
prove that, if $\ell$ is smaller than
all elements in $\II \cup \JJ \cup \KK$, the result is
true for $(\{\ell\}, \JJ, \KK)$, and it is true for
$(\{\ell\} \cup \II, \JJ, \KK)$, $(\II, \{\ell\} \cup \JJ,
\KK)$, and $(\II, \JJ, \{\ell\} \cup \KK)$ whenever it is for
$(\II, \JJ, \KK)$ and
$\II$ is nonempty.
For the case of $(\{\ell\}, \JJ, \KK)$, we find
\begin{align*}
\PP{\{\ell\}}{\JJ \cup \KK} \opp \PPP\JJ\KK
&= \ss1 \pp \ss{q+r-1} \cc{q+r} \opp \PPP\JJ\KK \\
& \eq{\eqref{E:Cox5*}}
\dd{\PP\JJ\KK} \opp \ss1 \pp \ss{q+r-1}
\cc{q+r}\\
&= \dd{\PP\JJ\KK} \opp \ss1 \pp \ss r
\opp \dddp{r}{\ss1 \pp \ss{q-1} \cc q}
= \PP{\{\ell\} \cup \JJ}\KK 
\opp \ddd{r}{\PP{\{\ell\}}\JJ}.
\end{align*}
Assume $\II \not= \emptyset$. For $(\{\ell\} \cup \II, \JJ,
\KK)$, using the induction hypothesis, we find
\begin{align*}
\PP{\II\cup \{\ell\}}{\JJ  \cup \KK}
\PPP\JJ\KK
&= \dd{\PP\II{\JJ \cup \KK}} \opp \ss1 \ss2
\pp \ss{q + r} \opp \PPP\JJ\KK
\eq{\eqref{E:Cox5}} \dd{\PP\II{\JJ \cup \KK}}
\opp 
\dd{\PPP\JJ\KK} \opp \ss1 \ss2
\pp \ss{q + r} \\
&\eq{(IH)} \dd{\PPP{\II\cup\JJ}\KK}
\opp  \ddd{{r+1}}{\PP\II\JJ} \opp \ss1 \ss2
\pp \ss{q + r} \\
&\eq{\eqref{E:Racs2}} \dd{\PPP{\II\cup\JJ}\KK}
\opp  \ss1 \ss2 \pp \ss r \opp
\ddd{r+1}{\PP\II\JJ} \opp \ss{r+1}
\pp \ss{q + r} \\
&= \dd{\PPP{\II\cup\JJ}\KK}
\opp  \ss1 \ss2 \pp \ss r \opp
\dddp{r}{\dd{\PP\II\JJ} \opp \ss1 \ss2
\pp \ss q} 
= \PP{\{\ell\} \cup \II \cup \JJ}\KK \opp
\ddd{r}{\PP{\II\cup \{\ell\}}\JJ},
\end{align*}
The remaining cases are easy:
\begin{align*}
\PP\II{\{\ell\} \cup \JJ \cup \KK}
\opp \PPP{\{\ell\} \cup \JJ}\KK
&= \dd{\PP\II{\JJ \cup \KK}} \opp
\dd{\PPP\JJ\KK} \opp \ss1 \ss2 \pp
\ss r 
\eq{(IH)} \dd{\PPP{\II\cup\JJ}\KK} \opp
\ddd{{r+1}}{\PP\II\JJ} \opp  \ss1 \ss2 \pp
\ss r \\
&\eq{\eqref{E:Racs2}} \dd{\PPP{\II\cup\JJ}\KK} 
\opp \ss1 \ss2 \pp \ss r \opp
\ddd{{r+1}}{\PP\II\JJ}
= \PP{\II \cup \{\ell\} \cup \JJ}\KK \opp
\ddd{r}{\PP\II{\{\ell\} \cup \JJ}};\\
\PP\II{\JJ \cup \{\ell\} \cup \KK}
\PPP\JJ{\{\ell\} \cup \KK}
&= \dd{\PP\II{\JJ \cup \KK}} \opp
\dd{\PPP\JJ\KK}
\eq{(IH)} \dd{\PP{\II \cup \JJ}\KK \opp
\ddd{r}{\PP\II\JJ}}
= \PP{\II \cup \JJ}{\{\ell\} \cup \KK} \opp
\ddd{{r+1}}{\PP\II\JJ}, 
\end{align*}
and the proof is complete.
\end{proof}

\begin{defi}
For $p, q \ge 1$, we put 
$\PP p q= \PP{\{q+1, \pp, q+p\}}{\{1, \pp, q\}}$
and $\PPP p q= \PPP{\{q+1, \pp, q+p\}}{\{1, \pp,
q\}}$.
\end{defi}   

So $\PPP p q$ is the iterated transposition that switches
two blocks of~$p$ and $q$~elements respectively, putting
the $p$~elements on the top. For instance, we have
$\PPP p 1 = \ss1 \pp \ss{p-1}$, and $\PPP 1 q = \ss{q-1}
\pp \ss1$.

\begin{lemm} 
For all~$p, q, r$, we have
\begin{gather}
\label{E:PPP1}
\PP{p+q}r 
\eqacs \PPP p r \opp \ddd{p}{\PP q r}
\text{\quad and \quad}
\PPP{p+q}r 
\eqacs \PPP p r \opp \ddd{p}{\PPP q r},\\
\label{E:PPP2}
\PP p{q+r} 
\eqacs \ddd{q}{\PP p r} \opp \PPP p q
\text{\quad and \quad}
\PPP p{q+r} 
\eqacs \ddd{q}{\PPP p r} \opp \PPP p q,\\
\label{E:PPP4}
\aa{q+1} \opp \PPP{p+1}q
\eqacs \PPP{p+2}q \opp \aa1.
\end{gather}
\end{lemm}

\begin{proof}
Relation~\eqref{E:PPP1} and~\eqref{E:PPP2} follow
from~\eqref{E:PPPP} by taking
$\II = \{r+1, \pp, r+p\}$,
$\JJ = \{r+p+1, \pp, r+p+q\}$,
$\KK = \{1, \pp, r\}$, and
$\II = \{q+r+1, \pp, q+r+p\}$,
$\JJ = \{1, \pp, q\}$, 
$\KK = \{q+1, \pp, q+r\}$,
respectively. In the first case, we have $\PP\II{\JJ \cup
\KK} = \PPP p r$, and, in the second one, we have $\PP{\II
\cup \JJ}\KK = \ddd{q}{\PP p r}$.
For~\eqref{E:PPP4}, we use induction.
For $q = 0$, the result is clear; for~$q \ge
1$, we find
\begin{align*}
 \aa{q+1} \PPP{p+1}q 
&\eq{\eqref{E:PPP1}} 
\aa{q+1} \opp \dd{\PPP{p+1}{q-1}} \opp \PPP{p+1}1
= \ddp{\aa q \PPP{p-1}{q-1}} \opp \PPP{p+1}1\\
&\eq{(IH)} \dd{\PPP{p+2}{q-1} \aa1} \opp \PPP{p+1}1
= \dd{\PPP{p+2}{q-1}} \opp \aa2 \ss1 \pp \ss{p+1}\\
&\eq{\eqref{E:Racs3}} \dd{\PPP{p+2}{q-1}} \opp 
\ss1 \ss2 \aa1 \ss2 \pp \ss{p+1}\\
&\eqacs \dd{\PPP{p+2}{q-1}} \opp 
\ss1 \ss2 \ss3 \pp \ss{p+2} \aa1
\eq{\eqref{E:PPP1}}  \PPP{p+2}q \aa1,
\end{align*}
which completes the computation.
\end{proof}

\begin{lemm}
Assume that $\tt$ is a size~$n$ tree. Then, for~$p, q \ge
0$, we have
\begin{gather}
\label{E:ComW}
\bb{i+n} \opp \HHH\tt
\eqacs \HHH\tt \opp \bb{i+1},
\text{\qquad for~$\bb¥ = \aa¥, \cc¥$
or~$\ss¥$,}\\
\label{E:HHP1}
\ddd{q}{\HH\tt} \opp \PP1q 
\eqacs \PP nq \opp \HHH\tt
\text{ \quad and \quad}
\ddd{q}{\HHH\tt} \opp \PPP{p+1}q
\eqacs\PPP{p+n}q \opp \HHH\tt,\\
\label{E:HHP2}
\HHH\tt \opp \PP p{q+1}
\eqacs \PP p{q+n} \opp \ddd{q}{\HH\tt},
\text{ \quad and \quad}
\HHH\tt \opp \PPP p{q+1}
\eqacs \PPP p{q+n} \opp \ddd{q}{\HHH\tt}.
\end{gather}
\end{lemm}

\begin{proof}
We use induction on~$n$. For $n = 1$, the
words~$\HH\tt$ and~$\HHH\tt$ are empty, and all
relations are equalities. Otherwise, assume $\tt =
\tti
\op
\ttii$, with, as usual, $\nn_k$ the size
of~$\tt_k$ and $\II_k$ its set of labels.
For~\eqref{E:ComW}, we find
\begin{align*}
\bb{i+n} \opp \HHH\tt 
&= \bb{i+n} \opp \PPP\IIi\IIii \opp \HHH\tti \opp
\dd{\HHH\ttii} \opp \AA¥ 
\eqG \PPP\IIi\IIii \opp \bb{i+n} \opp
\HHH\tti \opp
\dd{\HHH\ttii} \opp \AA¥ \\
&\eqIH \PPP\IIi\IIii \opp \HHH\tti \opp
\bb{i+\nnii}
\opp \dd{\HHH\ttii} \opp \AA¥ \\
&\eqIH \PPP\IIi\IIii \opp \HHH\tti \opp \dd{\HHH\ttii} 
\opp \bb{i+2} \opp \AA¥
\eqG \PPP\IIi\IIii \opp \HHH\tti \opp
\dd{\HHH\ttii} 
\opp \AA¥ \opp \bb{i+1}
= \HHH\tt \opp \bb{i+1}.
\end{align*}
(The first equivalence holds because we
consider~$\PPP\IIi\IIii$, which consists of~$\ss i$'s only.)

We turn to the second relation
in~\eqref{E:HHP1}. Then the expected
relation follows from the commutativity of the
following diagram
$$\begin{CD}
\vine{q,\! \II\!,\! p,\! \ttt}
@>{\ddd{q}{\PPP\IIi\IIii}}>>
\vine{q,\! \IIi\!,\! \IIii\!,\! p,\! \ttt}
@>{\ddd{q}{\HHH\tti}}>>
\vine{q,\! \tti\!,\! \IIii\!,\! p,\! \ttt}
@>{\ddd{q+1}{\HHH\ttii}}>>
\vine{q,\! \tti\!,\! \IIii\!,\! p,\! \ttt}
@>{\aa{q+1}}>>
\vine{q,\! \tti \op \ttii\!,\! p,\! \ttt}\\
@VV{\PPP{p+\nn}q}V
@VV{\PPP{p+\nn}q}V
@VV{\PPP{p+\nnii+1}q}V
@VV{\PPP{p+2}q}V
@VV{\PPP{p+1}q}V\\
\vine{\II\!,\! p, \!q,\! \ttt}
@>{\PPP\IIi\IIii}>>
\vine{\IIi\!,\! \IIii\!,\! p, \!q,\! \ttt}
@>{\HHH\tti}>>
\vine{\tti\!,\! \IIii\!,\! p, \!q,\! \ttt}
@>{\dd{\HHH\ttii}}>>
\vine{\tti\!,\! \IIii\!,\! p, \!q,\! \ttt}
@>{\aa1}>>
\vine{\tti \op \ttii\!,\! p, \!q,\! \ttt}
\end{CD}$$
The first (leftmost) square is commutative
by~\eqref{E:PPP1}. The second one is commutative by
induction hypothesis. For the third, \eqref{E:PPP1}
tells us that $\PPP{p+ \nnii + 1}q$ is $\Racs$-equivalent to
$\PPP1q \opp \dd{\PPP{p+\nnii}q}$, and that
$\PPP{p+2}q$ is $\Racs$-equivalent to
$\PPP1q \opp \dd{\PPP{p+1}q}$. As
$\ddd{q+1}{\HHH\ttii}$ $\Racs$-commutes with $\PPP1q$
by geometric relations, we are left with proving the
$\Racs$-equivalence of
$\ddd{q}{\HHH\ttii} \opp \PPP{p+1}q$ and
$\PPP{p\nnii}q \opp \dd{\HHH\ttii}$, which is the
induction hypothesis. Finally, the commutativity of the last
square follows from~\eqref{E:PPP4}.

The verification of the other three formulas is similar.
\end{proof}

We are now in position for proving the counterpart
to Lemma~\ref{L:BluF}: 

\begin{lemm}\label{L:BluV} 
Assume $\ttt = \tt \act \BB\a$, where
$\BB¥$ is $\AA¥, \CC¥$, or~$\SS¥$. Then we have
\begin{equation} \label{E:BluV}
\HH\ttt \eqACS \HH\tt \opp \BB\a
\mbox{\quad and \quad}
\HHH\ttt \eqACS \HHH\tt \opp \BB{\0\a}.
\end{equation}
\end{lemm}

\begin{proof}
Clearly, it suffices to consider the cases of~$\AA\a$
and~$\CC\a$, as $\SS\a$ is defined from the latter. As for
Lemma~\ref{L:BluF}, we use induction on the length
of~$\a$ as a sequence of~$0$'s and~$1$'s. So assume
first that $\a$ is the empty address. Let us
consider the case of~$\AA¥$. The hypothesis
$\ttt = \tt \act \AA¥$ implies that there  exist
trees~$\tti, \ttii, \ttiii$ such that $\tt$ is
$(\tti \op \ttii) \op \ttiii$, and $\ttt$ is $\tti \op
(\ttii \op \ttiii)$. We write $\IIi$ (\resp $\IIii$,
$\IIiii$) for the labels in~$\tti$ (\resp
$\ttii, \ttiii$), and $\nni$ (\resp $\nnii,
\nniii$) for their size. We obtain
\begin{gather}
\label{E:BluV1}
\HH\ttt 
= \PP{\IIi \cup \IIii}\IIiii \opp \PPP\IIi\IIii \opp
\HHH\tti \opp \dd{\HHH\ttii} \opp \AA¥ \opp
\dd{\HH\ttiii}\\
\label{E:BluV2}
\HH\tt \opp \AA¥
= \PPP\IIi{\IIii \cup \IIiii} \opp \HHH\tti \opp
\dd{\PP\IIii\IIiii} \opp \dd{\HHH\ttii} \opp
\ddd 2{\HH\ttiii} \opp \AA¥\\
\label{E:BluV3}
\HHH\ttt 
= \PPP{\IIi \cup \IIii}\IIiii \opp \PPP\IIi\IIii \opp
\HHH\tti \opp \dd{\HHH\ttii} \opp \AA¥ \opp
\dd{\HHH\ttiii} \opp \AA¥ \\
\label{E:BluV4}
\HHH\tt \opp \AA\0
= \PPP\IIi{\IIii \cup \IIiii} \opp \HHH\tti \opp
\dd{\PPP\IIii\IIiii} \opp \dd{\HHH\ttii} \opp
\ddd 2{\HHH\ttiii} \opp \AA\1 \AA¥ \AA\0.
\end{gather}
Using geometric relations, we may move the
factor~$\AA¥$ to the right in~\eqref{E:BluV1}, while, 
 in~\eqref{E:BluV2}, we may replace $\HHH\tti \opp
\dd{\PP\IIii\IIiii}$ with $\ddd p{\PP\IIii\IIiii} \opp
\HHH\tti$ using \eqref{E:ComW}. Then,
applying~\eqref{E:PPP1} gives the equivalence
of~$\HH\ttt$ and~$\HH\tt \opp \AA¥$. The argument is
similar for~\eqref{E:BluV3} and~\eqref{E:BluV4}, the only
difference being an additional pentagon relation for
replacing $\AA¥ \AA¥$ by~$\AA\1 \AA¥ \AA\0$ on the
right.

For~$\CC¥$, with similar notation, we have $\tt
= \tti \op \ttii$ and $\ttt = \ttii
\op \tti$, and we find now
\begin{gather}
\label{E:BluV5}
\HH\ttt  = \PP\IIii\IIi \opp \HHH\ttii \opp \dd{\HH\tti},\\
\label{E:BluV6}
\HH\tt \opp \CC¥ = \PP\IIi\IIii \opp \HHH\tti \opp
\dd{\HH\ttii} \opp \CC¥,\\
\label{E:BluV7}
\HHH\ttt  = \PPP\IIii\IIi \opp \HHH\ttii \opp
\dd{\HHH\tti} \opp \AA¥,\\
\label{E:BluV8}
\HHH\tt \opp \CC\0 = \PPP\IIi\IIii \opp \HHH\tti \opp
\dd{\HHH\ttii} \opp \AA¥ \CC\0.
\end{gather}
By~\eqref{E:ComW}, we have $\HHH\ttii \opp
\dd{\HH\tti} \eqacs \ddd{\nnii}{\HH\tti} \opp
\HHH\ttii$, and the $\Racs$-equivalence of~$\HH\ttt$
and~$\HH\tt \opp \CC¥$ follows from the commutativity of
the diagram
$$\begin{CD}
\vine{\II}
@>\Q{\PP\IIi\IIii}>>
\vine{\IIi, \IIii}
@>\Q{\HHH\tti}>>
\vine{\tti, \IIii}
@>\Q{\dd{\HH\ttii}}>>
\tti \op \ttii \\
@|
@V{\Racs+\text{torsion~}\quad}V{\PP\nnii\nni}V
@V{\eqref{E:HHP2}\qqquad}V{\PP\nnii1}V
@V{\eqref{E:HHP1}\qquad}V{\PP11 = \CC¥}V\\
\vine{\II}
@>\Q{\PP\IIii\IIi}>>
\vine{\IIii, \IIi}
@>\Q{{\ddd\nnii{\HH\tti}}}>>
\vine{\IIii, \tti}
@>\Q{{\HHH\ttii}}>>
\ttii \op \tti
\end{CD}$$
The commutativity of the left square follow from
the fact that both $\PP\IIi\IIii \opp
\PP\nnii\nni$ and $\PP\IIii\IIi$ induce the
same permutation of the labels: it follows that
these words must be equivalent with respect to
any family of relations that makes a presentation
of the symmetric group, and, therefore, they
are equivalent under the Coxeter
relations of~$\Racs$ completed with the torsion
relations~$\cc i^2=\ss i^2= 1$.

 The argument is similar for~$\HHH\ttt$. First
\eqref{E:ComW} gives
$\HHH\ttii \opp
\dd{\HHH\tti} \eqacs \ddd{\nnii}{\HHH\tti} \opp
\HHH\ttii$, and the rest is the commutativity of 
$$\begin{CD}
\vine{\II, \ttt}
@>{\PPP\IIi\IIii}>>
\vine{\IIi, \IIii, \ttt}
@>{\HHH\tti}>>
\vine{\tti, \IIii, \ttt}
@>{\dd{\HHH\ttii}}>>
\vine{\tti, \ttii, \ttt}
@>{\AA¥}>>
\vine{\tti \op \ttii, \ttt} \\
@|
@V{\Racs+\text{torsion}\quad}V{\PPP\nnii\nni}V
@V{\eqref{E:HHP2}\qquad}V{\PPP\nnii1}V
@V{\eqref{E:HHP1}\qquad}V{\PPP11=\SS¥}V
@V{\eqref{E:RACS1}\qquad}V{\CC\0}V\\
\vine{\II,\ttt}
@>{\PPP\IIii\IIi}>>
\vine{\IIii, \IIi, \ttt}
@>{{\ddd\nnii{\HHH\tti}}}>>
\vine{\IIii, \tti, \ttt}
@>{{\HHH\ttii}}>>
\vine{\ttii, \tti, \ttt}
@>{\AA¥}>>
\vine{\ttii \op \tti, \ttt}
\end{CD}$$
The induction is now easy, and there is no need to
consider the case of~$\AA\a$ and~$\CC\a$
separately. So we use~$\BB\a$ to represent the
two cases simultaneously. Assume $\a = \0\b$.
Then we have $\ttt = \ttti \op \ttii$ with $\ttti = \tti
\act \BB\a$, and we find
\begin{align*}
\HH\ttt 
= \PP\IIi\IIii \opp \HHH\ttti \opp \dd{\HH\ttii}
&\eqIH
\PP\IIi\IIii \opp \HHH\tti \opp \BB{\0\b} \opp
\dd{\HH\ttii}\\
&\eqG
\PP\IIi\IIii \opp \HHH\tti \opp
\dd{\HH\ttii} \opp \BB{\0\b}
= \HH\tt \opp \BB\a,\\
\HHH\ttt 
= \PPP\IIi\IIii \opp \HHH\ttti \opp \dd{\HHH\ttii}
\opp \AA¥
&\eqIH
\PPP\IIi\IIii \opp \HHH\tti \opp \BB{\0\b} \opp
\dd{\HHH\ttii} \opp \AA¥ \\
&\eqG
\PPP\IIi\IIii \opp \HHH\tti \opp
\dd{\HH\ttii} \opp \BB{\0\b} \opp \AA¥\\
& \eqG
\PPP\IIi\IIii \opp \HHH\tti \opp
\dd{\HH\ttii} \opp \AA¥ \opp \BB{\0\0\b}
= \HHH\tt \opp \BB{\0\a}.
\end{align*}
The argument is symmetric (and simpler:
no commutation is needed) in the case~$\a = \1\b$.
\end{proof}

Applying Proposition~\ref{P:Criterion}, we obtain

\begin{prop} \label{P:Prev}
The relations~$\RACS$ completed with the
torsion relations $\CC\a^2 = \SS\a^2 = 1$, make
a presentation of the group~$\GGAC$, \ie, of
Thompson's group~$V$, in terms of the
generators~$\AA\a, \CC\a$ and~$\SS\a$.
\end{prop}

As the relations~$\RACS$ follow from those
of~$\RAC$ and the definition of~$\SS\a$, we
immediately deduce:

\begin{prop} \label{P:PreV}
The relations~$\RAC$, \ie, the geometric relations, completed
with the pentagon and hexagon relations, and the
torsion relations $\CC\a^2 = 1$, make
a presentation of~$V$ in terms of the
generators~$\AA\a$ and~$\CC\a$.
\end{prop}

As in the case of the group~$F$, we can
restrict to the generators~$\aa i, \cc i$
and~$\ss i$. By looking at the proof of
Lemma~\ref{L:BluV}, we see that, if
$\ttt = \tt \act \bb i$ holds with
$\bb¥$ is $\aa¥, \cc¥$, or~$\ss¥$, then we
have
\begin{equation} \label{E:BluV}
\HH\ttt \eqacs \HH\tt \opp \bb i.
\end{equation}
Applying Proposition~\ref{P:Criterion} once more, we
deduce

\begin{prop}
The relations~$\Racs$ completed with the torsion
relations~$\cc i^2 = \ss i^2 = 1$ make a
presentation of the group~$\GGAC$,
\ie, of~$V$, in terms of the generators~$\aa i,
\cc i$ and~$\ss i$. 
\end{prop}

Finally, as all relations in~$\Racs$ follow
from~$\Rac$, we also obtain

\begin{prop}
The relations~$\Rac$ completed with the torsion
relations~$\cc i^2 = 1$ make a
presentation of the group~$\GGAC$,
\ie, of~$V$, in terms of the generators~$\aa
i$ and~$\cc i$. 
\end{prop}

\section{Semi-commutativity and the
group~$\Sys$} \label{S:PreU}

We have seen how to naturally connect 
Thompson's group~$V$ with the associativity
and commutativity laws. Inspecting
the computations of Section~\ref{S:PreV}, we see
that the main technical role is played
by the elements~$\SS\a$. This suggests to 
introduce the subgroup of~$\GGAC$ generated by
the elements~$\AA\a$ and~$\SS\a$. We
shall see now that the latter naturally arises as
a geometry group, namely that of
associativity together with a weak form of
commutativity.

\subsection{The semi-commutativity law}

\begin{defi}
We define {\it (left) semi-commutativity} to be
the law
\begin{equation}
\tag{$\Sem$}
x(yz) = y(xz).
\end{equation}
\end{defi}

As associativity and semi-commutativity are linear
laws in the sense of Section~\ref{S:PreV}, they
give rise to a geometry group~$\GGAS$.

\begin{prop}
The group~$\GGAS$ is (isomorphic to) the
subgroup~$\Sys$ of~$V$ generated by the
elements~$\AA\a$ and~$\SS\a$, \ie, $\Sys$ is the
geometry group of associativity and
semi-commutativity.
\end{prop}

\begin{proof}
Figure~\ref{F:ActS} shows that  the
operators associated with the
semi-commuta\-tivity law are the
operators~$\SSp\a$ of Section~\ref{S:PreV}, so
the geometry monoid~$\GMAS$ is the submonoid
of~$\GMAC$ generated by the
operators~$\AApm\a$ and~$\SSpm\a$.
Quotienting under near-equality gives a similar
relation for the geometry groups.
\end{proof}

So, in particular, if we extract from the relations
established for~$\GGAC$ those that involve the
generators~$\AA\a$ and~$\SS\a$ only, the latter
have to be satisfied in the group~$\GGAS$.

\begin{defi} 
We define~$\RAS$ to consist of the translated
copies of
\begin{gather}
\tag{$\RGperp$}
\BB{\0\a} \opp \BBB{\1\b} = \BBB{\1\b}
\opp \BB{\0\a},\\
\tag{$\RGA$}
\BB{\1\1\a} \opp \AA¥ = \AA¥ \opp\BB{\1\a},
\qquad
\BB{\1\0\a} \opp\AA¥ = \AA¥ \opp\BB{\0\1\a},
\qquad
\BB{\0\a} \opp\AA¥ = \AA¥ \opp\BB{\0\0\a},\\
\tag{$\RGS$}
\BB{\1\1\a} \opp\SS¥ = \SS¥ \opp\BB{\1\1\a},
\qquad
\BB{\1\0\a} \opp\SS¥ = \SS¥ \opp\BB{\0\a},
\qquad
\BB{\0\a} \opp\SS¥ = \SS¥ \opp\BB{\1\0\a},
\end{gather}
with $\BB¥, \BBB¥ = \AA¥, \SS¥$, plus the
translated copies of
\begin{gather}
\tag{$\RP$}
\AA¥ \AA¥ =\AA\1 \AA¥ \AA\0,\\
\label{E:RAS1}
\SS¥ \AA\1 \AA¥ = \AA\1 \AA¥ \SS\0, 
\quad
\SS\1 \SS¥ \AA\1 = \AA¥ \SS¥, 
\quad
\SS¥ \SS\1 \AA¥ = \AA\1 \SS¥,
\quad 
\SS¥ \SS\1 \SS¥ = \SS\1 \SS¥ \SS\1.
\end{gather}
\end{defi}

\begin{prop} \label{P:RelU}
All relations of~$\RAS$, as well as~$\SS\a^2 =
1$, are satisfied in~$\GGAS$, \ie, in~$\Sys$.
\end{prop}

We also consider the subfamily of~$\RAS$
associated with the elements of~$\Ga$ and~$\Gs$.

\begin{defi}
We define~$\Ras$ to consist of the following
relations:
\begin{gather*}
\aa i \bb{j-1}
= \bb j \aa i
\text{\quad and \quad}
\ss i \bb j 
= \bb j \ss i
\text{\quad for $j \ge i + 2$ and $\bb¥ =
\aa¥$ or~$\ss¥$},\\
\ss i \ss{i+1} \ss i 
= \ss{i+1} \ss i \ss{i+1},
\text{\qquad}
\ss{i+1} \ss i \aa{i+1}
= \aa i \ss i,
\text{\qquad}
\ss i \ss{i+1} \aa i 
= \aa{i+1} \ss i.
\end{gather*}
\end{defi}

\begin{lemm} 
All relations of~$\Ras$ are satisfied in~$\GGAS$,
\ie, in~$\Sys$.
\end{lemm}

Actually, it is easy to check that the relations
of~$\Ras$ follow from those of~$\RAS$ plus the
definitions $\ss i = \SS{\1^{i-1}}$.

\subsection{Presentations of~$\Sys$}

Our aim is to prove:

\begin{prop}
The family~$\RAS$ completed with the torsion
relations $\SS\a^2 = 1$ make a presentation
of~$\Sys$, in terms of the generators~$\AA\a$
and~$\SS\a$.
\end{prop}

\begin{proof}
The method should be clear: we select a
family of trees containing one element in each
$\Sys$-orbit, then define distinguished words
in~$\WW{\GA,\GS}$ describing how to construct
a tree starting from the distinguished element of
its orbit, and, finally, check that there are enough
relations in~$\RAS$ to witness for the
relations~\eqref{E:Criterion} of
Proposition~\ref{P:Criterion}.

The construction is a slight modification of the
one used in Section~\ref{S:PreV}. The difference
between commutativity and semi-commutativity
is that the latter cannot change the rightmost label
of a tree. To keep the same conventions as in
Section~\ref{S:PreV}, let~$\Tii$ denote the
subset of~$\Ti$ made by coloured trees in
which the rightmost leaf wears the maximal
label. Then every tree in~$\Tii$ is equivalent
up to associativity and semi-commutativity to
some right vine~$\vine{\II}$. For
such a tree~$\tt$, the word~$\HH\tt$
maps~$\vine{\II}$ to~$\tt$, and, by
construction, $\HH\tt$ consists of letters~$\aa i$
and~$\ss j$ exclusively, since the rightmost leaf
is never changed. Indeed, the only letter~$\cc i$
possibly occurring in~$\HH\tt$ comes from the
factors~$\PP\II\JJ$ in the inductive construction,
and this happens only when $\II$ contains the
largest element of~$\II
\cup \JJ$. We can therefore use the
words~$\HH\tt$ and~$\HHH\tt$ without change.
Then the only point is to check that $\ttt = \tt
\act \BB\a$ implies
\begin{equation}
\HH\ttt \eqas \HH\tt \opp \BB\a
\text{\quad and \quad}
\HHH\ttt \eqas \HHH\tt \opp \BB{\0\a}
\end{equation}
both in the case $\BB¥ = \AA¥$ and $\BB¥ = \SS¥$.
For the case of~$\AA\a$, it suffices to look at the
proof of Lemma~\ref{L:BluV}. The case of~$\SS¥$ 
has not been considered in Section~\ref{S:PreV},
and we consider it now. So we assume $\tt = \tti
\op (\ttii \op \ttiii)$ and $\ttt = \ttii \op (\tti
\op \ttiii)$. We obtain
\begin{gather}
\label{E:BluS5}
\HH\ttt  = \PPP\IIii{\IIi\cup\IIiii} \opp \HHH\ttii
\opp  \dd{\PPP\IIi\IIiii} \opp \dd{\HHH\tti} \opp
\ddd2{\HH\ttiii},\\
\label{E:BluS6}
\HH\tt \opp \SS¥ = \PPP\IIi{\IIii\cup\IIiii} \opp
\HHH\tti
\opp  \dd{\PPP\IIii\IIiii} \opp \dd{\HHH\ttii} \opp
\ddd2{\HH\ttiii}\opp \SS¥.
\end{gather}
By~\eqref{E:ComW}, we have $\HHH\ttii
\opp  \dd{\PPP\IIi\IIiii} \eqas
\ddd{\nnii}{\PPP\IIi\IIiii} \opp \HHH\ttii$,
$\HHH\tti
\opp  \dd{\PPP\IIii\IIiii} \eqas
\ddd{\nni}{\PPP\IIii\IIiii} \opp \HHH\tti$, and 
$\HHH\ttii \opp
\dd{\HH\tti} \eqas \ddd{\nnii}{\HH\tti} \opp
\HHH\ttii$. Then the $\RAS$-equivalence
of~$\HH\ttt$ and~$\HH\tt \opp \SS¥$ follows
from the commutativity of the diagram
$$\hspace{-4mm}\begin{CD}
\vine{\II}
@>{\PPP\IIi{\IIii\cup \IIiii} \opp
\ddd\nni{\PPP\IIii\IIiii}}>> 
\vine{\IIi, \IIii, \IIiii}
@>{\HHH\tti}>>
\vine{\tti, \IIii, \IIiii}
@>{\dd{\HHH\ttii}}>>
\vine{\tti, \ttii, \IIiii}
@>{\ddd2{\HHH\ttiii}}>>
\vine{\tti, \ttii, \ttiii} \\
@|
@V{\Ras+torsion\qqquad}V{\PPP\nnii\nni}V
@V{\eqref{E:HHP2}\qquad}V{\PPP\nnii1}V
@V{\eqref{E:HHP1}\qquad}V{\PPP11}V
@V{\eqref{E:RAS1}\qquad}V{\PPP11=\SS¥}V\\
\vine{\II}
@>{\PPP\IIii{\IIi\cup \IIiii} \opp
\ddd\nnii{\PPP\IIi\IIiii}}>> 
\vine{\IIii, \IIi, \IIiii}
@>{\ddd\nnii{\HHH\tti}}>>
\vine{\IIii, \tti, \IIiii}
@>{\HHH\ttii}>>
\vine{\ttii, \tti, \IIiii}
@>{\ddd2{\HHH\ttiii}}>>
\vine{\ttii, \tti, \ttiii}
\end{CD}$$
The relations of~$\RAS$ are sufficient to obtain
the commutativity of the last three squares. As for
the first square, the associated permutations are
equal, so the relations of~$\Ras$ completed with
the torsion relations $\ss i^2 = 1$ must give the
result.

The argument is similar for the words~$\HHH\tt$,
with an associated diagram coinciding with the
above one up to an additional square on the right
whose commutativity is provided by the
relation $\SS¥ \AA\1 \AA¥ = \AA\1 \AA¥ \SS\0$.
The induction along addresses is similar to the
one we used for the groups~$\GGA$
and~$\GGAC$, \ie,
for~$F$ and~$V$.
\end{proof}

As in Section~\ref{S:PreF} and~\ref{S:PreV},
we deduce that there are enough relations in
the list~$\Ras$ to generate all needed
equivalences, and we conclude:

\begin{prop} \label{P:PresVLas}
The group~$\Sys$ is generated by~$\Ga$
and~$\Gs$, and the relations~$\Ras$ completed
with $\ss i^2 = 1$ make a presentation
of~$\Sys$ in terms of these generators.
\end{prop}

\begin{coro}
The group~$\Sys$ is isomorphic to the
group~$\widehat{V}$ of~\cite{Bri2}.
\end{coro}

\section{The group~$\Bs$ and its connection to
twisted semi-commutativity}
\label{S:PreB}

The presentation of the group~$\Sys$ in terms of
the~$\aa i$'s and the~$\ss i$'s given in
Proposition~\ref{P:PresVLas} includes the Coxeter
presentation of the symmetric group~$\Syi$ in
terms of the~$\ss i$'s. Following the example of
Artin's braid group~$B_\infty$, which can be
defined by removing the torsion relations~$\ss i^2
= 1$ in the Coxeter presentation of~$\Syi$, or,
more generally, of Artin--Tits groups, we
introduce the group obtained from~$\Sys$ by
removing the torsion relations. This is specially
natural as we can see that the torsion relations
play a very small role in the computations of the
previous sections. This new group, here
denoted~$\Bs$, has rich properties, investigated
in~\cite{Dhe} and~\cite{Bri1, Bri2}. In this
paper, we study~$\Bs$ from the point of
view of geometry groups only. The main result is
that $\Bs$ is the geometry group of associativity
together with some twisted version of
semi-commutativity. This in particular provides a
concrete realization of~$\Bs$ as a group of
partial operators on coloured trees.

\subsection{The group~$\Bs$}

As is usual with permutations and braids, we
use~$\sss i$ for the torsion free lifting of the
generator~$\ss i$. Accordingly, we use~$\Gss$ for
the infinite family $\sss1, \sss2, \pp$, and
$\Rass$ for a copy of~$\Ras$ with~$\sss i$
replacing~$\ss i$ everywhere. 

\begin{defi}
We define~$\Bs$ to be the group $\Gr(\Ga,\Gss;
\Rass)$, \ie, the group generated by two
infinite sequences $\aa1, \aa2, \pp$, $\sss1,
\sss2, \pp$ with the relations
\begin{equation}\label{E:PresG}
\begin{cases}
\quad \aa i \bb{j-1} = \bb j \aa i
\text{\quad and \quad}
\sss i \bb j  = \bb j \sss i
\text{ \quad for $j \ge i + 2$ and $\bb¥ =
\aa¥$ or~$\sss¥$},\\
\quad \sss i \sss{i+1} \sss i 
= \sss{i+1} \sss i \sss{i+1},
\text{\qquad}
\sss{i+1} \sss i \aa{i+1}
= \aa i \sss i,
\text{\qquad}
\sss i \sss{i+1} \aa i 
= \aa{i+1} \sss i.
\end{cases}
\end{equation}
\end{defi}

Our current notation is chosen to emphasize the
similarity between~$\Bs$ and Artin's braid
group~$\Bi$: as shown in~\cite{Dhe}, the elements
of~$\Bs$ admit a natural realization
in terms of parenthesized braid diagrams, which are
analogous to ordinary braid diagrams but with
non-uniform distances between the strands. In
this framework, $\sss i$ correspond to a standard
crossing, while $\aa i$ corresponds to a rescaling
operator that shrinks the distances around the
$i$th position. The explicit presentation also
shows:

\begin{prop}
The group~$\Bs$ is isomorphic to the
group~$\widehat{BV}$ of~\cite{Bri2}.
\end{prop}

The group~$\Bs$ is a sort of twisted product
of Thompson's group~$F$ and Artin's braid
group~$\Bi$, and it is not surprising that it can be
investigated by the same methods as~$F$
and~$\Bi$. In particular, $\Bs$ is a group of
left fractions for the monoid with the
same presentation~\cite{Bri1, Dhe} and, as least
left common multiples exist in this monoid, the
group~$\Bs$ is torsion free. 

\subsection{Twisted commutation and
semi-commutation}

We turn to the realization of~$\Bs$ as a geometry
group, as we did for~$F, V$, and~$\Sys$. Applying
(semi)-commutativity is an involutive operation,
while $\Bs$ is torsion-free. So we are led
to considering non-involutive variants of
(semi)-commutativity. A natural way for making
commutativity operators non-involutive is to
assume that subtrees are changed when they are
switched. The simplest case is when only one
subtree is changed, and the new subtree depends on
the two subtrees that have been exchanged only.
This amounts to assuming that there exists a
binary operation on trees. 

\begin{figure} [htb]
\begin{picture}(70, 23)(0, 0)
\put(0,3){\includegraphics{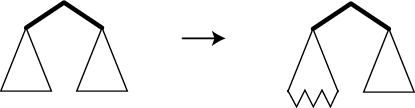}}
\put(4,2){$\tti$}
\put(16.5,2){$\ttii$}
\put(49,0){$\LD\tti\ttii$}
\put(65.5,2){$\tti$}
\put(33,17){$\CCtp¥$}
\end{picture}
\caption{\smaller The twisted commutation
operator~$\CCtp¥$}
\label{F:ActCt}
\end{figure}

\begin{defi}
(Figure~\ref{F:ActCt})
Assume that $\Tree$ is a set of trees equipped
with a binary operation~$\ld$. Then we
define the {\it $\Tree$-twisted commutation}
operator~$\CCtp¥$ by
\begin{equation}
\CCtp¥ : \tti \opp \ttii \longmapsto
\LD\tti\ttii \opp \tti.
\end{equation}
\end{defi}

So we still switch the left and the right
subtrees but, in the transformation, the right
subtree is (possibly) changed when it crosses the
left subtree. The bracket notation is chosen to
emphasize that $\LD\tti\ttii$ is the image
of~$\ttii$ under the action of~$\tti$. Note that
the standard commutation operator~$\CCp¥$
corresponds to using the trivial operation
$\LD\tti\ttii = \ttii$.

As in the case of the operators~$\CCp\a$, we
define~$\CCtp\a$ to be the translated
operator~$\DD\a{\CCtp¥}$, \ie,
$\CCtp¥$ acting on the $\a$-subtree. As for
inverses, the operators~$\CCtp\a$ need not be
injective in general, but we have the following
criterion:

\begin{lemm}
Assume that $\Tree$ is a set of trees equipped
with a bracket operation. Then the 
operators~$\CCtp\a$ are injective if and only the
bracket on~$\Tree$ is left cancellative,
\ie, 
\begin{equation} \label{E:Canc}
\LD\tt\tti = \LD\tt\ttii
\text{\quad implies \quad}
\tti = \ttii.
\end{equation}
\end{lemm}

Under such an hypothesis, the inverse
operator of~$\CCtp\a$ is still a partial operator
on~$\Tree$. 

\begin{figure} [htb]
\begin{picture}(78,25)(0, 0)
\put(0,3){\includegraphics*{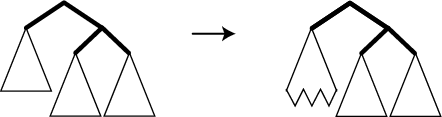}}
\put(3,4){$\tti$}
\put(12,0){$\ttii$}
\put(21,0){$\ttiii$}
\put(48.5,1.5){$\LD\tti\ttii$}
\put(60,0){$\tti$}
\put(70,0){$\ttiii$}
\put(34,19){$\SStp¥$}
\end{picture}
\caption{\smaller The twisted semi-commutation
operator~$\SStp¥$}
\label{F:ActSt}
\end{figure}

As we chosed to investigate the torsion-free
version~$\Bs$ of~$\Sys$ rather than that
of~$V$, we are led to considering a twisted
version of semi-commutation too. We keep the
definition of Section~\ref{S:PreV}, \ie, we define
the twisted version~$\SStp¥$ of~$\SSp¥$ by
$\SStp¥ = \CCtp¥ \AAm¥ \CCtm\1$, which corresponds
to:

\begin{defi}
(Figure~\ref{F:ActSt})
Assume that $\Tree$ is a set of trees equipped
with a binary operation~$\ld$. Then we
define the {\it $\Tree$-twisted semi-commutation}
operator~$\SStp¥$  by
\begin{equation}
\SStp¥ : \tti \opp (\ttii \opp \ttiii)
\longmapsto 
\LD\tti\ttii \opp (\tti \opp \ttiii).
\end{equation}
\end{defi}

We naturally define~$\SStp\a$ to be the
$\a$-translated copy of~$\SStp¥$. Under the
hypothesis that the bracket on~$\Tree$ is left
cancellative, the operator~$\SStp\a$ is injective,
and its inverse~$\SStm\a$ is a partial operator. 
The (semi)-commutation operators correspond to
no algebraic law, but we still have a family of
partial self-injections of a set of trees, and it
is natural to consider the monoids they generate:

\begin{defi}
Assume that $\Tree$ is a family of trees equipped
with a left cancellative bracket operation. Then we
define $\GMACT$ (\resp $\GMAST$) to be the
monoid generated by the operators~$\AApm\a$
and~$\CCtpm\a$ (\resp the
operators~$\AApm\a$ and~$\SStpm\a$) acting
on~$\Tree$.
\end{defi}

Our aim is now to investigate the monoids
($\GMACT$ and) $\GMAST$ for
appropriate choices of the bracket operation. When
$\Tree$ is equipped with the trivial bracket
$\LD\tti\ttii = \ttii$, we find
$$\GMACT =  \GMAC
\text{\quad and \quad}
\GMAST =  \GMAS,$$
\ie, we come back to the framework of
Sections~\ref{S:PreV} and~\ref{S:PreU}.

\subsection{LD-systems}

In general, the twisted operators~$\CCtp\a$
and~$\SStp\a$ need not satisfy the same relations
as their standard versions. However it is easy to list
the requirements needed for the relations
of~$\RACS$ to be valid in the monoid~$\GMAST$.

\begin{prop} \label{P:Cond}
(i) The relations $\AAp\1 \SStp¥ = \SStp¥
\SStp\1
\AAp¥$, 
$\AAp¥ \SStp¥ = \SStp\1 \SStp¥ \AAp\1$, and 
$\SStp¥ \SStp\1 \SStp¥ = \SStp\1 \SStp¥ \SStp\1$
hold in the monoid $\GMAST$ if and only if, for all trees~$\tti,
\ttii,
\ttiii$ in~$\Tree$, we have
\begin{gather}
\label{E:Ide1}
\LD\tti{\ttii \op \ttiii}
= \LD{\tti}{\ttii} \opp \LD{\tti}{\ttii},\\
\label{E:Ide2}
\LD{(\tti \op \ttii)}{\ttiii}
= \LD{\tti}{\LD{\ttii}{\ttiii}},\\
\label{E:Ide3}
\LD\tti{\LD\ttii\ttiii}
= \LD{\LD\tti\ttii}{\LD\tti\ttiii}.
\end{gather}

(ii) Assume that $\Tree$ is the set of all
$\LDS$-coloured trees for some set~$\LDS$. 
Then the conditions of~$(i)$ are satisfied if
and only if there exists a left cancellative
left self-distributive bracket operation
on~$\LDS$ such that, for all trees~$\tti, \ttii$
in~$\Tree_\LDS$, the tree~$\LD\tti\ttii$ is
obtained by replacing each label~$\yy$
in~$\ttii$ with $\LD{\xx_1}{\LD{\xx_2}{\pp
\LD{\xx_n}{\yy}\pp}}$, where $(\xx_1, \pp,
\xx_n)$ is the left-to-right enumeration of the
labels in~$\tti$. 

(iii) In this case, all relations of~$\RACS$ are
satisfied by the operators~$\AAp\a$, $\CCtp\a$,
and~$\SStp\a$, and the torsion relations
$\CCtp\a{}^2 \approx \SStp\a{}^2 \approx \id$
are satisfied if and only  if, for all trees~$\tti,
\ttii$, we have
\begin{equation}
\label{E:Ide4}
\LD\tti{\LD\tti\ttii} = \ttii.
\end{equation}
\end{prop}

\begin{figure} [htb]
$$\includegraphics*{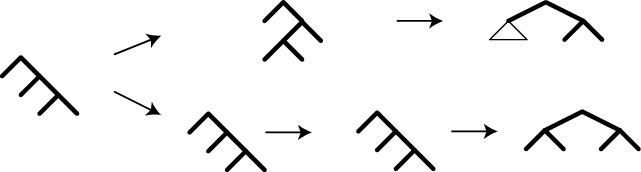}$$
\begin{picture}(0, 0)(60, 0)
\put(24, 30){$\AAp\1$}
\put(74, 33){$\SStp¥$}
\put(24, 13){$\SStp¥$}
\put(51, 15){$\SStp\1$}
\put(83, 15){$\AAp¥$}
\put(3, 18){$\tti$}
\put(7, 15){$\ttii$}
\put(11, 12){$\ttiii$}
\put(17, 12){$\ttiiii$}
\put(47, 29){$\tti$}
\put(50, 21){$\ttii$}
\put(55, 21){$\ttiii$}
\put(59, 24){$\ttiiii$}
\put(84, 24){$\LD\tti{\ttii\op\ttiii}$}
\put(100, 24){$\tti$}
\put(106, 24){$\ttiiii$}
\put(31, 9){$\LD\tti\ttii$}
\put(39, 6){$\tti$}
\put(43, 3){$\ttiii$}
\put(49, 3){$\ttiiii$}
\put(59, 9.5){$\LD\tti\ttii$}
\put(63, 5.5){$\LD\tti\ttiii$}
\put(71, 3){$\tti$}
\put(77, 3){$\ttiiii$}
\put(89, 6){$\LD\tti\ttii$}
\put(97.5, 6){$\LD\tti\ttiii$}
\put(106, 6){$\tti$}
\put(112, 6){$\ttiiii$}
\end{picture}
\caption{\smaller The relation $\AAp\1 \SStp¥ = 
\SStp¥ \SStp\1 \AAp¥$ requires  $\LD{\tti}{\ttii
\op \ttiii}  = \LD{\tti}{\ttii} \opp
\LD{\tti}{\ttii}$}
\label{F:Rel1}
\end{figure}

\begin{figure} [htb]
$$\includegraphics*{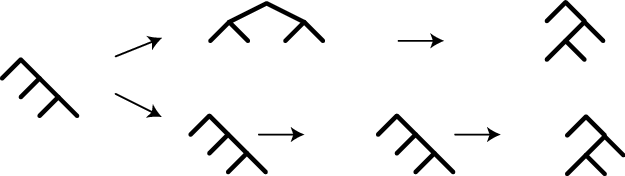}$$
\setlength{\unitlength}{1mm}
\begin{picture}(0, 0)(60, 0)
\put(27,29){$\AAp¥$}
\put(27,13){$\SStp\1$}
\put(53,14){$\SStp¥$}
\put(85,14){$\AAp¥$}
\put(75,30){$\SStp¥$}

\put(5,19){$\tti$}
\put(9,15.5){$\ttii$}
\put(13,12){$\ttiii$}
\put(18,12){$\ttiiii$}

\put(41,25){$\tti$}
\put(47.5,25){$\ttii$}
\put(54,25){$\ttiii$}
\put(60.5,25){$\ttiiii$}

\put(37,10){$\tti$}
\put(36,6){$\LD\ttii\ttiii$}
\put(45,3){$\ttii$}
\put(50,3){$\ttiiii$}

\put(61,9.2){$\LD\tti{\LD\ttii\ttiii}$}
\put(73,6){$\tti$}
\put(77,3){$\ttii$}
\put(82,3){$\ttiiii$}

\put(92,9.2){$\LD\tti{\LD\ttii\ttiii}$}
\put(102,3){$\tti$}
\put(108,3){$\ttii$}
\put(112,6){$\ttiiii$}

\put(86,28.5){$\LD{(\tti\op\ttii)}\ttiii$}
\put(99,22){$\tti$}
\put(105,22){$\ttii$}
\put(109,25){$\ttiiii$}
\end{picture}
\vspace{0mm}
\caption{\smaller The relation $\AAp¥ \SStp¥ 
= \SStp\1 \SStp¥ \AAp\1$ requires 
$\LD{(\tti\op\ttii)}\ttiii
= \LD\tti{\LD\ttii\ttiii}$}
\label{F:Rel2}
\end{figure}

\begin{figure} [htb]
$$\includegraphics*{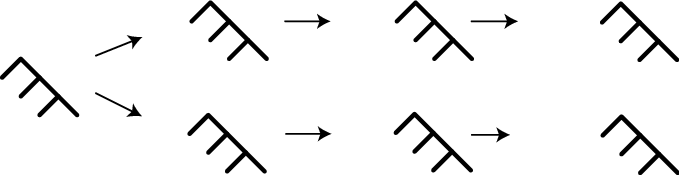}$$
\setlength{\unitlength}{1mm}
\begin{picture}(0, 0)(60, 0)
\put(20,30){$\SStp\1$}
\put(20,13){$\SStp¥$}
\put(52, 34){$\SStp¥$}
\put(52,14.5){$\SStp\1$}
\put(83,34){$\SStp\1$}
\put(83,14.5){$\SStp¥$}

\put(1,19){$\tti$}
\put(4.5,15.5){$\ttii$}
\put(8,12){$\ttiii$}
\put(14,12){$\ttiiii$}

\put(33,28.5){$\tti$}
\put(32,25){$\LD\ttii\ttiii$}
\put(40,22){$\ttii$}
\put(46,22){$\ttiiii$}

\put(59,28.5){$\LD\tti{\LD\ttii\ttiii}$}
\put(71,25){$\tti$}
\put(75,22){$\ttii$}
\put(81,22){$\ttiiii$}

\put(28.5,9.5){$\LD\tti\ttii$}
\put(36,6){$\tti$}
\put(40,3){$\ttiii$}
\put(46,3){$\ttiiii$}

\put(62,9.5){$\LD\tti\ttii$}
\put(67,6){$\LD\tti\ttiii$}
\put(75,3){$\tti$}
\put(81,3){$\ttiiii$}

\put(90,9.5){$\LD{\LD\tti\ttii}{\LD\ttii\ttiii}$}
\put(102,6){$\LD\tti\ttii$}
\put(110,3){$\tti$}
\put(116,3){$\ttiiii$}

\put(94,28.5){$\LD\tti{\LD\ttii\ttiii}$}
\put(101,25){$\LD\tti\ttii$}
\put(110,22){$\tti$}
\put(116,22){$\ttiiii$}
\end{picture}
\vspace{0mm}
\caption{\smaller The relation
$\SStp\1 \SStp¥  \SStp\1=\SStp¥ 
\SStp\1 \SStp¥$ requires 
$\LD\tti{\LD\ttii\ttiii} =
\LD{\LD\tti\ttii}{\LD\tti\ttiii}$}
\label{F:Rel3}
\end{figure}

\begin{figure} [htb]
$$\includegraphics*{Hexagon.eps}$$
\setlength{\unitlength}{1mm}
\begin{picture}(0, 0)(60, 0)
\put(23,29){$\AAp¥$}
\put(23,11){$\CCtp\1$}
\put(58,34){$\CCtp¥$}
\put(58,12){$\AAp¥$}
\put(92,29){$\AAp¥$}
\put(92,11){$\CCtp\0$}

\put(5.5,18){$\scriptstyle1$}
\put(11.5,15){$\scriptstyle2$}
\put(17,15){$\scriptstyle3$}

\put(36,26){$\scriptstyle1$}
\put(42,26){$\scriptstyle2$}
\put(48,30){$\scriptstyle3$}

\put(36,7){$\scriptstyle1$}
\put(40.5,3){$\scriptstyle\LD23$}
\put(48,3){$\scriptstyle2$}

\put(67,30){$\scriptstyle\LD1{\LD23}$}
\put(76,26){$\scriptstyle1$}
\put(82,26){$\scriptstyle2$}

\put(71,3){$\scriptstyle1$}
\put(75,3){$\scriptstyle\LD23$}
\put(82,7){$\scriptstyle2$}

\put(101,15){$\scriptstyle1$}
\put(105,15){$\scriptstyle\LD23$}
\put(112.5,18){$\scriptstyle2$}
\end{picture}
\vspace{0mm}
$$\includegraphics*{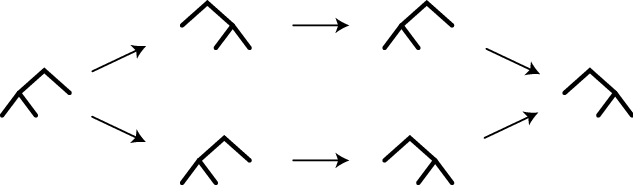}$$
\setlength{\unitlength}{1mm}
\begin{picture}(0, 0)(60, 0)
\put(23,29){$\AAm¥$}
\put(23,11){$\CCtp\0$}
\put(58,34){$\CCtp¥$}
\put(58,12){$\AAm¥$}
\put(92,29){$\AAm¥$}
\put(92,11){$\CCtm\1$}

\put(5.5,15){$\scriptstyle1$}
\put(11.5,15){$\scriptstyle2$}
\put(17,18){$\scriptstyle3$}

\put(36,30){$\scriptstyle1$}
\put(42,26){$\scriptstyle2$}
\put(48,26){$\scriptstyle3$}

\put(35,3){$\scriptstyle\LD12$}
\put(42,3){$\scriptstyle1$}
\put(48,7){$\scriptstyle3$}

\put(69,26){$\scriptstyle\LD12$}
\put(75.5,26){$\scriptstyle\LD13$}
\put(82,30){$\scriptstyle1$}

\put(69,7){$\scriptstyle\LD12$}
\put(76,3){$\scriptstyle1$}
\put(82,3){$\scriptstyle3$}

\put(99.5,18){$\scriptstyle\LD12$}
\put(105,15){$\scriptstyle\LD13$}
\put(112.5,15){$\scriptstyle1$}
\end{picture}
\vspace{-3mm}
\caption{\smaller The twisted
hexagon relations}
\label{F:Hext}
\end{figure}

\begin{figure} [htb]
\begin{picture}(0, 38)(50, 0)
\put(1,3){\includegraphics{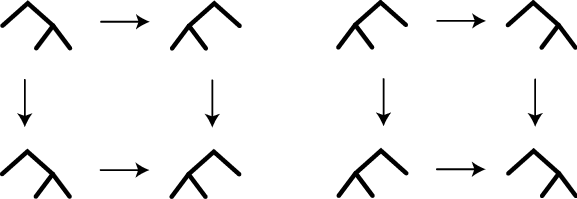}}
\put(20, 35){$\CCtp¥$}
\put(20, 10){$\CCtp¥$}
\put(77, 35){$\CCtp¥$}
\put(77, 10){$\CCtp¥$}
\put(7, 19){$\CCtp\1$}
\put(39, 19){$\CCtp\0$}
\put(68, 19){$\CCtp\0$}
\put(93.5, 19){$\CCtp\1$}
\put(0.5,29.5){$\scriptstyle1$}
\put(6,26){$\scriptstyle2$}
\put(12,26){$\scriptstyle3$}
\put(0.5,4){$\scriptstyle1$}
\put(4.5,0.5){$\scriptstyle\LD23$}
\put(12,0.5){$\scriptstyle2$}
\put(28,26){$\scriptstyle\LD13$}
\put(34,26){$\scriptstyle\LD12$}
\put(41,29.5){$\scriptstyle1$}
\put(25,0.5){$\scriptstyle\LD1{\LD23}$}
\put(34,0.5){$\scriptstyle\LD23$}
\put(41,4){$\scriptstyle2$}
\put(57.5,26){$\scriptstyle1$}
\put(63,26){$\scriptstyle2$}
\put(69,29.5){$\scriptstyle3$}
\put(57.5,0.5){$\scriptstyle1$}
\put(62.53,0.5){$\scriptstyle\LD23$}
\put(69,4){$\scriptstyle2$}
\put(83,29.5){$\scriptstyle\LD1{\LD23}$}
\put(92,26){$\scriptstyle1$}
\put(98,26){$\scriptstyle2$}
\put(81,4){$\scriptstyle\LD{\LD12}{\LD13}$}
\put(92,0.5){$\scriptstyle1$}
\put(98,0.5){$\scriptstyle2$}
\end{picture}
\caption{\smaller The twisted
$\CC¥$-geometric relations}
\label{F:GeoCt}
\end{figure}

\begin{proof}
For~$(i)$, the verifications are given in
Figures~\ref{F:Rel1}, \ref{F:Rel2},
and~\ref{F:Rel3}, respectively. Then $(ii)$ follows
from an induction on the size of the trees~$\tti$
and~$\ttii$. Finally, in order to establish~$(iii)$,
it suffices to check the $\CC¥$-geometric
relations and the hexagon relations, which is done
in Figures~\ref{F:Hext} and~\ref{F:GeoCt}.
\end{proof}

\begin{rema}
As we are mostly interested in the group~$\Bs$,
we concentrated on the constraints guaranteeing
that the relations of~$\RAS$ are satisfied, and
we saw that all relations of~$\RACS$ are then
valid. If we start with the operators~$\CCtp\a$
and require that the relations of~$\RAC$ be
satisfied, we come up with exactly the same
constraints, as can be read in
Figures~\ref{F:Hext} and~\ref{F:GeoCt}.
\end{rema}

We shall therefore be interested in the sequel with
sets equipped with a left self-distributive
operation, \ie, a binary operation that satisfies the
algebraic law
\begin{equation}
\tag{$LD$}
\LD\xx{\LD\yy\zz} = \LD{\LD\xx\yy}{\LD\xx\zz}
\end{equation}
---or $\xx(\yy\zz) = (\xx\yy)(\xx\zz)$ when the
operation symbol is omitted.

\begin{defi}
An algebraic system consisting of a set equipped
with a left self-distributive operation is called an
{\it LD-system}. An LD-system is said to be left
cancellative if its left translations are injective, \ie,
if \eqref{E:Canc} holds; it is called an {\it
LD-quasigroup} (in~\cite{Dgd}) or a {\it rack}
(in~\cite{FeR}) if its left translations are bijective.
An LD-system is said to be {\it involutory} if
\eqref{E:Ide4} holds. Note that an involutory
LD-system is necessarily an LD-quasigroup.
\end{defi}

\begin{exam} \label{X:Conj}
Any set~$\LDS$ equipped with $\LD\xx\yy=\yy$ is
a (trivial) involutory LD-system. If $G$ is a
group, then $G$ equipped with
$\LD\xx \yy = \xx \yy \xx\inv$ is an
LD-quasigroup, denoted~$\conj(G)$ in the sequel.
\end{exam}

From now on, we always restrict to the
context of Proposition~\ref{P:Cond}$(ii)$, \ie,
consider the twisted (semi)-commutation
operators on the set~$\Tree_\LDS$ that stem from
some left cancellative LD-system~$\LDS$.
Accordingly, we shall simplify our notation, and
write $\GMASL$ for $\GM{\Ass,
\Sem^{\scriptscriptstyle \Tree_\LDS}}$, and,
similarly, $\GMACL$ for $\GM{\Ass,
\Com^{\scriptscriptstyle \Tree_\LDS}}$.

\subsection{Making groups}

As in the case of associativity and
semi-commutativity, and for each fixed left
cancellative LD-system~$\LDS$, one can derive a
group from the monoid~$\GMASL$ by identifying
near-equal operators. However, controlling a
possible collapsing is not trivial, as we are not
in the framework of linear algebraic laws. 

The problem is to show that the near-equality
relation~$\approx$ defines a congruence on the
monoid~$\GMASL$. As in Section~\ref{S:PreV},
the solution is to show that each operator
admits a convenient seed in order to deduce that
$\approx$ is transitive. Now the notions of a
substitution and, consequently, of a seed, have
to be adapted to our current context. At the
expense of considering coloured trees whose
labels are formal expressions containing
variables and bracket operations, one can show
that, in a convenient sense, the pair of
coloured trees $(\vine{1, 2}, \vine{\LD12,
1})$, \ie, $(\et_1 \op \et_2, \et_{\LD12} \op
\et_1)$, is a seed for the operator~$\CCtp¥$,
while $(\vine{1, 2, 3}, \vine{\LD12, 1, 3})$ is
a a seed for the operator~$\SStp¥$. The details
are easy in the case of an LD-quasigroup; in
the more general case of a left cancellative
LD-system, more care is needed, but all required
techniques are explained in Chapter~VIII
of~\cite{Dgd}. The key ingredient is the result
that, if two self-distributivity operators
(analogous to the current operators~$\AAp\a$
or~$\SSp\a$ but for the left self-distributivity
law) agree on some tree, then they agree
everywhere. All we need in the sequel is the
following result:

\begin{lemm}
Assume that $\LDS$ is a left cancellative
LD-system. Then near-equality is a congruence
on the monoid~$\GMASL$, and the action of the
latter on~$\Tree_\LDS$ induces a partial
action of the associated
quotient-group~$\GGASL$. 
\end{lemm}

The group~$\GGASL$ will naturally be called the
geometry group of associativity and
$\LDS$-twisted semi-commutativity. 
Proposition~\ref{P:Cond} directly
implies:

\begin{prop}
For each left cancellative LD-system~$\LDS$,
the group~$\GGASL$ is a quotient of~$\Bs$.
\end{prop}

The group~$\GGASL$ depends on the
considered LD-system~$\LDS$. For instance, when
$\LDS$ is any infinite set equipped with the
trivial operation $\LD\xx\yy=\yy$, then $\GGASL$
coincides with~$\GGAS$, \ie, with~$\Sys$. On
the other hand, we can expect that non-trivial
LD-systems give rise to larger geometry
groups, and we can in particular raise:

\begin{ques} \label{Q:Main}
Does there exist a left cancellative
LD-system~$\LDS$ satisfying $\GGASL = \Bs$? 
\end{ques}

A positive answer would correspond to
what can be called a geometric realization
of~$\Bs$, \ie, a realization of~$\Bs$ as the
geometry group of associativity and twisted
semi-commutativity.  

\subsection{$\Bs$-twisted semi-commutativity}

In order to answer Question~\ref{Q:Main} in the
positive, we have to exhibit a convenient
LD-system. Several solutions are possible, but
the quickest and maybe most interesting one
involves a self-distributive structure
on~$\Bs$ itself.
 
\begin{defi}
For~$\bx, \by$ in~$\Bs$, we set
\begin{gather}
\LD\bx\by = \bx \opp \dd\by \opp \sss1
\opp \dd\bx\inv,\\
\bx \OP \by = \bx \opp \dd\by \opp \aa1.
\end{gather}
\end{defi}

\begin{prop}
The set~$\Bs$ equipped with the bracket
operation is a left cancellative LD-system.
Moreover, the following mixed relations are
satisfied
\begin{equation} \label{E:Enhanced}
\LD\bx{\LD\by\bz} = \LD{(\bx \OP \by)}\bz,
\qquad
\LD\bx{\by \OP \bz} = \LD\bx\by \OP \LD\bx\bz,
\end{equation}
where $\dd$ denotes the endomorphism of~$\Bs$ that
maps~$\sss i$ to~$\sss{i+1}$ and $\aa i$
to~$\aa{i+1}$ for every~$i$.
\end{prop}

The self-distributivity of the bracket operation
and the relations~\eqref{E:Enhanced} follow from
the relations of~$\Rass$ using easy
verifications; proving that the bracket operation
is left cancellativity requires to know that the
endomorphism~$\dd$ is injective, which in turn
uses the decomposition of~$\Bs$ as a group of
fractions. As the arguments appear in~\cite{Dhe},
we shall not repeat them here.

\begin{prop} \label{P:PreB}
The group~$\GG{\Ass, \Sem^{\Bs}}$ is
(isomorphic to)~$\Bs$, \ie, $\Bs$ is the
geometry group of associativity and
$\Bs$-twisted semi-commutativity. 
\end{prop}

Proving Proposition~\ref{P:PreB} amounts to
proving that the relations~$\Rass$ make a
presentation of the group~$\GG{\Ass, \Sem^{\Bs}}$ in terms of
the generators~$\aa i$ and~$\sss i$, \ie,
equivalently, that the canonical surjective
hommorphism of~$\Bs$ onto~$\GG{\Ass,
\Sem^{\Bs}}$ is an isomorphism. We use
Proposition~\ref{P:InjCrit}. To this end, we
associate with every $\Bs$-coloured tree a
distinguished element of~$\Bs$ in such a way
that the (external) action of~$\Bs$ on trees
corresponds to an (internal) multiplication
inside~$\Bs$. We proceed in two steps. 

\begin{defi}
For~$t$ a $\Bs$-coloured tree, we define
$\ee(t)$ to be the $\OP$-evaluation of~$t$,
\ie, to be the element of~$\Bs$ inductively
defined by~$\ee(\et_x) = x$ and $\ee(t_1 \op
t_2) = \ee(t_1) \OP \ee(t_2)$.
Then we put
$$\ev(\tt) = \ee(t_1) \opp \dd \ee(t_2) \opp
\pp \opp \dd^{n-1} \ee(t_n)$$
where $\vine{t_1, \pp, t_n, \et_x}$ is the
decomposition of~$\tt$ along its right branch.
\end{defi}

The element~$\ev(t)$ is also defined by the
inductive rules $\ev(\et_x) = 1$ and $\ev(t_1
\op t_2) = \ee(t_1) \opp \dd\ev(t_2)$. Observe
that the definitions of~$\ee(t)$ and~$\ev(t)$
are parallel to those of~$\HHH t$ and~$\HH t$
in Section~\ref{S:PreF}. The key point is the
following computation:

\begin{lemm}
Assume that $\tt$ is an $\Bs$-coloured tree. 
Then we have
\begin{equation}\label{E:MainTreeEval}
\ev(\tt \act \aa i) = \ev(\tt) \opp \aa i, 
\qquad
\ev(\tt \act \sss i) = \ev(\tt) \opp \sss i, 
\end{equation}
whenever the involved trees are defined.
\end{lemm}

\begin{proof}
First, we observe that, for all $\Bs$-coloured
trees~$\tti, \ttii$, we have
\begin{equation} \label{E:Evaluation}
\ee(\LD\tti\ttii) =
\LD{\ee(\tti)}{\ee(\ttii)},
\end{equation}
as follows from an induction on the sizes
of~$\tti$ and~$\ttii$, using the relations
of~\eqref{E:Enhanced}

Now, for~$\tt = t_1 \op \pp \op t_n \op \et_x$,
let~$D(t)$ denote the sequence $(t_1, \pp,
t_n)$, and let $\eee(t)$ be the sequence
$(\ee(t_1), \pp, \ee(t_n))$. By definition, we
have
\begin{gather*}
D(t \act \aa i) = (t_1, \pp, t_{i-1}, t_i \op t_{i+1}, t_{i+2},
\pp, t_n),\\
D(t \act \sss i) = (t_1, \pp, t_{i-1}, 
\LD{t_i}{t_{i+1}}, t_i, t_{i+2}, \pp, t_n).
\end{gather*}
Hence, assuming $\eee(t) = (x_1, \pp, x_n)$, 
and using~\eqref{E:Evaluation} for
the second relation, we obtain 
\begin{gather*}
\eee(t \act \aa i) = (x_1, \pp, x_{i-1}, x_i \OP x_{i+1},
x_{i+2}, \pp, x_n),\\
\eee(t \act \sss i) = (x_1, \pp, x_{i-1},
\LD{x_i}{x_{i+1}},
x_i, x_{i+2}, \pp, x_n)
\end{gather*}
whenever the involved terms are defined. Using
the explicit definition of~$\ev(t \act \aa i)$
and~$\ev(t \act \sss i)$ from $\eee(t \act \aa
i)$ and $\eee(t \act \sss i)$ then easily
gives~\eqref{E:MainTreeEval} using the relations
of~$\Ras$.
\end{proof}

We are now able to conclude.

\begin{proof}[Proof of Proposition~\ref{P:PreB}]
We are in position for applying
Proposition~\ref{P:InjCrit}. Indeed, we have a
surjective homomorphism $\Bs \to \GGASB$
together with a map $\ev: \Tree_{\Bs} \to \Bs$
satisfying~\eqref{E:MainTreeEval}, which are
the relations~\eqref{E:InjCrit} corresponding
to the generating subset~$\Ga \cup \Gss$
of~$\Bs$.
\end{proof}

\subsection{The group of general twisted
semi-commutativity}

To conclude with a simple statement,
let~$\SSps\a$ denote the union of all
operators~$\SSpL\a$ (considered as sets of
pairs) for all possible sets~$\Tree_\LDS$
associated with a left cancellative LD-system,
and define~$\GMASs$ to be the monoid generated
by all operators~$\AAp\a$, $\SSps\a$ and their
inverses. By construction, each specific
monoid~$\GMASL$ is a quotient of~$\GMASs$. Then
near-equality is still a congruence on~$\GMASs$,
and the corresponding group~$\GGASs$ naturally
appears as the geometry group of associativity
and (general) twisted semi-commutativity. We can
state:

\begin{prop} \label{P:GeomGrAbstractSC}
The group~$\GGASs$ is (isomorphic) to~$\Bs$,
\ie, $\Bs$ is the geometry group of
associativity and twisted semi-commuta\-tivity.
\end{prop}

\begin{proof}
By construction, the group~$\GG{\Ass, \Sem^{\Bs}}$ is
a quotient of the general group~$\GGASs$. By
Proposition~\ref{P:Cond}$(iii)$, the
group~$\GGASs$ is a quotient of~$\Bs$.
Now, Proposition~\ref{P:PreB} shows that the
canonical mapping of~$\Bs$ to~$\GG{\Ass, \Sem^{\Bs}}$ is an
isomorphism, so the two surjective
homomorphisms of which the latter is the product
must be isomorphisms as well.
\end{proof}

We mentioned above that group conjugacy
provides examples of left cancellative
LD-systems. Therefore, we obtain for each
particular group~$G$ a notion of
$\conj(G)$-twisted (semi)-commutati\-vity, with an
associated inverse monoid~$\GMASG$ and the
associated group $\GGASG$. The
latter group depends on the group~$G$: if $G$
is abelian, conjugacy is trivial on~$G$, and the
geometry group~$\GGASG$ is therefore~$\Sys$, as
was proved in Section~\ref{S:PreU}. On the other
hand, if $G$ is a non-abelian free group,
conjugacy is not trivial, and we raise the
question of recognizing the corresponding geometry
group.

\begin{prop} \label{P:GeomGroupConj}
If $\FG$ be a free group of rank at least~$2$, the
group~$\GGASF$ is (isomorphic to)~$\Bs$, \ie,
$\Bs$ is the geometry group of associativity and
$\conj(\FG)$-twisted semi-commutativity.
\end{prop}

\begin{proof}[Proof (sketch)]
Without loss of generality, we can assume that
$\FG$ is a free group based on a family of
generators~$x_\a$ indexed by binary addresses,
\ie, finite sequences of~$0$'s and~$1$'s. The
problem is to show that, if $\ww$ is a word
in~$\WW{\Ga,\Gss}$, then the image~$\cl\ww$
of~$\ww$ in~$\Bs$ can be recovered from the
operator of~$\GMASF$ associated with~$\ww$. 

Now, it is shown in~\cite{Dhe} that Artin's
representation of the braid group~$\Bi$ extends
to~$\Bs$: there exists an injective
morphism~$\psi$ of~$\Bs$ into~$\Aut(\FG)$. 
Hence, it suffices to prove that $\psi(\cl\ww)$
is determined by the operator of~$\GMASF$
associated with~$\ww$. We claim that there exists a
$\FG$-coloured tree~$\tt$ such that $\tt \act
\ww$ exists and
$\psi(\cl\ww)$ can be recovered from the
pair~$(\tt, \tt \act \ww)$, hence {\it a
fortiori} from the operator of~$\GMASF$
associated with~$\ww$.

Let us say that a $\FG$-coloured tree~$t$ is {\it
natural} if the labels of~$t$ of each leaf with
address~$\a\0\1^k$ is $x_{\a\0\1^{k-1}}\inv
\pp x_{\a\0\1}\inv x_{\a\0}\inv x_\a$ and the one
of the leaf with address~$\1^k$
is~$x_{\1^{k-1}}\inv \pp
x_{\1}\inv x_{\eadd}\inv$. Proposition~5.4
of~\cite{Dhe} shows (with different notation)
that, if $\tt$ is a natural $\FG$-coloured tree,
and $\ww$ is a word in~$\WW{\Ga, \Gss}$ such that
$\tt \act \ww$ is defined, then, for each
address~$\g$ such that $\g\0$ is the address of a
leaf in~$t \act w$, the image of~$x_\g$
under~$\psi(\cl w)$ is the label at~$\g\0$
in~$\tt \act \ww$: the property can be checked
for~$\aa i$ and~$\sss i$ directly, and, then, one
uses an induction on the length of~$\ww$. As we
can choose~$\tt$ as large as we wish, this shows
that $\psi(\cl\ww)$, hence~$\cl\ww$, is determined
by the action of~$\ww$ on
$\FG$-coloured trees.
\end{proof}

Proposition~\ref{P:GeomGroupConj} gives an
alternative proof of
Proposition~\ref{P:GeomGrAbstractSC}.

As a final remark, let us observe that the
above treatment of twisted semi-commutativity
and its connection with the group~$\Bs$ can
be repeated for twisted commutativity and its
connection with the group obtained by removing
the torsion relations $\cc i^2 = 1$ in the
presentation of~$V$ described in
Section~\ref{S:PreV}. The latter group is
(isomorphic to) the group denoted~$BV$
in~\cite{Bri1, Bri2}, and it also identifies with
the subgroup of~$\Bs$ generated by the elements
$\aa1\inv \pp \aa i\inv \aa{i+1}\inv \aa i \aa i
\pp \aa 1$ and $\aa1\inv \pp \aa i\inv \sss i \aa i
\pp \aa 1$ corresponding to the elements that,
under the action by associativity and twisted
semi-commutativity, act trivially outside the
$0$-subtree.

\end{document}